\documentclass{amsart}
\usepackage{amssymb}
\usepackage{amsmath}
\usepackage{amsfonts}

\newtheorem{theorem}{Theorem}[section]
\theoremstyle{plain}

\newtheorem{corollary}[theorem]{Corollary}

\newtheorem{definition}[theorem]{Definition}

\newtheorem{lemma}[theorem]{Lemma}

\newtheorem{proposition}[theorem]{Proposition}
\newtheorem{remark}[theorem]{Remark}

\numberwithin{equation}{section}
\input{tcilatex}

\begin{document}
\title[On Schatten-class properties of pseudo-differential operators]{On
Schatten-von Neumann class properties of pseudo-differential operators. The
Cordes-Kato method.}
\author{Gruia Arsu}
\address{Institute of Mathematics of The Romanian Academy\\
P.O. Box 1-174\\
RO-70700\\
Bucharest \\
Romania}
\email{agmil@home.ro}
\maketitle

\begin{abstract}
We investigate the Schatten-class properties of pseudo-differential
operators with the (revisted) method of Cordes and Kato. As symbol classes
we use classes similar to those of Cordes in which the \ $L^{\infty }$%
-conditions are replaced by\ $L^{p}$-conditions, $1\leq p<\infty $.
\end{abstract}

\section{Introduction}

In two classical papers \cite{Cordes} and \cite{Kato}, H.O. Cordes and T.
Kato develop an elegant method to deal with pseudo-differential operators.
In \cite{Cordes}, H.O. Cordes shows, among others, that if a symbol $a\left(
x,\xi \right) $\ defined on $%
%TCIMACRO{\U{211d} }%
%BeginExpansion
\mathbb{R}
%EndExpansion
^{n}\times 
%TCIMACRO{\U{211d} }%
%BeginExpansion
\mathbb{R}
%EndExpansion
^{n}$ has bounded derivatives $D_{x}^{\alpha }D_{\xi }^{\beta }a$\ for $%
\left\vert \alpha \right\vert ,\left\vert \beta \right\vert \leq $\ $\left[
n/2\right] +1$, then the associated pseudo-differential operator\ $A=a\left(
x,D\right) $ is $L^{2}$-bounded. This result, known as Calder\'{o}%
n-Vaillancourt theorem, appears for the first time in \cite{Calderon},
except that the number of the required derivatives is different. Cordes-Kato
method can be used to obtain Calder\'{o}n-Vaillancourt theorem with a
minimal H\"{o}lder continuity assumption on the symbol of a
pseudo-differential operator. (See \cite{Childs} and \cite{Bourdaud} for
another approach).

The main purpose of the present paper is to describe an extension of the
Cordes-Kato method which can be used to prove trace-class properties of
pseudo-differential operators. Among others, we prove that $A\in \mathcal{B}%
_{p}\left( L^{2}\left( 
%TCIMACRO{\U{211d} }%
%BeginExpansion
\mathbb{R}
%EndExpansion
^{n}\right) \right) $ if $D_{x}^{\alpha }D_{\xi }^{\beta }a\in L^{p}\left( 
%TCIMACRO{\U{211d} }%
%BeginExpansion
\mathbb{R}
%EndExpansion
^{n}\times 
%TCIMACRO{\U{211d} }%
%BeginExpansion
\mathbb{R}
%EndExpansion
^{n}\right) $ for $\left\vert \alpha \right\vert ,\left\vert \beta
\right\vert \leq $\ $\left[ n/2\right] +1$ and $1\leq p<\infty $, where $%
\mathcal{B}_{p}\left( L^{2}\left( 
%TCIMACRO{\U{211d} }%
%BeginExpansion
\mathbb{R}
%EndExpansion
^{n}\right) \right) $ denote the Schatten ideal of compact operators whose
singular values lie in $l^{p}$. It is remarkable that this method can be
used for $\left( X,\tau \right) $-quantization, in particular, for both the
Weyl quantization and Kohn-Nirenberg quantization.

For other approaches to the $L^{2}$-boundedness and Schatten-class
properties of pseudo-differential operators we refer to the works \cite%
{Grochenig}, \cite{Heil}, \cite{Toft1}, \cite{Toft2}, \cite{Toft3}, \cite%
{Toft4}.

\section{Weyl calculus in the Schr\"{o}dinger representation}

We find it convenient to present the results in the Schr\"{o}dinger
representation formalism (or Weyl systems formalism).

A symplectic space is a real finite dimensional vector space $\mathfrak{S}$
equipped with a real antisymmetric nondegenerate bilinear form $\sigma $.

We define the \textit{Fourier measure} $\func{d}\xi $ as the unique
translation invariant, Borel regular measure on $\mathfrak{S}$ such that the
(symplectic) Fourier transform%
\begin{equation*}
\left( \mathcal{F}_{\mathfrak{S}}a\right) \left( \xi \right) \equiv \widehat{%
a}\left( \xi \right) =\int_{\mathfrak{S}}e^{-i\sigma \left( \xi ,\eta
\right) }a\left( \eta \right) \func{d}\eta
\end{equation*}%
is involutive (i.e.$\mathcal{F}_{\mathfrak{S}}^{2}=1$) and unitary on $L^{2}(%
\mathfrak{S})$. We use the same notation $\mathcal{F}_{\mathfrak{S}}\mathcal{%
\ }$for the extension to $\mathcal{S}^{\ast }(\mathfrak{S})$ of this Fourier
transform.

We recall some facts in connection with the theory of canonical commutation
relations (see \cite{Georgescu}, \cite{Kastler}). We denote by $\mathcal{M}(%
\mathfrak{S})$ the space of integrable measures on $\mathfrak{S}$ equipped
with \ a unital $\ast $-algebra structure defined by the twisted convolution%
\begin{equation*}
(\mu \times \nu )(\xi )=\int_{\mathfrak{S}}e^{\frac{i}{2}\sigma (\xi ,\eta
)}\mu (\xi -\eta )\nu (\eta )\func{d}\eta
\end{equation*}%
as product, $\mu ^{\ast }(\xi )=\overline{\mu (-\xi )}$ as involution and $%
\delta $ the Dirac measure at $0$ as unit. These definition must be
interpreted in the sense of distributions, i.e. for $f\in \mathcal{C}_{0}(%
\mathfrak{S})$:%
\begin{eqnarray*}
\int_{\mathfrak{S}}f(\xi )\func{d}(\mu \times \upsilon )(\xi ) &=&\iint e^{%
\frac{i}{2}\sigma (\xi ,\eta )}f(\xi +\eta )\limfunc{d}\mu \left( \xi
\right) \func{d}\upsilon (\eta ), \\
\int_{\mathfrak{S}}f(\xi )\func{d}\mu ^{\ast }(\xi ) &=&\overline{\int_{%
\mathfrak{S}}\overline{f(-\xi )}\limfunc{d}\mu \left( \xi \right) }.
\end{eqnarray*}

Let $\widehat{\mathcal{M}}\left( \mathfrak{S}\right) $ be the space of
Fourier transforms of measures in $\mathcal{M}(\mathfrak{S})$, provided with
the product $a\circ b=\mathcal{F}_{\mathfrak{S}}\left( \widehat{a}\times 
\widehat{b}\right) $ and with the usual conjugation $a^{\ast }=\overline{a}$(%
$=\mathcal{F}_{\mathfrak{S}}\left( \widehat{a}^{\ast }\right) )$. $\widehat{%
\mathcal{M}}\left( \mathfrak{S}\right) $ becomes a unital $\ast $-algebra,
with the function identically equal to $1$ ($\mathcal{F}_{\mathfrak{S}%
}\left( \delta \right) =1$) as unit. Observe that $\widehat{\mathcal{M}}%
\left( \mathfrak{S}\right) $ consists of bounded continuous functions on $%
\mathfrak{S}$. The product $a\circ b$ is called composition product.

By representation of a symplectic space $\mathfrak{S}$ on a Hilbert space $%
\mathcal{H}$ we understand a strongly continuous map $\mathcal{W}$ from $%
\mathfrak{S}$ to the set of unitary operators on $\mathcal{H}$ satisfying:%
\begin{equation}
\mathcal{W}\left( \xi \right) \mathcal{W}\left( \eta \right) =e^{\frac{i}{2}%
\sigma \left( \xi ,\eta \right) }\mathcal{W}\left( \xi +\eta \right) \quad 
\text{\textit{for all }}\xi ,\eta \in \mathfrak{S}.  \label{SR1}
\end{equation}%
This implies $\mathcal{W}\left( 0\right) =1$, $\mathcal{W}\left( \xi \right)
^{\ast }=\mathcal{W}\left( -\xi \right) $ and%
\begin{equation}
\mathcal{W}\left( \xi \right) \mathcal{W}\left( \eta \right) =e^{i\sigma
\left( \xi ,\eta \right) }\mathcal{W}\left( \eta \right) \mathcal{W}\left(
\xi \right) \quad \text{\textit{for all }}\xi ,\eta \in \mathfrak{S}.
\label{SR2}
\end{equation}%
The couple $\left( \mathcal{H},\mathcal{W}\right) $ is also called Weyl
system associated to the symplectic space $\mathfrak{S}$.

For integrable Borel measures $\mu $ on $\mathfrak{S}$, we denote 
\begin{equation}
\mathcal{W}\left( \mu \right) =\int_{\mathfrak{S}}\mathcal{W}(\xi )\func{d}%
\mu (\xi ).  \label{SR3}
\end{equation}%
One can easily check then that 
\begin{equation*}
\mathcal{W}\left( \mu \times \upsilon \right) =\mathcal{W}\left( \mu \right) 
\mathcal{W}\left( \upsilon \right) \quad \text{\textit{and}\quad }\mathcal{W}%
\left( \mu ^{\ast }\right) =\mathcal{W}\left( \mu \right) ^{\ast },
\end{equation*}%
i.e. $\mathcal{W}:\mathcal{M}(\mathfrak{S})\rightarrow \mathcal{B(H)}$ is a $%
\ast $-representation of $\mathcal{M}(\mathfrak{S})$ on the Hilbert space $%
\mathcal{H}$. Note that for any $\xi \in \mathfrak{S}$, $\mathcal{W}\left(
\delta _{\xi }\right) =\mathcal{W}(\xi )$, where $\delta _{\xi }$ is the
Dirac measure at $\xi $.

There is a bijective corespondence betweeen faithful (irreducible) $\ast $%
-representations $\limfunc{Op}:\widehat{\mathcal{M}}\left( \mathfrak{S}%
\right) \rightarrow \mathcal{B(H)}$ and strongly continuous representations $%
\mathcal{W}:\mathfrak{S}\rightarrow \mathcal{U}\left( \mathcal{H}\right) $
(which act irreducibly in $\mathcal{H}$). The corespondence is specified by%
\begin{equation}
\limfunc{Op}\left( a\right) =a\left( R\right) =\mathcal{W}\left( \widehat{a}%
\right) =\int_{\mathfrak{S}}\mathcal{W}\left( \xi \right) \func{d}\widehat{a}%
\left( \xi \right)  \label{SR4}
\end{equation}%
for $a\in \widehat{\mathcal{M}}\left( \mathfrak{S}\right) $. The application 
$\limfunc{Op}$ is called Weyl calculus, $\limfunc{Op}\left( a\right) $ being
the operator associated to the symbol $a$.

Let $\mathcal{S}$ be the dense linear subspace of $\mathcal{H}$ consisting
of the $\mathcal{C}^{\infty }$ vectors of the representation $\mathcal{W}$%
\begin{equation}
\mathcal{S}=\mathcal{S}(\mathcal{H},\mathcal{W})=\left\{ \varphi \in 
\mathcal{H}:\mathfrak{S}\ni \xi \rightarrow \mathcal{W}(\xi )\varphi \in 
\mathcal{H}\text{\textit{\ is a }}\mathcal{C}^{\infty }\text{\textit{\ map}}%
\right\} .  \label{SR5}
\end{equation}%
For each $\xi \in \mathfrak{S}$ the family $\left\{ \mathcal{W}(t\xi
)\right\} _{t\in 
%TCIMACRO{\U{211d} }%
%BeginExpansion
\mathbb{R}
%EndExpansion
}$ is a strongly continuous unitary representation of $%
%TCIMACRO{\U{211d} }%
%BeginExpansion
\mathbb{R}
%EndExpansion
$ in $\mathcal{H}$ which leaves $\mathcal{S}$ invariant. We denote $\sigma
(\xi ,R)$ the infinitesimal generator of this group, so that 
\begin{equation}
\mathcal{W}(t\xi )=e^{it\sigma (\xi ,R)},\quad \text{\textit{for all} }t\in 
%TCIMACRO{\U{211d} }%
%BeginExpansion
\mathbb{R}
%EndExpansion
.  \label{SR6}
\end{equation}%
Clearly $\mathcal{S}\subset D\left( \sigma (\xi ,R)\right) $, $\mathcal{S}$
is stable under $\sigma (\xi ,R)$, and $\sigma (\xi ,R)\mid \mathcal{S}$ is
essentially self-adjoint by Nelson's lemma (Theorem VIII.11 in \cite{Simon}%
). From (\ref{SR1}) and (\ref{SR2}) we get%
\begin{gather}
\sigma (\xi +\eta ,R)=\sigma (\xi ,R)+\sigma (\eta ,R),  \notag \\
\mathcal{W}(\xi )\sigma (\eta ,R)\mathcal{W}(-\xi )=\sigma (\eta ,R)+\sigma
(\xi ,\eta ),\quad i\left[ \sigma (\xi ,R),\sigma (\eta ,R)\right] =\sigma
(\xi ,\eta )  \label{SR7}
\end{gather}%
on $\mathcal{S}$ for all $\xi ,\eta \in \mathfrak{S}$.

Note that $R$ in $\sigma (\xi ,R)$ or $a\left( R\right) $ is the pair $%
\left( Q,P\right) $, where $Q$ is position and $P$ is momentum observable in
physics and the pair $\left( x,D\right) $ for users of pseudo-differential
calculus as presented , for example, in ch. XVIII of \cite{Hormander2} or in 
\cite{Folland}.

The space $\mathcal{S}$ can be described in terms of the subspaces $D\left(
\sigma (\xi ,R)\right) $%
\begin{eqnarray*}
\mathcal{S} &=&\bigcap_{k\in 
%TCIMACRO{\U{2115} }%
%BeginExpansion
\mathbb{N}
%EndExpansion
}\bigcap_{\xi _{1},...,\xi _{k}\in \mathfrak{S}}D\left( \sigma (\xi
_{1},R)...\sigma (\xi _{k},R)\right) \\
&=&\bigcap_{k\in 
%TCIMACRO{\U{2115} }%
%BeginExpansion
\mathbb{N}
%EndExpansion
}\bigcap_{\xi _{1},...,\xi _{k}\in \mathcal{B}}D\left( \sigma (\xi
_{1},R)...\sigma (\xi _{k},R)\right) ,
\end{eqnarray*}%
where $\mathcal{B}$ is a (symplectic) basis. The topology in $\mathcal{S}$
defined by the family of seminorms $\left\{ \left\Vert \cdot \right\Vert
_{k,\xi _{1},...,\xi _{k}}\right\} _{k\in 
%TCIMACRO{\U{2115} }%
%BeginExpansion
\mathbb{N}
%EndExpansion
,\xi _{1},...,\xi _{k}\in \mathfrak{S}}$ 
\begin{equation*}
\left\Vert \varphi \right\Vert _{k,\xi _{1},...,\xi _{k}}=\left\Vert \sigma
(\xi _{1},R)...\sigma (\xi _{k},R)\varphi \right\Vert _{\mathcal{H}},\quad
\varphi \in \mathcal{S}
\end{equation*}%
makes $\mathcal{S}$ a Fr\'{e}chet space. We denote by $\mathcal{S}^{\ast }$
the space of all continuous, antilinear (semilinear) mappings $\mathcal{S}%
\rightarrow \mathbb{%
%TCIMACRO{\U{2102} }%
%BeginExpansion
\mathbb{C}
%EndExpansion
}$ equipped with the weak topology $\sigma (\mathcal{S}^{\ast },\mathcal{S})$%
. Since $\mathcal{S}\hookrightarrow \mathcal{H}$ continuously and densely,
and since $\mathcal{H}$ is always identified with its adjoint $\mathcal{H}%
^{\ast }$, we obtain a scale of dense inclusions%
\begin{equation*}
\mathcal{S}\hookrightarrow \mathcal{H}\hookrightarrow \mathcal{S}^{\ast }
\end{equation*}%
such that, if $\left\langle \cdot ,\cdot \right\rangle :\mathcal{S}\times 
\mathcal{S}^{\ast }\rightarrow \mathbb{%
%TCIMACRO{\U{2102} }%
%BeginExpansion
\mathbb{C}
%EndExpansion
}$ is the antiduality between $\mathcal{S}$ and $\mathcal{S}^{\ast }$
(antilinear in the first and linear in the second argument), then for $%
\varphi \in \mathcal{S}$ and $u\in \mathcal{H}$, if $u$ is considered as an
element of $\mathcal{S}^{\ast }$, the number $\left\langle \varphi
,u\right\rangle $ is just the scalar product in $\mathcal{H}$. For this
reason we do not distinguish between the the scalar product in $\mathcal{H}$
and the antiduality between $\mathcal{S}$ and $\mathcal{S}^{\ast }$.

For $\xi \in \mathfrak{S}$, $\mathcal{W}(\xi )$ and $\mathcal{W}(-\xi )$ are
topological isomorphisms of $\mathcal{S}$ such that $\mathcal{W}(-\xi )%
\mathcal{W}(\xi )=\mathcal{W}(\xi )\mathcal{W}(-\xi )=1_{\mathcal{S}}$.
Hence we can extend $\mathcal{W}(\xi )$ to $\mathcal{S}^{\ast }$ as the
adjoint of the mapping of $\mathcal{W}(-\xi )|\mathcal{S}$. Then $\mathcal{W}%
(\xi )$ and $\mathcal{W}(-\xi )$ are isomorphisms of $\mathcal{S}^{\ast }$
onto itself such that $\mathcal{W}(-\xi )=\mathcal{W}(\xi )^{-1}$.

The first formula of (\ref{SR7}) and an induction argument give 
\begin{equation*}
\mathcal{W}(-\xi )\sigma (\eta _{1},R)...\sigma (\eta _{k},R)\mathcal{W}(\xi
)=\left( \sigma (\eta _{1},R)+\sigma (\eta _{1},\xi )\right) ...\left(
\sigma (\eta _{k},R)+\sigma (\eta _{k},\xi )\right)
\end{equation*}%
for all $\xi ,\eta _{1},...\eta _{k}\in \mathfrak{S}$. Let $\mathcal{B}%
=\left\{ \varepsilon _{1},...,\varepsilon _{n},\varepsilon
^{1},...,\varepsilon ^{n}\right\} $ be a symplectic basis, i.e. we have for $%
j,k=1,...,n$%
\begin{equation*}
\sigma (\varepsilon _{j},\varepsilon _{k})=\sigma (\varepsilon
^{j},\varepsilon ^{k})=\sigma (\varepsilon ^{j},\varepsilon _{k})-\delta
_{jk}=0
\end{equation*}%
where $\delta _{jk}$ is the Kroneker delta, equal to $1$ when $j=k$ and $0$
when $j\neq k$. By making a suitable choise of the vectors $\eta
_{1},...\eta _{k}$ in $\mathcal{B}$, we obtain 
\begin{equation*}
e^{i\left\langle x,P\right\rangle }Q^{\alpha }e^{-i\left\langle
x,P\right\rangle }=\left( Q+x\right) ^{\alpha },\quad e^{-i\left\langle
Q,p\right\rangle }P_{\beta }e^{i\left\langle Q,p\right\rangle }=\left(
P+p\right) _{\beta }.
\end{equation*}%
Here we used the following notations%
\begin{gather*}
\mathfrak{S}=X+X^{\ast },\quad X=%
%TCIMACRO{\U{211d} }%
%BeginExpansion
\mathbb{R}
%EndExpansion
\varepsilon _{1}+...+%
%TCIMACRO{\U{211d} }%
%BeginExpansion
\mathbb{R}
%EndExpansion
\varepsilon _{n},\quad X^{\ast }=%
%TCIMACRO{\U{211d} }%
%BeginExpansion
\mathbb{R}
%EndExpansion
\varepsilon ^{1}+...+%
%TCIMACRO{\U{211d} }%
%BeginExpansion
\mathbb{R}
%EndExpansion
\varepsilon ^{n}, \\
\xi =x+p,\quad x=x^{1}\varepsilon _{1}+...+x^{n}\varepsilon _{n},\quad
p=p_{1}\varepsilon ^{1}+...+p_{n}\varepsilon ^{n}, \\
x^{1}=\sigma (\varepsilon ^{1},\xi ),...,x^{n}=\sigma (\varepsilon ^{n},\xi
),\quad p_{1}=\sigma (-\varepsilon _{1},\xi ),...,p_{n}=\sigma (-\varepsilon
_{n},\xi ), \\
\left\langle x,p\right\rangle =\sigma (p,x)=x^{1}p_{1}+...+x^{n}p_{n}, \\
x^{\alpha }=\left( x^{1}\right) ^{\alpha _{1}}...\left( x^{n}\right)
^{\alpha _{n}},\quad p_{\beta }=\left( p_{1}\right) ^{\beta _{1}}....\left(
p_{n}\right) ^{\beta _{n}},\quad \alpha ,\beta \in 
%TCIMACRO{\U{2115} }%
%BeginExpansion
\mathbb{N}
%EndExpansion
^{n}
\end{gather*}%
and 
\begin{gather*}
\left\langle Q,p\right\rangle =\sigma (p,R),\quad p\in X^{\ast },\quad
\left\langle x,P\right\rangle =\sigma (-x,R),\quad x\in X, \\
\mathcal{W}(p)=e^{i\left\langle Q,p\right\rangle },\quad p\in X^{\ast
},\quad \mathcal{W}(x)=e^{-i\left\langle x,P\right\rangle },\quad x\in X, \\
Q^{j}=\left\langle Q,\varepsilon ^{j}\right\rangle =\sigma (\varepsilon
^{j},R),\quad P_{j}=\left\langle \varepsilon _{j},P\right\rangle =\sigma
(-\varepsilon _{j},R),\quad j=1,...,n, \\
Q^{\alpha }=\left( Q^{1}\right) ^{\alpha _{1}}...\left( Q^{n}\right)
^{\alpha _{n}},\quad P_{\beta }=\left( P_{1}\right) ^{\beta _{1}}....\left(
P_{n}\right) ^{\beta _{n}},\quad \alpha ,\beta \in 
%TCIMACRO{\U{2115} }%
%BeginExpansion
\mathbb{N}
%EndExpansion
^{n}.
\end{gather*}%
These formulas together with the binomial formula give 
\begin{eqnarray*}
x^{\alpha } &=&\left( Q+x-Q\right) ^{\alpha }=\sum_{0\leq \gamma \leq \alpha
}\left( -1\right) ^{\left\vert \gamma \right\vert }\left( 
\begin{array}{c}
\alpha \\ 
\gamma%
\end{array}%
\right) \left( Q+x\right) ^{\alpha -\gamma }Q^{\gamma } \\
&=&\sum_{0\leq \gamma \leq \alpha }\left( -1\right) ^{\left\vert \gamma
\right\vert }\left( 
\begin{array}{c}
\alpha \\ 
\gamma%
\end{array}%
\right) e^{i\left\langle x,P\right\rangle }Q^{\alpha -\gamma
}e^{-i\left\langle x,P\right\rangle }Q^{\gamma }, \\
p_{\beta } &=&\left( P-\left( P-p\right) \right) _{\beta }=\sum_{0\leq
\gamma \leq \beta }\left( -1\right) ^{\left\vert \gamma \right\vert }\left( 
\begin{array}{c}
\beta \\ 
\gamma%
\end{array}%
\right) P_{\beta -\gamma }\left( P-p\right) _{\gamma } \\
&=&\sum_{0\leq \gamma \leq \beta }\left( -1\right) ^{\left\vert \gamma
\right\vert }\left( 
\begin{array}{c}
\beta \\ 
\gamma%
\end{array}%
\right) P_{\beta -\gamma }e^{i\left\langle Q,p\right\rangle }P_{\gamma
}e^{-i\left\langle Q,p\right\rangle }.
\end{eqnarray*}%
Hence on $\mathcal{S}$ we have 
\begin{multline*}
x^{\alpha }p_{\beta }e^{i\left\langle Q,p\right\rangle }e^{-i\left\langle
x,P\right\rangle }=p_{\beta }e^{i\left\langle Q,p\right\rangle
}e^{-i\left\langle x,P\right\rangle }x^{\alpha } \\
=\sum_{0\leq \gamma \leq \beta ,0\leq \delta \leq \alpha }\left( -1\right)
^{\left\vert \gamma \right\vert +\left\vert \delta \right\vert }\left( 
\begin{array}{c}
\beta \\ 
\gamma%
\end{array}%
\right) \left( 
\begin{array}{c}
\alpha \\ 
\delta%
\end{array}%
\right) P_{\beta -\gamma }Q^{\alpha -\delta }e^{i\left\langle
Q,p\right\rangle }e^{-i\left\langle x,P\right\rangle }Q^{\delta }P_{\gamma }.
\end{multline*}

We can now prove prove the following result.

\begin{lemma}
\label{SR8}If $\varphi ,\psi \in \mathcal{S}$, then the map $\mathfrak{S}\ni
\xi \rightarrow \left\langle \psi ,\mathcal{W}(\xi )\varphi \right\rangle _{%
\mathcal{S},\mathcal{S}^{\ast }}\in 
%TCIMACRO{\U{2102} }%
%BeginExpansion
\mathbb{C}
%EndExpansion
$ belongs to $\mathcal{S}(\mathfrak{S})$. Moreover, for each continuous
seminorm $p$ on $\mathcal{S}(\mathfrak{S})$ there are continuous seminorms $%
q $ and $q^{\prime }$ on $\mathcal{S}$ such that 
\begin{equation*}
p\left( \left\langle \psi ,\mathcal{W}(\cdot )\varphi \right\rangle _{%
\mathcal{S},\mathcal{S}^{\ast }}\right) \leq q(\psi )q^{\prime }(\varphi ).
\end{equation*}
\end{lemma}

\begin{proof}
With the above notations for $\xi =x+p$ we have 
\begin{equation*}
\mathcal{W}(\xi )=\mathcal{W}(x+p)=e^{-\frac{i}{2}\sigma (p,x)}\mathcal{W}(p)%
\mathcal{W}(x)=e^{-\frac{i}{2}\left\langle x,p\right\rangle
}e^{i\left\langle Q,p\right\rangle }e^{-i\left\langle x,P\right\rangle }
\end{equation*}%
so it suffices to study the behaviour of the map $\left( x,p\right)
\rightarrow e^{i\left\langle Q,p\right\rangle }e^{-i\left\langle
x,P\right\rangle }$ at infinity. We may use the decomposition of $x^{\alpha
}p_{\beta }e^{i\left\langle Q,p\right\rangle }e^{-i\left\langle
x,P\right\rangle }$ to obtain that 
\begin{multline*}
\left\vert x^{\alpha }p_{\beta }\left\langle \psi ,e^{i\left\langle
Q,p\right\rangle }e^{-i\left\langle x,P\right\rangle }\varphi \right\rangle
\right\vert \\
\leq \sum_{0\leq \gamma \leq \beta ,0\leq \delta \leq \alpha }\left( 
\begin{array}{c}
\beta \\ 
\gamma%
\end{array}%
\right) \left( 
\begin{array}{c}
\alpha \\ 
\delta%
\end{array}%
\right) \left\vert \left\langle Q^{\alpha -\delta }P_{\beta -\gamma }\psi
,e^{i\left\langle Q,p\right\rangle }e^{-i\left\langle x,P\right\rangle
}Q^{\delta }P_{\gamma }\varphi \right\rangle \right\vert \\
\leq \sum_{0\leq \gamma \leq \beta ,0\leq \delta \leq \alpha }\left( 
\begin{array}{c}
\beta \\ 
\gamma%
\end{array}%
\right) \left( 
\begin{array}{c}
\alpha \\ 
\delta%
\end{array}%
\right) \left\Vert Q^{\alpha -\delta }P_{\beta -\gamma }\psi \right\Vert
\left\Vert Q^{\delta }P_{\gamma }\varphi \right\Vert <\infty .
\end{multline*}%
Since 
\begin{eqnarray*}
\partial _{\beta }^{\alpha }\left[ e^{i\left\langle Q,p\right\rangle
}e^{-i\left\langle x,P\right\rangle }\varphi \right] &=&i^{\left\vert \alpha
\right\vert }\left( -i\right) ^{\left\vert \beta \right\vert
}e^{i\left\langle Q,p\right\rangle }Q^{\alpha }e^{-i\left\langle
x,P\right\rangle }P_{\beta }\varphi \\
&=&i^{\left\vert \alpha \right\vert }\left( -i\right) ^{\left\vert \beta
\right\vert }e^{i\left\langle Q,p\right\rangle }e^{-i\left\langle
x,P\right\rangle }e^{i\left\langle x,P\right\rangle }Q^{\alpha
}e^{-i\left\langle x,P\right\rangle }P_{\beta }\varphi \\
&=&i^{\left\vert \alpha \right\vert }\left( -i\right) ^{\left\vert \beta
\right\vert }e^{i\left\langle Q,p\right\rangle }e^{-i\left\langle
x,P\right\rangle }\left( Q+x\right) ^{\alpha }P_{\beta }\varphi \\
&=&i^{\left\vert \alpha \right\vert }\left( -i\right) ^{\left\vert \beta
\right\vert }e^{i\left\langle Q,p\right\rangle }e^{-i\left\langle
x,P\right\rangle }\sum_{0\leq \gamma \leq \alpha }\left( 
\begin{array}{c}
\alpha \\ 
\gamma%
\end{array}%
\right) x^{\alpha -\gamma }Q^{\gamma }P_{\beta }\varphi
\end{eqnarray*}%
and $Q^{\gamma }P_{\beta }\varphi \in \mathcal{S}$, we conclude that all
derivatives of $\left( x,p\right) \rightarrow \left\langle \psi
,e^{i\left\langle Q,p\right\rangle }e^{-i\left\langle x,P\right\rangle
}\varphi \right\rangle $ have similar estimates. The proof is complete.
\end{proof}

Lemma \ref{SR8} makes it natural to extend Weyl calculus from symbols in $%
\widehat{\mathcal{M}}\left( \mathfrak{S}\right) $ to symbols in $\mathcal{S}%
^{\ast }(\mathfrak{S}).$

\begin{corollary}
\label{SR10}Let $\left( \mathcal{H},\mathcal{W}\right) $ be a Weyl system
associated to the symplectic space $\left( \mathfrak{S},\sigma \right) $.
Then for each $\mu ,a\in \mathcal{S}^{\ast }\left( \mathfrak{S}\right) $ we
can define the operators $\mathcal{W}(\mu ),\limfunc{Op}(a):\mathcal{S}%
\rightarrow \mathcal{S}^{\ast }$ by 
\begin{equation*}
\mathcal{W}(\mu )=\int_{\mathfrak{S}}\mathcal{W}(\xi )\mu (\xi )\func{d}\xi
,\quad \limfunc{Op}(a)=\int_{\mathfrak{S}}\mathcal{W}(\xi )\widehat{a}(\xi )%
\func{d}\xi .
\end{equation*}%
The above integrals make sense if they are taken in the weak sense, i.e. for 
$\varphi ,\psi \in \mathcal{S}$%
\begin{eqnarray*}
\left\langle \varphi ,\mathcal{W}(\mu )\psi \right\rangle _{\mathcal{S},%
\mathcal{S}^{\ast }} &=&\left\langle \overline{\left\langle \varphi ,%
\mathcal{W}\left( \cdot \right) \psi \right\rangle }_{\mathcal{S},\mathcal{S}%
^{\ast }},\mu \right\rangle _{\mathcal{S}\left( \mathfrak{S}\right) ,%
\mathcal{S}^{\ast }\left( \mathfrak{S}\right) }, \\
\left\langle \varphi ,\limfunc{Op}(a)\psi \right\rangle _{\mathcal{S},%
\mathcal{S}^{\ast }} &=&\left\langle \overline{\left\langle \varphi ,%
\mathcal{W}\left( \cdot \right) \psi \right\rangle }_{\mathcal{S},\mathcal{S}%
^{\ast }},\widehat{a}\right\rangle _{\mathcal{S}\left( \mathfrak{S}\right) ,%
\mathcal{S}^{\ast }\left( \mathfrak{S}\right) }.
\end{eqnarray*}%
Moreover, from Lemma \ref{SR8} one obtains that 
\begin{equation*}
\left\vert \left\langle \varphi ,\mathcal{W}(\mu )\psi \right\rangle _{%
\mathcal{S},\mathcal{S}^{\ast }}\right\vert +\left\vert \left\langle \varphi
,\limfunc{Op}(a)\psi \right\rangle _{\mathcal{S},\mathcal{S}^{\ast
}}\right\vert \leq p\left( \left\langle \varphi ,\mathcal{W}(\cdot )\psi
\right\rangle _{\mathcal{S},\mathcal{S}^{\ast }}\right) \leq q\left( \varphi
\right) q^{\prime }\left( \psi \right) ,
\end{equation*}%
where $p$ is a continuous seminorm on $\mathcal{S}(\mathfrak{S})$ and $q$
and $q^{\prime }$ are continuous seminorms on $\mathcal{S}$.

If on $\mathcal{S}^{\ast }\left( \mathfrak{S}\right) $ we shall consider the
weak topology $\sigma (\mathcal{S}^{\ast }\left( \mathfrak{S}\right) ,%
\mathcal{S}\left( \mathfrak{S}\right) )$ and on $\mathcal{B}(\mathcal{S},%
\mathcal{S}^{\ast })$ the topology defined by the seminorms $\left\{
p_{\varphi ,\psi }\right\} _{\varphi ,\psi \in \mathcal{S}}$ 
\begin{equation*}
p_{\varphi ,\psi }\left( A\right) =\left\vert \left\langle \varphi ,A\psi
\right\rangle \right\vert ,\quad A\in \mathcal{B}(\mathcal{S},\mathcal{S}%
^{\ast }),
\end{equation*}%
then the mappings%
\begin{eqnarray*}
\mathcal{W} &:&\mathcal{S}^{\ast }\left( \mathfrak{S}\right) \rightarrow 
\mathcal{B}\left( \mathcal{S},\mathcal{S}^{\ast }\right) ,\quad \mu
\rightarrow \mathcal{W}(\mu ), \\
\limfunc{Op} &:&\mathcal{S}^{\ast }\left( \mathfrak{S}\right) \rightarrow 
\mathcal{B}\left( \mathcal{S},\mathcal{S}^{\ast }\right) ,\quad a\rightarrow 
\limfunc{Op}(a)
\end{eqnarray*}%
are well defined linear and continuous.
\end{corollary}

A subspace $E\subset \mathfrak{S}$ is called \textit{isotropic} if $E\subset
E^{\sigma }$ and \textit{involutive} if $E^{\sigma }\subset E$. If both are
valid, i.e. $E=E^{\sigma }$, then $E$ is \textit{lagrangian}. An isotropic
subspace $X\subset \mathfrak{S}$ is lagrangian if and only if $2\dim X=\dim 
\mathfrak{S}$.

Next we shall make the connection with Weyl calculus familiar to users of
pseudo-differential calculus as presented , for example, in ch. XVIII of 
\cite{Hormander2} or in \cite{Folland}.

Let $X$ be an $n$ dimensional vector space over $%
%TCIMACRO{\U{211d} }%
%BeginExpansion
\mathbb{R}
%EndExpansion
$ and $X^{\ast }$ its dual. Denote $x,y,...$ the elements af $X$ and $%
k,p,... $ those of $X^{\ast }$. Let $\left\langle \cdot ,\cdot \right\rangle
:X\times X^{\ast }\rightarrow 
%TCIMACRO{\U{211d} }%
%BeginExpansion
\mathbb{R}
%EndExpansion
$ be the duality form, i.e. a non-degenerate bilinear form. The symplectic
space is defined by $\mathfrak{S}=T^{\ast }(X)=X\times X^{\ast }$ the
symplectic form being $\sigma \left( \left( x,p\right) ,\left( x^{\prime
},p^{\prime }\right) \right) =\left\langle x^{\prime },p\right\rangle
-\left\langle x,p^{\prime }\right\rangle $. Observe that $X$ and $X^{\ast }$
are lagrangian subspaces of $\mathfrak{S}$. Let us mention that there is a
kind of converse to this construction. Let $\left( X,X^{\ast }\right) $ be a
couple of lagrangian subspaces of $\mathfrak{S}$ such that $X\cap X^{\ast
}=0 $ or, equivalently, $X+X^{\ast }=\mathfrak{S}$. If for $x\in X$ and $%
p\in X^{\ast }$ we define $\left\langle x,p\right\rangle =\sigma \left(
p,x\right) $, then we get a non-degenerate bilinear form on $X\times X^{\ast
}$ which allows us to identify $X^{\ast }$ with the dual of $X$. A couple $%
\left( X,X^{\ast }\right) $ of subspaces of $\mathfrak{S}$ with the
preceding properties is called a \textit{holonomic decomposition} of $%
\mathfrak{S}$.

We define the (Fourier) transforms 
\begin{eqnarray*}
\mathcal{F}_{X},\ \overline{\mathcal{F}}_{X} &:&\mathcal{S}^{\ast }\left(
X\right) \rightarrow \mathcal{S}^{\ast }\left( X^{\ast }\right) , \\
\mathcal{F}_{X^{\ast }},\ \overline{\mathcal{F}}_{X^{\ast }} &:&\mathcal{S}%
^{\ast }\left( X^{\ast }\right) \rightarrow \mathcal{S}^{\ast }\left(
X\right) ,
\end{eqnarray*}%
by 
\begin{eqnarray*}
(\mathcal{F}_{X}u)(p) &=&\int_{X}e^{-i\left\langle x,p\right\rangle }u(x)%
\func{d}x,\quad (\overline{\mathcal{F}}_{X}u)(p)=\int_{X}e^{+i\left\langle
x,p\right\rangle }u(x)\func{d}x, \\
(\mathcal{F}_{X^{\ast }}v)(x) &=&\int_{X^{\ast }}e^{-i\left\langle
x,p\right\rangle }v(p)\func{d}p,\quad (\overline{\mathcal{F}}_{X^{\ast
}}v)(x)=\int_{X^{\ast }}e^{+i\left\langle x,p\right\rangle }v(p)\func{d}p.
\end{eqnarray*}%
Here $\func{d}x$ is a Haar measure in $X$ and $\func{d}p$ is the dual one in 
$X^{\ast }$ such that Fourier's inversion formulas $\overline{\mathcal{F}}%
_{X^{\ast }}\circ \mathcal{F}_{X}=1_{\mathcal{S}^{\ast }\left( X\right) }$, $%
\mathcal{F}_{X}\circ \overline{\mathcal{F}}_{X^{\ast }}=1_{\mathcal{S}^{\ast
}\left( X^{\ast }\right) }$ hold. Replacing $\func{d}x$ by $c\func{d}x$ one
must change $\func{d}p$ to $c^{-1}\func{d}p$ so $\func{d}\xi =\func{d}%
x\otimes \func{d}p$ is invariantly defined and it is exactly the Fourier
measure on $\mathfrak{S}$. Then the symplectic Fourier transform is given by 
\begin{equation*}
\left( \mathcal{F}_{\mathfrak{S}}a\right) \left( x,p\right) =\iint_{X\times
X^{\ast }}e^{-i\left[ \left\langle y,p\right\rangle -\left\langle
x,k\right\rangle \right] }a\left( y,k\right) \func{d}y\func{d}k=\left( 
\mathcal{F}_{X}\otimes \overline{\mathcal{F}}_{X^{\ast }}\right) a\left(
p,x\right) ,
\end{equation*}%
so $\mathcal{F}_{\mathfrak{S}}=\mathcal{I}\circ (\mathcal{F}_{X}\otimes 
\overline{\mathcal{F}}_{X^{\ast }})$, where $\mathcal{I}:\mathcal{S}^{\ast
}(X^{\ast }\times X)\rightarrow \mathcal{S}^{\ast }(X\times X^{\ast })$ is
given by $(\mathcal{I}b)(x,p)=b(p,x)$.

To each finite dimensional vector space $X$ over $%
%TCIMACRO{\U{211d} }%
%BeginExpansion
\mathbb{R}
%EndExpansion
$ and to each Haar measure $\func{d}x$ on $X$, one associates a
representation of the symplectic space $\mathfrak{S}=T^{\ast }(X)=X\times
X^{\ast }$, the Schr\"{o}dinger representation, defined as follows: $%
\mathcal{H}\left( X\right) =L^{2}\left( X,\limfunc{d}x\right) $ and for $\xi
=\left( x,p\right) $ and $\psi \in \mathcal{H}\left( X\right) $%
\begin{equation*}
\mathcal{W}\left( \xi \right) \psi \left( \cdot \right) =e^{i\langle \cdot
-x/2,p\rangle }\psi \left( \cdot -x\right) .
\end{equation*}%
Equivalently, 
\begin{equation*}
\mathcal{W}\left( \xi \right) =e^{-\frac{i}{2}\langle x,p\rangle
}e^{i\left\langle Q,p\right\rangle }e^{-i\left\langle x,P\right\rangle }=e^{%
\frac{i}{2}\langle x,p\rangle }e^{-i\left\langle x,P\right\rangle
}e^{i\left\langle Q,p\right\rangle }=e^{i\left( \left\langle
Q,p\right\rangle -\left\langle x,P\right\rangle \right) }.
\end{equation*}

\begin{remark}
$(\func{a})$ If $\left\{ \varepsilon _{1},...,\varepsilon _{n}\right\} $ is
a basis in $X$ and $\left\{ \varepsilon ^{1},...,\varepsilon ^{n}\right\} $
is the dual basis in $X^{\ast }$, then $\mathcal{B}=\left\{ \varepsilon
_{1},...,\varepsilon _{n},\varepsilon ^{1},...,\varepsilon ^{n}\right\} $ is
a symplectic basis and 
\begin{equation*}
P_{j}=-i\frac{\partial }{\partial x^{j}},\quad Q^{j}=M_{x^{j}},\quad
j=1,...,n,
\end{equation*}%
in the Schr\"{o}dinger representation. Here $M_{f}$ is the multiplication
operator by the function $f$.

$(\func{b})$ The space $\mathcal{S}=\mathcal{S}(\mathcal{H}\left( X\right) ,%
\mathcal{W})$ of the $\mathcal{C}^{\infty }$ vectors of the Schr\"{o}dinger
representation $(\mathcal{H}\left( X\right) ,\mathcal{W})$ is the space $%
\mathcal{S}(X)$ of tempered test functions.

$(\func{c})$ The Schr\"{o}dinger representation $(\mathcal{H}\left( X\right)
,\mathcal{W})$ is irreducible (this is just Lemma 7.1.4 in \cite{Hormander2}%
).
\end{remark}

We recall that a symplectic space has only one irreducible representation
(modulo unitary equivalence) and that each of its representations is a
multiple of this one (see Theorem 15 in \cite{Kastler}, Theorem C.38 in \cite%
{Raeburn}).

In the rest of the section we shall work in the representation described
above.

For a function $f$ on $\mathfrak{S}=T^{\ast }(X)$, we denote by $f\left( P_{%
\mathfrak{S}}\right) $ the operator $\mathcal{F}_{\mathfrak{S}}M_{f}\mathcal{%
F}_{\mathfrak{S}}$ on $\mathcal{H}\left( X\right) $. Let $\tau $ be an
endomorphism of $X$. If for $\left( x,p\right) \in T^{\ast }(X)$ we set $%
\theta _{X,\tau }\left( x,p\right) =\left\langle \tau x,p\right\rangle $
then we get a quadratic form on $T^{\ast }(X)$ which allows us to introduce
the following definition.

\begin{definition}
Let $a\in \mathcal{S}^{\ast }\left( X\times X^{\ast }\right) $. We can
define the operator $a_{X}^{\tau }\left( R\right) =a_{X}^{\tau }\left(
Q,P\right) :\mathcal{S}\left( X\right) \rightarrow \mathcal{S}^{\ast }\left(
X\right) $ by 
\begin{equation*}
a_{X}^{\tau }\left( R\right) =\int \mathcal{W}\left( \xi \right) \widehat{%
a_{X}^{\tau }}\left( \xi \right) \limfunc{d}\xi ,
\end{equation*}%
where $a_{X}^{\tau }=e^{i\theta _{X,\frac{1}{2}-\tau }\left( P_{\mathfrak{S}%
}\right) }a$. The above integral make sense if it is taken in the weak
sense, i.e. for $\varphi ,\psi \in \mathcal{S}\left( X\right) $ 
\begin{equation*}
\left\langle \varphi ,a_{X}^{\tau }\left( R\right) \psi \right\rangle _{%
\mathcal{S}\left( X\right) ,\mathcal{S}^{\ast }\left( X\right)
}=\left\langle \overline{\left\langle \varphi ,\mathcal{W}(\cdot )\psi
\right\rangle }_{\mathcal{S}\left( X\right) ,\mathcal{S}^{\ast }\left(
X\right) },\widehat{a_{X}^{\tau }}\right\rangle _{\mathcal{S}\left( 
\mathfrak{S}\right) ,\mathcal{S}^{\ast }\left( \mathfrak{S}\right) }.
\end{equation*}%
$a\in \mathcal{S}^{\ast }\left( X\times X^{\ast }\right) $ is called $\left(
X,\tau \right) $-symbol of $a_{X}^{\tau }\left( Q,P\right) $ and $%
a_{X}^{\tau }$ is the Weyl symbol of this operator.
\end{definition}

Let $a,b\in \widehat{\mathcal{M}}\left( \mathfrak{S}\right) $. Then the Weyl
symbol of the operator $a_{X}^{\tau }\left( Q,P\right) b_{X}^{\tau }\left(
Q,P\right) $ is $a_{X}^{\tau }\circ b_{X}^{\tau }$ while the $\left( X,\tau
\right) $-symbol is $e^{-i\left( \frac{1}{2}-\tau \right) \theta _{X}\left(
P_{\mathfrak{S}}\right) }a_{X}^{\tau }\circ b_{X}^{\tau }$ denoted by $a%
\overset{\tau }{\underset{X}{\circ }}b$ and called $\left( X,\tau \right) $%
-composition product.

We shall now check the action of $a_{X}^{\tau }\left( R\right) $ when $a\in 
\mathcal{S}\left( X\times X^{\ast }\right) $. For $\varphi \in \mathcal{S}%
\left( X\right) $ we have%
\begin{eqnarray*}
\left( a_{X}^{\tau }\left( R\right) \varphi \right) \left( x\right)
&=&\iint\limits_{X\times X^{\ast }}e^{i\theta _{X,\frac{1}{2}-\tau }\left(
y,p\right) }\mathcal{W}\left( y,p\right) \varphi \left( x\right) \widehat{a}%
\left( y,p\right) \limfunc{d}y\limfunc{d}p \\
&=&\iint\limits_{X\times X^{\ast }}e^{i\langle \left( \frac{1}{2}-\tau
\right) y,p\rangle }e^{i\langle x-\frac{1}{2}y,p\rangle }\varphi \left(
x-y\right) (\mathcal{F}_{X}\otimes \mathcal{F}_{X}^{-1})a\left( p,y\right) 
\limfunc{d}y\limfunc{d}p \\
&=&\iint\limits_{X\times X^{\ast }}e^{i\langle x-\tau y,p\rangle }\varphi
\left( x-y\right) (\mathcal{F}_{X}\otimes \mathcal{F}_{X}^{-1})a\left(
p,y\right) \limfunc{d}y\limfunc{d}p \\
&=&\int_{X}\varphi \left( x-y\right) \left( \mathcal{F}_{X}^{-1}\otimes 
\limfunc{id}\right) \left( \left( \mathcal{F}_{X}\otimes \mathcal{F}%
_{X}^{-1}\right) a\right) \left( x-\tau y,y\right) \limfunc{d}y \\
&=&\int_{X}\varphi \left( x-y\right) \left( \limfunc{id}\otimes \mathcal{F}%
_{X}^{-1}\right) a\left( x-\tau y,y\right) \limfunc{d}y \\
&=&\int_{X}\left( \limfunc{id}\otimes \mathcal{F}_{X}^{-1}\right) a\left(
\left( 1-\tau \right) x+\tau y,x-y\right) \varphi \left( y\right) \limfunc{d}%
y.
\end{eqnarray*}%
It follows that the kernel of $a_{X}^{\tau }\left( R\right) $ is given by 
\begin{equation}
\mathcal{K}_{a_{X}^{\tau }\left( R\right) }=\left( \left( \limfunc{id}%
\otimes \mathcal{F}_{X}^{-1}\right) a\right) \circ C_{\tau }  \label{SR9}
\end{equation}%
where $C_{\tau }$ is the map 
\begin{equation*}
C_{\tau }:X\times X\rightarrow X\times X,\quad C_{\tau }\left( x,y\right)
=\left( \left( 1-\tau \right) x+\tau y,x-y\right) .
\end{equation*}%
Let us note that (\ref{SR9}) is true for $a\in \mathcal{S}^{\ast }\left(
X\times X^{\ast }\right) $ because $\mathcal{S}\left( X\times X^{\ast
}\right) $ is dense in $\mathcal{S}^{\ast }\left( X\times X^{\ast }\right) $
in the weak topology and the mappings 
\begin{gather*}
\mathcal{S}^{\ast }\left( X\times X^{\ast }\right) \longrightarrow \mathcal{S%
}^{\ast }\left( X\times X^{\ast }\right) \longrightarrow \mathcal{B}(%
\mathcal{S}\left( X\right) ,\mathcal{S}^{\ast }\left( X\right) ), \\
a\rightarrow a_{X}^{\tau }=e^{i\theta _{X,\frac{1}{2}-\tau }\left( P_{%
\mathfrak{S}}\right) }a\rightarrow \limfunc{Op}\left( a_{X}^{\tau }\right)
=a_{X}^{\tau }\left( R\right) , \\
\mathcal{S}^{\ast }\left( X\times X^{\ast }\right) \longrightarrow \mathcal{S%
}^{\ast }\left( X\times X^{\ast }\right) ,\quad a\rightarrow \left( \left( 
\limfunc{id}\otimes \mathcal{F}_{X}^{-1}\right) a\right) \circ C_{\tau }
\end{gather*}%
are continuous if on $\mathcal{S}^{\ast }\left( X\times X^{\ast }\right) $
we consider the $\sigma (\mathcal{S}^{\ast }\left( X\times X^{\ast }\right) ,%
\mathcal{S}\left( X\times X^{\ast }\right) )$ topology and on $\mathcal{B}(%
\mathcal{S}\left( X\right) ,\mathcal{S}^{\ast }\left( X\right) )$ the
topology defined by the family of seminorms $\left\{ p_{\varphi ,\psi
}\right\} _{\varphi ,\psi \in \mathcal{S}\left( X\right) }$ 
\begin{equation*}
p_{\varphi ,\psi }\left( A\right) =\left\vert \left\langle \varphi ,A\psi
\right\rangle \right\vert ,\quad A\in \mathcal{B}(\mathcal{S}\left( X\right)
,\mathcal{S}^{\ast }\left( X\right) ).
\end{equation*}%
Since the equation in $a\in \mathcal{S}^{\ast }\left( X\times X^{\ast
}\right) $, $\left( \left( \limfunc{id}\otimes \mathcal{F}_{X}^{-1}\right)
a\right) \circ C_{\tau }=\mathcal{K}$, has a unique solution for each $%
\mathcal{K}\in \mathcal{S}^{\ast }\left( X\times X^{\ast }\right) $, a
consequence of the kernel theorem is the fact that the map 
\begin{equation*}
\mathcal{S}^{\ast }\left( X\times X^{\ast }\right) \rightarrow \mathcal{B}(%
\mathcal{S}\left( X\right) ,\mathcal{S}^{\ast }\left( X\right) ),\quad
a\rightarrow a_{X}^{\tau }\left( R\right)
\end{equation*}%
is linear, continuous and bijective. Hence to each $A\in \mathcal{B}(%
\mathcal{S}\left( X\right) ,\mathcal{S}^{\ast }\left( X\right) )$ we
associate a distribution $a\in \mathcal{S}^{\ast }\left( X\times X^{\ast
}\right) $ such that $A=a_{X}^{\tau }(R)$. This distribution is called $%
\left( X,\tau \right) $-symbol of $A$ and we shall use the notation $%
a=\sigma _{X}^{\tau }\left( A\right) $. When $\tau =\frac{1}{2}$ then $%
\sigma _{X}^{\frac{1}{2}}\left( A\right) $ is just the Weyl symbol of $A$.

\section{Kato's identity}

In this section we state and prove an extension of a formula due to T. Kato
which is a basic tool in this paper. For a finite dimensional vector space $%
E $ over $%
%TCIMACRO{\U{211d} }%
%BeginExpansion
\mathbb{R}
%EndExpansion
$, we shall use the notation $\mathcal{C}_{pol}^{\infty }\left( E\right) $
for the subalgebra of $\mathcal{C}^{\infty }\left( E\right) $ consisting of
functions whose derivatives have at most polynomial groth at infinity. We
shall need the following auxiliary result.

\begin{lemma}
Let $\left( \mathcal{H},\mathcal{W}\right) $ be a Weyl system associated to
the symplectic space $\mathfrak{S}$ and let $\varphi ,\psi \in \mathcal{S}$.

$(\func{a})$ If $a\in \mathcal{S}^{\ast }\left( \mathfrak{S}\right) $, then
the map 
\begin{equation*}
\mathfrak{S}\ni \xi \rightarrow \left\langle \varphi ,\mathcal{W}\left( \xi
\right) a\left( R\right) \mathcal{W}\left( -\xi \right) \psi \right\rangle _{%
\mathcal{S},\mathcal{S}^{\ast }}\in 
%TCIMACRO{\U{2102}}%
%BeginExpansion
\mathbb{C}%
%EndExpansion
\end{equation*}%
belongs to $\mathcal{C}_{pol}^{\infty }\left( \mathfrak{S}\right) $.

$(\func{b})$ If $a\in \mathcal{S}\left( \mathfrak{S}\right) $, then the map 
\begin{equation*}
\mathfrak{S}\ni \xi \rightarrow \left\langle \varphi ,\mathcal{W}\left( \xi
\right) a\left( R\right) \mathcal{W}\left( -\xi \right) \psi \right\rangle _{%
\mathcal{S},\mathcal{S}^{\ast }}\in 
%TCIMACRO{\U{2102}}%
%BeginExpansion
\mathbb{C}%
%EndExpansion
\end{equation*}%
belongs to $\mathcal{S}(\mathfrak{S})$.

$(\func{c})$ If $a\in \mathcal{S}^{\ast }\left( \mathfrak{S}\right) $ and $%
\xi \in \mathfrak{S}$, then 
\begin{equation*}
\mathcal{W}\left( \xi \right) a\left( R\right) \mathcal{W}\left( -\xi
\right) =\left( T_{\xi }a\right) \left( R\right) ,
\end{equation*}%
where $T_{\xi }a$ denote the translate by $\xi $ of the distribution $a$,
i.e. $\left( T_{\xi }a\right) \left( \eta \right) =a\left( \eta -\xi \right) 
$.
\end{lemma}

\begin{proof}
We know that $w=w_{\varphi ,\psi }=\left\langle \varphi ,\mathcal{W}(\cdot
)\psi \right\rangle _{\mathcal{S},\mathcal{S}^{\ast }}\in \mathcal{S}(%
\mathfrak{S})$. Assume that $a\in \mathcal{S}^{\ast }\left( \mathfrak{S}%
\right) $. Then from Corollary \ref{SR10}, (\ref{SR1}) and (\ref{SR2}) we get%
\begin{multline*}
\left\langle \varphi ,\mathcal{W}\left( \xi \right) a\left( R\right) 
\mathcal{W}\left( -\xi \right) \psi \right\rangle _{\mathcal{S},\mathcal{S}%
^{\ast }}=\left\langle \mathcal{W}\left( -\xi \right) \varphi ,a\left(
R\right) \mathcal{W}\left( -\xi \right) \psi \right\rangle _{\mathcal{S},%
\mathcal{S}^{\ast }} \\
=\left\langle \overline{\left\langle \mathcal{W}\left( -\xi \right) \varphi ,%
\mathcal{W}\left( \cdot \right) \mathcal{W}\left( -\xi \right) \psi
\right\rangle }_{\mathcal{S},\mathcal{S}^{\ast }},\widehat{a}\right\rangle _{%
\mathcal{S}\left( \mathfrak{S}\right) ,\mathcal{S}^{\ast }\left( \mathfrak{S}%
\right) } \\
=\left\langle \overline{\left\langle \varphi ,\mathcal{W}\left( \xi \right) 
\mathcal{W}\left( \cdot \right) \mathcal{W}\left( -\xi \right) \psi
\right\rangle }_{\mathcal{S},\mathcal{S}^{\ast }},\widehat{a}\right\rangle _{%
\mathcal{S}\left( \mathfrak{S}\right) ,\mathcal{S}^{\ast }\left( \mathfrak{S}%
\right) } \\
=\left\langle \overline{\left\langle \varphi ,e^{i\sigma \left( \xi ,\cdot
\right) }\mathcal{W}\left( \cdot \right) \psi \right\rangle }_{\mathcal{S},%
\mathcal{S}^{\ast }},\widehat{a}\right\rangle _{\mathcal{S}\left( \mathfrak{S%
}\right) ,\mathcal{S}^{\ast }\left( \mathfrak{S}\right) } \\
=\left\langle e^{-i\sigma \left( \xi ,\cdot \right) }\overline{\left\langle
\varphi ,\mathcal{W}\left( \cdot \right) \psi \right\rangle }_{\mathcal{S},%
\mathcal{S}^{\ast }},\widehat{a}\right\rangle _{\mathcal{S}\left( \mathfrak{S%
}\right) ,\mathcal{S}^{\ast }\left( \mathfrak{S}\right) } \\
=\left\langle e^{-i\sigma \left( \xi ,\cdot \right) }\overline{w},\widehat{a}%
\right\rangle _{\mathcal{S}\left( \mathfrak{S}\right) ,\mathcal{S}^{\ast
}\left( \mathfrak{S}\right) }=\left\langle \mathcal{F}_{\mathfrak{S}}\left(
e^{-i\sigma \left( \xi ,\cdot \right) }\overline{w}\right) ,a\right\rangle _{%
\mathcal{S}\left( \mathfrak{S}\right) ,\mathcal{S}^{\ast }\left( \mathfrak{S}%
\right) } \\
=\left\langle \widehat{\overline{w}}\left( \cdot +\xi \right)
,a\right\rangle _{\mathcal{S}\left( \mathfrak{S}\right) ,\mathcal{S}^{\ast
}\left( \mathfrak{S}\right) }=\left\langle \overline{\widehat{w}}\left(
-\cdot -\xi \right) ,a\right\rangle _{\mathcal{S}\left( \mathfrak{S}\right) ,%
\mathcal{S}^{\ast }\left( \mathfrak{S}\right) }=\left( a\ast \widehat{w}%
\right) \left( -\xi \right) .
\end{multline*}%
Hence 
\begin{equation}
\left\langle \varphi ,\mathcal{W}\left( \xi \right) a\left( R\right) 
\mathcal{W}\left( -\xi \right) \psi \right\rangle _{\mathcal{S},\mathcal{S}%
^{\ast }}=\left( a\ast \widehat{w}\right) \left( -\xi \right) ,\quad \xi \in 
\mathfrak{S}  \label{K1}
\end{equation}%
and $(\func{a})$ and $(\func{b})$ follows at once from this equality.

$(\func{c})$ We shall prove the equality for $a\in \mathcal{S}\left( 
\mathfrak{S}\right) $, then the general case follows by continuity. For $%
a\in \mathcal{S}\left( \mathfrak{S}\right) $ we have 
\begin{eqnarray*}
\mathcal{W}\left( \xi \right) a\left( R\right) \mathcal{W}\left( -\xi
\right) &=&\mathcal{W}\left( \delta _{\xi }\right) \mathcal{W}\left( 
\widehat{a}\right) \mathcal{W}\left( \delta _{-\xi }\right) =\mathcal{W}%
\left( \delta _{\xi }\times \widehat{a}\times \delta _{-\xi }\right) \\
&=&\mathcal{W}\left( e^{i\sigma \left( \xi ,\cdot \right) }\widehat{a}%
\right) =\mathcal{W}\left( \widehat{T_{\xi }a}\right) =\left( T_{\xi
}a\right) \left( R\right)
\end{eqnarray*}%
and the proof is complete.
\end{proof}

We shall now prove two extensions of some important identities due to T.
Kato.

\begin{lemma}[Kato]
\label{K3}Let $\left( \mathcal{H},\mathcal{W}\right) $ be a Weyl system
associated to the symplectic space $\mathfrak{S}$.

$(\func{a})$ If $b\in \mathcal{S}\left( \mathfrak{S}\right) $, $c\in 
\mathcal{S}^{\ast }\left( \mathfrak{S}\right) $, then $b\ast c\in \mathcal{S}%
^{\ast }\left( \mathfrak{S}\right) $ and 
\begin{eqnarray*}
\left( b\ast c\right) \left( R\right) &=&\int_{\mathfrak{S}}b\left( \xi
\right) \mathcal{W}\left( \xi \right) c\left( R\right) \mathcal{W}\left(
-\xi \right) \limfunc{d}\xi \\
&=&\int_{\mathfrak{S}}c\left( \xi \right) \mathcal{W}\left( \xi \right)
b\left( R\right) \mathcal{W}\left( -\xi \right) \limfunc{d}\xi ,
\end{eqnarray*}%
where the first integral is weakly absolutely convergent while the second
one must be interpreted in the sese of distributions and represents the
operator defined by%
\begin{multline*}
\left\langle \varphi ,\left( \int_{\mathfrak{S}}c\left( \xi \right) \mathcal{%
W}\left( \xi \right) b\left( R\right) \mathcal{W}\left( -\xi \right) 
\limfunc{d}\xi \right) \psi \right\rangle _{\mathcal{S},\mathcal{S}^{\ast }}
\\
=\left\langle \overline{\left\langle \varphi ,\mathcal{W}\left( \cdot
\right) b\left( R\right) \mathcal{W}\left( -\cdot \right) \psi \right\rangle 
}_{\mathcal{S},\mathcal{S}^{\ast }},c\right\rangle _{\mathcal{S}\left( 
\mathfrak{S}\right) ,\mathcal{S}^{\ast }\left( \mathfrak{S}\right) },
\end{multline*}%
for all $\varphi $, $\psi \in \mathcal{S}$.

$(\func{b})$ Let $h\in \mathcal{C}_{pol}^{\infty }\left( \mathfrak{S}\right) 
$. If $b\in L^{p}\left( \mathfrak{S}\right) $ and $c\in L^{q}\left( 
\mathfrak{S}\right) $, where $1\leq p,q\leq \infty $ and $p^{-1}+q^{-1}\geq
1 $, then $b\ast c\in L^{r}\left( \mathfrak{S}\right) $, $%
r^{-1}=p^{-1}+q^{-1}-1$ and 
\begin{equation}
\left( h\left( P_{\mathfrak{S}}\right) \left( b\ast c\right) \right) \left(
R\right) =\int_{\mathfrak{S}}b\left( \xi \right) \mathcal{W}\left( \xi
\right) \left( h\left( P_{\mathfrak{S}}\right) c\right) \left( R\right) 
\mathcal{W}\left( -\xi \right) \limfunc{d}\xi ,  \label{K2}
\end{equation}%
where the integral is weakly absolutely convergent.
\end{lemma}

\begin{proof}
$(\func{a})$ Let $\varphi $, $\psi \in \mathcal{S}$. Then $w=w_{\varphi
,\psi }=\left\langle \varphi ,\mathcal{W}(\cdot )\psi \right\rangle _{%
\mathcal{S},\mathcal{S}^{\ast }}\in \mathcal{S}(\mathfrak{S})$. First we
consider the case when $b$, $c\in \mathcal{S}\left( \mathfrak{S}\right) $.
Then 
\begin{eqnarray*}
\left\langle \varphi ,\left( b\ast c\right) \left( R\right) \psi
\right\rangle _{\mathcal{S},\mathcal{S}^{\ast }} &=&\int_{\mathfrak{S}}%
\widehat{b\ast c}\left( \eta \right) \left\langle \varphi ,\mathcal{W}\left(
\eta \right) \psi \right\rangle _{\mathcal{S},\mathcal{S}^{\ast }}\limfunc{d}%
\eta =\int_{\mathfrak{S}}\widehat{b\ast c}\left( \eta \right) w\left( \eta
\right) \limfunc{d}\eta \\
&=&\int_{\mathfrak{S}}\left( b\ast c\right) \left( \eta \right) \widehat{w}%
\left( -\eta \right) \limfunc{d}\eta \\
&=&\int_{\mathfrak{S}}\left( \int_{\mathfrak{S}}b\left( \xi \right) c\left(
\eta -\xi \right) \limfunc{d}\xi \right) \widehat{w}\left( -\eta \right) 
\limfunc{d}\eta \\
&=&\int_{\mathfrak{S}}b\left( \xi \right) \left( \int_{\mathfrak{S}}c\left(
\eta -\xi \right) \widehat{w}\left( -\eta \right) \limfunc{d}\eta \right) 
\limfunc{d}\xi \\
&=&\int_{\mathfrak{S}}b\left( \xi \right) \left( \int_{\mathfrak{S}}\left(
T_{\xi }c\right) \left( \eta \right) \widehat{w}\left( -\eta \right) 
\limfunc{d}\eta \right) \limfunc{d}\xi \\
&=&\int_{\mathfrak{S}}b\left( \xi \right) \left( \int_{\mathfrak{S}}\widehat{%
T_{\xi }c}\left( \eta \right) w\left( \eta \right) \limfunc{d}\eta \right) 
\limfunc{d}\xi \\
&=&\int_{\mathfrak{S}}b\left( \xi \right) \left( \int_{\mathfrak{S}}\widehat{%
T_{\xi }c}\left( \eta \right) \left\langle \varphi ,\mathcal{W}\left( \eta
\right) \psi \right\rangle _{\mathcal{S},\mathcal{S}^{\ast }}\limfunc{d}\eta
\right) \limfunc{d}\xi \\
&=&\int_{\mathfrak{S}}b\left( \xi \right) \left\langle \varphi ,\left(
T_{\xi }c\right) \left( R\right) \psi \right\rangle _{\mathcal{S},\mathcal{S}%
^{\ast }}\limfunc{d}\xi \\
&=&\int_{\mathfrak{S}}b\left( \xi \right) \left\langle \varphi ,\mathcal{W}%
\left( \xi \right) c\left( R\right) \mathcal{W}\left( -\xi \right) \psi
\right\rangle _{\mathcal{S},\mathcal{S}^{\ast }}\limfunc{d}\xi ,
\end{eqnarray*}%
where in the last equality we used the formula $\left( T_{\xi }c\right)
\left( R\right) =\mathcal{W}\left( \xi \right) c\left( R\right) \mathcal{W}%
\left( -\xi \right) $.

Let $c\in \mathcal{S}^{\ast }\left( \mathfrak{S}\right) $ and let $\left\{
c_{j}\right\} \subset \mathcal{S}(\mathfrak{S})$ be such that $%
c_{j}\rightarrow c$ weakly in $\mathcal{S}^{\ast }\left( \mathfrak{S}\right) 
$. The uniform boundedness principle and Peetre's inequality imply that

\begin{itemize}
\item[-] $\left\langle \varphi ,\mathcal{W}\left( \xi \right) a_{j}\left(
R\right) \mathcal{W}\left( -\xi \right) \psi \right\rangle _{\mathcal{S},%
\mathcal{S}^{\ast }}\rightarrow \left\langle \varphi ,\mathcal{W}\left( \xi
\right) a\left( R\right) \mathcal{W}\left( -\xi \right) \psi \right\rangle _{%
\mathcal{S},\mathcal{S}^{\ast }},\quad \xi \in \mathfrak{S}.$

\item[-] There are $M\in 
%TCIMACRO{\U{2115} }%
%BeginExpansion
\mathbb{N}
%EndExpansion
$\textit{, }$C=C\left( M,w\right) >0$\textit{\ }such that\textit{\ }%
\begin{equation*}
\left\vert \left\langle \varphi ,\mathcal{W}\left( \xi \right) c_{j}\left(
R\right) \mathcal{W}\left( -\xi \right) \psi \right\rangle _{\mathcal{S},%
\mathcal{S}^{\ast }}\right\vert \leq C\left\langle \xi \right\rangle
^{M},\quad \xi \in \mathfrak{S},
\end{equation*}%
where $\left\langle \xi \right\rangle =\left( 1+\left\vert \xi \right\vert
^{2}\right) ^{\frac{1}{2}}$ and $\left\vert \cdot \right\vert $ is an
euclidean norm on $\mathfrak{S}$.
\end{itemize}

\noindent The general case can be deduced from the above case if we observe
that%
\begin{equation*}
\left\langle \varphi ,\left( b\ast c_{j}\right) \left( R\right) \psi
\right\rangle _{\mathcal{S},\mathcal{S}^{\ast }}\rightarrow \left\langle
\varphi ,\left( b\ast c\right) \left( R\right) \psi \right\rangle _{\mathcal{%
S},\mathcal{S}^{\ast }},
\end{equation*}
\begin{multline*}
\left\langle \overline{\left\langle \varphi ,\mathcal{W}\left( \cdot \right)
b\left( R\right) \mathcal{W}\left( -\cdot \right) \psi \right\rangle }_{%
\mathcal{S},\mathcal{S}^{\ast }},c_{j}\right\rangle _{\mathcal{S}\left( 
\mathfrak{S}\right) ,\mathcal{S}^{\ast }\left( \mathfrak{S}\right) } \\
\rightarrow \left\langle \overline{\left\langle \varphi ,\mathcal{W}\left(
\cdot \right) b\left( R\right) \mathcal{W}\left( -\cdot \right) \psi
\right\rangle }_{\mathcal{S},\mathcal{S}^{\ast }},c\right\rangle _{\mathcal{S%
}\left( \mathfrak{S}\right) ,\mathcal{S}^{\ast }\left( \mathfrak{S}\right) }
\end{multline*}%
and that the sequence $\left\{ b\left( \cdot \right) \left\langle \varphi ,%
\mathcal{W}\left( \cdot \right) c_{j}\left( R\right) \mathcal{W}\left(
-\cdot \right) \psi \right\rangle _{\mathcal{S},\mathcal{S}^{\ast }}\right\} 
$ converge dominated to $b\left( \cdot \right) \left\langle \varphi ,%
\mathcal{W}\left( \cdot \right) c\left( R\right) \mathcal{W}\left( -\cdot
\right) \psi \right\rangle _{\mathcal{S},\mathcal{S}^{\ast }}$.

$(\func{b})$ We recall the Young inequality. If $1\leq p,q\leq \infty $, $%
p^{-1}+q^{-1}\geq 1$, $r^{-1}=p^{-1}+q^{-1}-1$, $b\in L^{p}\left( \mathfrak{S%
}\right) $ and $c\in L^{q}\left( \mathfrak{S}\right) $, then $b\ast c\in
L^{r}\left( \mathfrak{S}\right) $ and 
\begin{equation*}
\left\Vert b\ast c\right\Vert _{L^{r}\left( \mathfrak{S}\right) }\leq
\left\Vert b\right\Vert _{L^{p}\left( \mathfrak{S}\right) }\left\Vert
c\right\Vert _{L^{q}\left( \mathfrak{S}\right) }.
\end{equation*}%
Let $p^{\prime },r^{\prime }\geq 1$ such that $p^{-1}+p^{\prime
-1}=r^{-1}+r^{\prime -1}=1$. Then $r^{\prime -1}+q^{-1}\geq 1$ and $%
p^{\prime -1}=r^{\prime -1}+q^{-1}-1$.

If $b\in \mathcal{S}\left( \mathfrak{S}\right) $ and $c\in L^{q}\left( 
\mathfrak{S}\right) $, then $g=h\left( P_{\mathfrak{S}}\right) c\in \mathcal{%
S}^{\ast }\left( \mathfrak{S}\right) $ and $h\left( P_{\mathfrak{S}}\right)
\left( b\ast c\right) =b\ast h\left( P_{\mathfrak{S}}\right) c=b\ast g$.
Using $(\func{a})$ it follows that%
\begin{equation*}
\left\langle \varphi ,\left( b\ast g\right) \left( R\right) \psi
\right\rangle _{\mathcal{S},\mathcal{S}^{\ast }}=\int_{\mathfrak{S}}b\left(
\xi \right) \left\langle \varphi ,\mathcal{W}\left( \xi \right) g\left(
R\right) \mathcal{W}\left( -\xi \right) \psi \right\rangle _{\mathcal{S},%
\mathcal{S}^{\ast }}\limfunc{d}\xi .
\end{equation*}%
Similarly, if $b\in L^{p}\left( \mathfrak{S}\right) \subset \mathcal{S}%
^{\ast }\left( \mathfrak{S}\right) $ and $c\in \mathcal{S}\left( \mathfrak{S}%
\right) $, then $g=h\left( P_{\mathfrak{S}}\right) c\in \mathcal{S}\left( 
\mathfrak{S}\right) $ and $h\left( P_{\mathfrak{S}}\right) \left( b\ast
c\right) =b\ast h\left( P_{\mathfrak{S}}\right) c=b\ast g$. Using $(\func{a}%
) $ again, it follows that 
\begin{eqnarray*}
\left\langle \varphi ,\left( b\ast g\right) \left( R\right) \psi
\right\rangle _{\mathcal{S},\mathcal{S}^{\ast }} &=&\left\langle \overline{%
\left\langle \varphi ,\mathcal{W}\left( \cdot \right) g\left( R\right) 
\mathcal{W}\left( -\cdot \right) \psi \right\rangle }_{\mathcal{S},\mathcal{S%
}^{\ast }},b\right\rangle _{\mathcal{S}\left( \mathfrak{S}\right) ,\mathcal{S%
}^{\ast }\left( \mathfrak{S}\right) } \\
&=&\int_{\mathfrak{S}}b\left( \xi \right) \left\langle \varphi ,\mathcal{W}%
\left( \xi \right) g\left( R\right) \mathcal{W}\left( -\xi \right) \psi
\right\rangle _{\mathcal{S},\mathcal{S}^{\ast }}\limfunc{d}\xi .
\end{eqnarray*}%
Hence we proved $(\func{b})$ in the case when either $b\in \mathcal{S}\left( 
\mathfrak{S}\right) $ or $c\in \mathcal{S}\left( \mathfrak{S}\right) $.

The general case can be obtained from these particular cases by an
approximation argument. Observe that condition $p^{-1}+q^{-1}\geq 1$ implies
that $p<\infty $ or $q<\infty $. For the convergence of the left-hand side
of (\ref{K2}) we use the Young inequality and the continuity of the map $%
\limfunc{Op}$ in Corollary \ref{SR10}. To estimate right-hand side of (\ref%
{K2}) we use (\ref{K1}) and the H\"{o}lder and Young inequalities. We have%
\begin{equation*}
\left\langle \varphi ,\mathcal{W}\left( \xi \right) \left( h\left( P_{%
\mathfrak{S}}\right) c\right) \left( R\right) \mathcal{W}\left( -\xi \right)
\psi \right\rangle _{\mathcal{S},\mathcal{S}^{\ast }}=\left( h\left( P_{%
\mathfrak{S}}\right) c\ast \widehat{w}\right) \left( -\xi \right) =\left(
c\ast \widehat{hw}\right) \left( -\xi \right)
\end{equation*}%
and%
\begin{multline*}
\left\vert \left\langle \varphi ,\left( \int_{\mathfrak{S}}b\left( \xi
\right) \mathcal{W}\left( \xi \right) \left( h\left( P_{\mathfrak{S}}\right)
c\right) \left( R\right) \mathcal{W}\left( -\xi \right) \limfunc{d}\xi
\right) \psi \right\rangle _{\mathcal{S},\mathcal{S}^{\ast }}\right\vert \\
=\left\vert \left( \int_{\mathfrak{S}}b\left( \xi \right) \left( c\ast 
\widehat{hw}\right) \left( -\xi \right) \limfunc{d}\xi \right) \right\vert \\
\leq \left\Vert b\right\Vert _{L^{p}\left( \mathfrak{S}\right) }\left\Vert
c\ast \widehat{hw}\right\Vert _{L^{p^{\prime }}\left( \mathfrak{S}\right) }
\\
\leq \left\Vert b\right\Vert _{L^{p}\left( \mathfrak{S}\right) }\left\Vert
c\right\Vert _{L^{q}\left( \mathfrak{S}\right) }\left\Vert \widehat{hw}%
\right\Vert _{L^{r^{\prime }}\left( \mathfrak{S}\right) },
\end{multline*}%
where $w=w_{\varphi ,\psi }=\left\langle \varphi ,\mathcal{W}(\cdot )\psi
\right\rangle _{\mathcal{S},\mathcal{S}^{\ast }}\in \mathcal{S}(\mathfrak{S}%
) $.
\end{proof}

\begin{corollary}
\label{K4}Let $\left( \mathcal{H},\mathcal{W}\right) $ be a Weyl system
associated to the symplectic space $\mathfrak{S}$. If $b\in L^{p}\left( 
\mathfrak{S}\right) $ and $c\in L^{q}\left( \mathfrak{S}\right) $, where $%
1\leq p,q\leq \infty $ and $p^{-1}+q^{-1}\geq 1$, then $b\ast c\in
L^{r}\left( \mathfrak{S}\right) $, $r^{-1}=p^{-1}+q^{-1}-1$ and 
\begin{equation*}
\left( b\ast c\right) \left( R\right) =\int_{\mathfrak{S}}b\left( \xi
\right) \mathcal{W}\left( \xi \right) c\left( R\right) \mathcal{W}\left(
-\xi \right) \limfunc{d}\xi ,
\end{equation*}%
where the integral is weakly absolutely convergent.

If $\mathfrak{S}=T^{\ast }(X)=X\times X^{\ast }$ and $\tau $ an endomorphism
of $X$, then we have%
\begin{equation*}
\left( b\ast c\right) _{X}^{\tau }\left( R\right) =\int_{\mathfrak{S}%
}b\left( \xi \right) \mathcal{W}\left( \xi \right) c_{X}^{\tau }\left(
R\right) \mathcal{W}\left( -\xi \right) \limfunc{d}\xi ,
\end{equation*}%
where the integral is weakly absolutely convergent.
\end{corollary}

\begin{proof}
We take $h=e^{i\theta _{X,\frac{1}{2}-\tau }\left( \cdot \right) }$ in the
previous lemma.
\end{proof}

\section{Kato's operator calculus}

In \cite{Cordes}, H.O. Cordes noticed that the $L^{2}$-boundedness of an
operator $a\left( x,D\right) $ in $OPS_{0,0}^{0}$ could be deduced by a
synthesis of $a\left( x,D\right) =a\left( R\right) $ from trace-class
operators. In \cite{Kato}, T. Kato extended this argument to the general
case $OPS_{\rho ,\rho }^{0}$, $0<\rho <1$, and abstracted the functional
analysis involved in Cordes' argument. This operator calculus can be
extended further to investigate the Schatten-class properties of
pseudo-differential operators in $OPS_{0,0}^{0}$.

Let $\mathcal{H}$ be a separable Hilbert space. For $1\leq p<\infty $, we
denote by $\mathcal{B}_{p}\left( \mathcal{H}\right) $ the Scatten ideal of
compact operators on $\mathcal{H}$ whose singular values lie in $l^{p}$ with
the associated norm $\left\Vert \cdot \right\Vert _{p}$. For $p=\infty $, $%
\mathcal{B}_{\infty }\left( \mathcal{H}\right) $ is the ideal of compact
operators on $\mathcal{H}$ with $\left\Vert \cdot \right\Vert _{\infty
}=\left\Vert \cdot \right\Vert $.

\begin{definition}
Let $T,A,B\in \mathcal{B}\left( \mathcal{H}\right) $, $A\geq 0$, $B\geq 0$.
We write 
\begin{equation*}
T\ll \left( A;B\right) \overset{def}{\Longleftrightarrow }\left\vert \left(
u,Tv\right) \right\vert ^{2}\leq \left( u,Au\right) \left( v,Bv\right)
,\quad \text{for }u,v\in \mathcal{H}.
\end{equation*}
\end{definition}

\begin{lemma}
\label{KOC1}Let $S,T,A,B\in \mathcal{B}\left( \mathcal{H}\right) $, $A\geq 0$%
, $B\geq 0$. Then

$(\func{i})$ $T\ll \left( \left| T^{\ast }\right| ;\left| T\right| \right) $.

$(\func{ii})$ $T\ll \left( A;B\right) \Rightarrow T^{\ast }\ll \left(
B;A\right) $.

$\left( \func{iii}\right) $ $T\ll \left( A;B\right) \Rightarrow S^{\ast
}TS\ll \left( S^{\ast }AS;S^{\ast }BS\right) .$

$\left( \func{iv}\right) $ Let $\left\{ T_{j}\right\} ,\left\{ A_{j}\right\}
,\left\{ B_{j}\right\} \subset \mathcal{B}\left( \mathcal{H}\right) $, $%
A_{j}\geq 0$, $B_{j}\geq 0$, $j=1,2,...$. If $T_{j}\ll \left(
A_{j};B_{j}\right) $, $j=1,2,...$, then%
\begin{equation*}
\sum T_{j}\ll \left( \sum A_{j};\sum B_{j}\right)
\end{equation*}%
in the sense that whenever the series $\sum A_{j}$ and $\sum B_{j}$ converge
in the strong sense, the same is true for the serie $\sum T_{j}$ and the
inequality holds.
\end{lemma}

\begin{proof}
$(\func{i})$ If $T\geq 0$, then $T=\left\vert T\right\vert =\left\vert
T^{\ast }\right\vert $ and 
\begin{equation*}
\left\vert \left( u,Tv\right) \right\vert ^{2}=\left\vert \left( T^{\frac{1}{%
2}}u,T^{\frac{1}{2}}v\right) \right\vert ^{2}\leq \left\Vert T^{\frac{1}{2}%
}u\right\Vert ^{2}\left\Vert T^{\frac{1}{2}}v\right\Vert ^{2}=\left(
u,Tu\right) \left( v,Tv\right) ,\quad u,v\in \mathcal{H}.
\end{equation*}

In the general case, $T\in \mathcal{B}\left( \mathcal{H}\right) $, we shall
use the polar decomposition of $A$. Let $T=V\left\vert T\right\vert $ with $%
\left\vert T\right\vert =\left( T^{\ast }T\right) ^{\frac{1}{2}}$and $V\in 
\mathcal{B}\left( \mathcal{H}\right) $ a partial isometry such that $%
\limfunc{Ker}V=\limfunc{Ker}T$. Then $V^{\ast }V$ is the projection onto the
initial space of $V$, $\left( \limfunc{Ker}V\right) ^{\bot }(=\left( 
\limfunc{Ker}T\right) ^{\bot }=\overline{\limfunc{Ran}T^{\ast }})$. It
follows that $\left\vert T^{\ast }\right\vert =V\left\vert T\right\vert
V^{\ast }$ since $T^{\ast }=V^{\ast }VT^{\ast }$, $V\left\vert T\right\vert
V^{\ast }\geq 0$ and 
\begin{equation*}
\left( V\left\vert T\right\vert V^{\ast }\right) ^{2}=V\left\vert
T\right\vert V^{\ast }V\left\vert T\right\vert V^{\ast }=TV^{\ast }VT^{\ast
}=TT^{\ast }=\left\vert T^{\ast }\right\vert ^{2}.
\end{equation*}%
Then we have%
\begin{eqnarray*}
\left\vert \left( u,Tv\right) \right\vert ^{2} &=&\left\vert \left(
u,V\left\vert T\right\vert v\right) \right\vert ^{2}=\left\vert \left(
V^{\ast }u,\left\vert T\right\vert v\right) \right\vert ^{2}\leq \left(
V^{\ast }u,\left\vert T\right\vert V^{\ast }u\right) \left( v,\left\vert
T\right\vert v\right) \\
&=&\left( u,\left\vert T^{\ast }\right\vert u\right) \left( v,\left\vert
T\right\vert v\right) ,\quad u,v\in \mathcal{H}.
\end{eqnarray*}

$(\func{ii})$, $\left( \func{iii}\right) $ are obvious.

$\left( \func{iv}\right) $ Assume that there are $A,B\in \mathcal{B}\left( 
\mathcal{H}\right) $ such that $A=s-\lim_{n\rightarrow \infty
}\sum_{j=1}^{n}A_{j}$, $B=s-\lim_{n\rightarrow \infty }\sum_{j=1}^{n}B_{j}$.
For $n\geq 1$, we set $A\left( n\right) =\sum_{j=1}^{n}A_{j}$, $B\left(
n\right) =\sum_{j=1}^{n}B_{j}$, $T\left( n\right) =\sum_{j=1}^{n}T_{j}$ and
for $n=0$, we set $A\left( 0\right) =B\left( 0\right) =T\left( 0\right) =0$.
Since $A=s-\lim_{n\rightarrow \infty }A\left( n\right) $ by the uniform
boundedness principle there is a constant $C>0$ such that $\left\Vert
A\left( n\right) \right\Vert \leq C$ for all $n$ in $%
%TCIMACRO{\U{2115} }%
%BeginExpansion
\mathbb{N}
%EndExpansion
$.

Let $m>n.$ Then 
\begin{eqnarray*}
\left\vert \left( u,\left( T\left( m\right) -T\left( n\right) \right)
v\right) \right\vert &\leq &\sum_{j=n+1}^{m}\left\vert \left(
u,T_{j}v\right) \right\vert \leq \sum_{j=n+1}^{m}\left( u,A_{j}u\right) ^{%
\frac{1}{2}}\left( v,B_{j}v\right) ^{\frac{1}{2}} \\
&\leq &\left( \sum_{j=n+1}^{m}\left( u,A_{j}u\right) \right) ^{\frac{1}{2}%
}\left( \sum_{j=n+1}^{m}\left( v,B_{j}v\right) \right) ^{\frac{1}{2}} \\
&=&\left( u,\left( A\left( m\right) -A\left( n\right) \right) u\right) ^{%
\frac{1}{2}}\left( v,\left( B\left( m\right) -B\left( n\right) \right)
v\right) ^{\frac{1}{2}} \\
&\leq &\left( u,A\left( m\right) u\right) ^{\frac{1}{2}}\left( v,\left(
B\left( m\right) -B\left( n\right) \right) v\right) ^{\frac{1}{2}} \\
&\leq &C^{\frac{1}{2}}\left\Vert u\right\Vert \left\Vert v\right\Vert ^{%
\frac{1}{2}}\left\Vert \left( B\left( m\right) -B\left( n\right) \right)
v\right\Vert ^{\frac{1}{2}},\quad u,v\in \mathcal{H}.
\end{eqnarray*}%
Thus 
\begin{equation*}
\left\Vert T\left( m\right) v-T\left( n\right) v\right\Vert \leq C^{\frac{1}{%
2}}\left\Vert v\right\Vert ^{\frac{1}{2}}\left\Vert \left( B\left( m\right)
-B\left( n\right) \right) v\right\Vert ^{\frac{1}{2}},\quad v\in \mathcal{H},
\end{equation*}%
which implies that $\left\{ T\left( n\right) v\right\} $ is a Cauchy
sequence for any $v\in \mathcal{H}$, so there is $T\in \mathcal{B}\left( 
\mathcal{H}\right) $ such that $T=s-\lim_{n\rightarrow \infty }T\left(
n\right) $. By passing to the limit in the estimate%
\begin{equation*}
\left\vert \left( u,T\left( m\right) v\right) \right\vert \leq \left(
u,A\left( m\right) u\right) ^{\frac{1}{2}}\left( v,B\left( m\right) v\right)
^{\frac{1}{2}},\quad u,v\in \mathcal{H},
\end{equation*}%
we conclude that $T\ll \left( A;B\right) $.
\end{proof}

\begin{lemma}
Let $Y$ be a measure space and $Y\ni y\rightarrow U\left( y\right) \in 
\mathcal{B}\left( \mathcal{H}\right) $ a weakly measurable map.

$(\func{a})$ Assume that there is $C>0$ such that 
\begin{equation*}
\int_{Y}\left\vert \left( \varphi ,U\left( y\right) \psi \right) \right\vert
^{2}\limfunc{d}y\leq C\left\Vert \varphi \right\Vert ^{2}\left\Vert \psi
\right\Vert ^{2},\quad \varphi ,\psi \in \mathcal{H}.
\end{equation*}%
If $b\in L^{\infty }\left( Y\right) $ and $G\in \mathcal{B}_{1}\left( 
\mathcal{H}\right) $, then the integral%
\begin{equation*}
b\left\{ G\right\} =\int_{Y}b\left( y\right) U\left( y\right) ^{\ast
}GU\left( y\right) \limfunc{d}y
\end{equation*}%
is weakly absolutely convergent and defines a bounded operator such that%
\begin{equation*}
\left\Vert b\left\{ G\right\} \right\Vert \leq C\left\Vert b\right\Vert
_{L^{\infty }}\left\Vert G\right\Vert _{1}.
\end{equation*}

$(\func{b})$ Assume that there is $C>0$ such that 
\begin{equation*}
\left\Vert U\left( y\right) \right\Vert \leq C^{\frac{1}{2}}\quad a.e.\ y\in
Y.
\end{equation*}%
If $b\in L^{1}\left( Y\right) $ and $G\in \mathcal{B}_{1}\left( \mathcal{H}%
\right) $, then the integral%
\begin{equation*}
b\left\{ G\right\} =\int_{Y}b\left( y\right) U\left( y\right) ^{\ast
}GU\left( y\right) \limfunc{d}y
\end{equation*}%
is absolutely convergent and defines a trace class operator such that%
\begin{equation*}
\left\Vert b\left\{ G\right\} \right\Vert _{1}\leq C\left\Vert b\right\Vert
_{L^{1}}\left\Vert G\right\Vert _{1}.
\end{equation*}

$(\func{c})$ Assume that there is $C>0$ such that 
\begin{eqnarray*}
\int_{Y}\left\vert \left( \varphi ,U\left( y\right) \psi \right) \right\vert
^{2}\limfunc{d}y &\leq &C\left\Vert \varphi \right\Vert ^{2}\left\Vert \psi
\right\Vert ^{2},\quad \varphi ,\psi \in \mathcal{H}\quad \text{and} \\
\left\Vert U\left( y\right) \right\Vert &\leq &C^{\frac{1}{2}}\quad a.e.\
y\in Y.
\end{eqnarray*}%
If $b\in L^{p}\left( Y\right) $ with $1\leq p<\infty $ and $G\in \mathcal{B}%
_{1}\left( \mathcal{H}\right) $, then the integral%
\begin{equation*}
b\left\{ G\right\} =\int_{Y}b\left( y\right) U\left( y\right) ^{\ast
}GU\left( y\right) \limfunc{d}y
\end{equation*}%
is weakly absolutely convergent and defines an operator $b\left\{ G\right\} $
in $\mathcal{B}_{p}\left( \mathcal{H}\right) $ which satisfies%
\begin{equation*}
\left\Vert b\left\{ G\right\} \right\Vert _{p}\leq C\left\Vert b\right\Vert
_{L^{p}}\left\Vert G\right\Vert _{1}.
\end{equation*}
\end{lemma}

\begin{proof}
$(\func{a})$ We do this in several steps.

\textit{Step }1. Suppose $G\geq 0$. Write $G=\sum_{j=1}^{\infty }\lambda
_{j}\left\vert \varphi _{j})(\varphi _{j}\right\vert =\sum_{j=1}^{\infty
}\lambda _{j}\left( \varphi _{j},\cdot \right) \varphi _{j}$ with $\left\{
\varphi _{j}\right\} $ an orthonormal basis of $\mathcal{H}$, $\lambda
_{j}\geq 0$, $\limfunc{Tr}\left( G\right) =\sum_{j=1}^{\infty }\lambda
_{j}=\sum_{j=1}^{\infty }\left\vert \lambda _{j}\right\vert =\left\Vert
G\right\Vert _{1}$. Then 
\begin{equation*}
\left( U\left( y\right) \varphi ,GU\left( y\right) \psi \right)
=\sum_{j=1}^{\infty }\lambda _{j}\left( U\left( y\right) \varphi ,\varphi
_{j}\right) \left( \varphi _{j},U\left( y\right) \psi \right) ,\quad y\in
Y,\ \varphi ,\psi \in \mathcal{H}
\end{equation*}%
and%
\begin{eqnarray*}
\left\vert \left( \varphi ,B\psi \right) \right\vert &=&\left\vert
\int_{Y}b\left( y\right) \left( U\left( y\right) \varphi ,GU\left( y\right)
\psi \right) \limfunc{d}y\right\vert \\
&\leq &\left\Vert b\right\Vert _{L^{\infty }}\sum_{j=1}^{\infty }\lambda
_{j}\int_{Y}\left\vert \left( U\left( y\right) \varphi ,\varphi _{j}\right)
\left( \varphi _{j},U\left( y\right) \psi \right) \right\vert \limfunc{d}y \\
&\leq &\left\Vert b\right\Vert _{L^{\infty }}\sum_{j=1}^{\infty }\lambda
_{j}\left( \int_{Y}\left\vert \left( U\left( y\right) \varphi ,\varphi
_{j}\right) \right\vert ^{2}\limfunc{d}y\right) ^{\frac{1}{2}}\left(
\int_{Y}\left\vert \left( \varphi _{j},U\left( y\right) \psi \right)
\right\vert ^{2}\limfunc{d}y\right) ^{\frac{1}{2}} \\
&\leq &\left\Vert b\right\Vert _{L^{\infty }}\sum_{j=1}^{\infty }\lambda
_{j}C^{\frac{1}{2}}\left\Vert \varphi \right\Vert \left\Vert \varphi
_{j}\right\Vert C^{\frac{1}{2}}\left\Vert \psi \right\Vert \left\Vert
\varphi _{j}\right\Vert \\
&=&C\left\Vert b\right\Vert _{L^{\infty }}\left( \sum_{j=1}^{\infty }\lambda
_{j}\right) \left\Vert \varphi \right\Vert \left\Vert \psi \right\Vert
=C\left\Vert b\right\Vert _{L^{\infty }}\left\Vert G\right\Vert
_{1}\left\Vert \varphi \right\Vert \left\Vert \psi \right\Vert ,\quad
\varphi ,\psi \in \mathcal{H}.
\end{eqnarray*}%
\textit{Step }2. The general case, $G\in \mathcal{B}_{1}\left( \mathcal{H}%
\right) $, can be reduced to the above case by using Lemma \ref{KOC1} and
the Schwarz inequality to evaluate the integrand. We have $U\left( y\right)
^{\ast }GU\left( y\right) \ll \left( U\left( y\right) ^{\ast }\left\vert
G^{\ast }\right\vert U\left( y\right) ;U\left( y\right) ^{\ast }\left\vert
G\right\vert U\left( y\right) \right) $ and 
\begin{equation*}
\left\vert b\left( y\right) \left( U\left( y\right) \varphi ,GU\left(
y\right) \psi \right) \right\vert \leq \left\Vert b\right\Vert _{L^{\infty
}}\left( U\left( y\right) \varphi ,\left\vert G^{\ast }\right\vert U\left(
y\right) \varphi \right) ^{\frac{1}{2}}\left( U\left( y\right) \psi
,\left\vert G\right\vert U\left( y\right) \psi \right) ^{\frac{1}{2}}.
\end{equation*}%
Thus 
\begin{multline*}
\left\vert \left( \varphi ,b\left\{ G\right\} \psi \right) \right\vert
=\left\vert \int_{Y}b\left( y\right) \left( U\left( y\right) \varphi
,GU\left( y\right) \psi \right) \limfunc{d}y\right\vert \\
\leq \left\Vert b\right\Vert _{L^{\infty }}\left( \int_{Y}\left( U\left(
y\right) \varphi ,\left\vert G^{\ast }\right\vert U\left( y\right) \varphi
\right) \limfunc{d}y\right) ^{\frac{1}{2}}\left( \int_{Y}\left( U\left(
y\right) \psi ,\left\vert G\right\vert U\left( y\right) \psi \right) 
\limfunc{d}y\right) ^{\frac{1}{2}} \\
\leq \left\Vert b\right\Vert _{L^{\infty }}C^{\frac{1}{2}}\left\Vert
\left\vert G^{\ast }\right\vert \right\Vert _{1}^{\frac{1}{2}}\left\Vert
\varphi \right\Vert C^{\frac{1}{2}}\left\Vert \left\vert G\right\vert
\right\Vert _{1}^{\frac{1}{2}}\left\Vert \psi \right\Vert \\
=C\left\Vert b\right\Vert _{L^{\infty }}\left\Vert G\right\Vert
_{1}\left\Vert \varphi \right\Vert \left\Vert \psi \right\Vert ,\quad
\varphi ,\psi \in \mathcal{H}.
\end{multline*}

$(\func{b})$ Let $\left\{ \psi _{\alpha }\right\} _{\alpha \in I}$ and $%
\left\{ \varphi _{\alpha }\right\} _{\alpha \in I}$ be two orthonormal
systems in $\mathcal{H}$. Then we have 
\begin{eqnarray*}
\sum_{\alpha \in I}\left\vert \left( \psi _{\alpha },b\left\{ G\right\}
\varphi _{\alpha }\right) \right\vert &\leq &\sum_{\alpha \in
I}\int_{Y}\left\vert b\left( y\right) \right\vert \left\vert \left( \psi
_{\alpha },U\left( y\right) ^{\ast }GU\left( y\right) \varphi _{\alpha
}\right) \right\vert \limfunc{d}y \\
&=&\int_{Y}\left\vert b\left( y\right) \right\vert \sum_{\alpha \in
I}\left\vert \left( \psi _{\alpha },U\left( y\right) ^{\ast }GU\left(
y\right) \varphi _{\alpha }\right) \right\vert \limfunc{d}y \\
&\leq &\int_{Y}\left\vert b\left( y\right) \right\vert \left\Vert U\left(
y\right) ^{\ast }GU\left( y\right) \right\Vert _{1}\limfunc{d}y \\
&\leq &\int_{Y}\left\vert b\left( y\right) \right\vert \left\Vert U\left(
y\right) \right\Vert ^{2}\left\Vert G\right\Vert _{1}\limfunc{d}y \\
&\leq &C\left\Vert G\right\Vert _{1}\int_{Y}\left\vert b\left( y\right)
\right\vert \limfunc{d}y=C\left\Vert b\right\Vert _{L^{1}}\left\Vert
G\right\Vert _{1}<+\infty .
\end{eqnarray*}%
It follows that 
\begin{equation*}
\sum_{\alpha \in I}\left\vert \left( \psi _{\alpha },b\left\{ G\right\}
\varphi _{\alpha }\right) \right\vert \leq C\left\Vert b\right\Vert
_{L^{1}}\left\Vert G\right\Vert _{1}<+\infty
\end{equation*}%
for any orthonormal systems $\left\{ \psi _{\alpha }\right\} _{\alpha \in I}$%
, $\left\{ \varphi _{\alpha }\right\} _{\alpha \in I}$. We obtain that $%
b\left\{ G\right\} \in \mathcal{B}_{1}\left( \mathcal{H}\right) $ and 
\begin{equation*}
\left\Vert b\left\{ G\right\} \right\Vert _{1}\leq C\left\Vert b\right\Vert
_{L^{1}}\left\Vert G\right\Vert _{1}.
\end{equation*}%
See \cite[p.246-247, Theorems 3 and 4]{Birman}

$(\func{c})$ is a consequence of $(\func{a})$ and $(\func{b})$ since the $p$%
-Schatten classes interpolate like $L^{p}$-spaces: $\left\lfloor \mathcal{B}%
_{1}\left( \mathcal{H}\right) ,\mathcal{B}_{\infty }\left( \mathcal{H}%
\right) \right\rfloor _{\theta }=\left\lfloor \mathcal{B}_{1}\left( \mathcal{%
H}\right) ,\mathcal{B}\left( \mathcal{H}\right) \right\rfloor _{\theta }=%
\mathcal{B}_{\frac{1}{1-\theta }}\left( \mathcal{H}\right) ,$ $0<\theta <1$,
(see \cite[p.147]{Triebel}).

We now give a direct elementary proof of this part. Let $\left\{ \psi
_{\alpha }\right\} _{\alpha \in I}$ and $\left\{ \varphi _{\alpha }\right\}
_{\alpha \in I}$ be two orthonormal systems in $\mathcal{H}$. Then by H\"{o}%
lder inequality and part $(\func{a})$:%
\begin{multline*}
\sum_{\alpha \in I}\left\vert \left( \psi _{\alpha },b\left\{ G\right\}
\varphi _{\alpha }\right) \right\vert ^{p}\leq \sum_{\alpha \in I}\left(
\int_{Y}\left\vert b\left( y\right) \right\vert \left\vert \left( \psi
_{\alpha },U\left( y\right) ^{\ast }GU\left( y\right) \varphi _{\alpha
}\right) \right\vert \limfunc{d}y\right) ^{p} \\
\leq \sum_{\alpha \in I}\left( \int_{Y}\left\vert \left( \psi _{\alpha
},U\left( y\right) ^{\ast }GU\left( y\right) \varphi _{\alpha }\right)
\right\vert \limfunc{d}y\right) ^{\frac{p}{q}}\int_{Y}\left\vert b\left(
y\right) \right\vert ^{p}\left\vert \left( \psi _{\alpha },U\left( y\right)
^{\ast }GU\left( y\right) \varphi _{\alpha }\right) \right\vert \limfunc{d}y
\\
\leq \left( C\left\Vert G\right\Vert _{1}\right) ^{\frac{p}{q}%
}\int_{Y}\left\vert b\left( y\right) \right\vert ^{p}\sum_{\alpha \in
I}\left\vert \left( \psi _{\alpha },U\left( y\right) ^{\ast }GU\left(
y\right) \varphi _{\alpha }\right) \right\vert \limfunc{d}y \\
\leq \left( C\left\Vert G\right\Vert _{1}\right) ^{\frac{p}{q}%
}\int_{Y}\left\vert b\left( y\right) \right\vert ^{p}\left\Vert U\left(
y\right) ^{\ast }GU\left( y\right) \right\Vert _{1}\limfunc{d}y \\
\leq \left( C\left\Vert G\right\Vert _{1}\right) ^{\frac{p}{q}%
+1}\int_{Y}\left\vert b\left( y\right) \right\vert ^{p}\limfunc{d}y=\left(
C\left\Vert G\right\Vert _{1}\right) ^{\frac{p}{q}+1}\left\Vert b\right\Vert
_{L^{p}}^{p}<+\infty .
\end{multline*}%
Hence 
\begin{equation*}
\left( \sum_{\alpha \in I}\left\vert \left( \psi _{\alpha },b\left\{
G\right\} \varphi _{\alpha }\right) \right\vert ^{p}\right) ^{\frac{1}{p}%
}\leq C\left\Vert G\right\Vert _{1}\left\Vert b\right\Vert _{L^{p}}<+\infty
\end{equation*}%
for any orthonormal systems $\left\{ \psi _{\alpha }\right\} _{\alpha \in I}$%
, $\left\{ \varphi _{\alpha }\right\} _{\alpha \in I}$. Thus, by using
Proposition 2.6 of Simon \cite{Simon1}, we conclude that $b\left\{ G\right\}
\in \mathcal{B}_{p}\left( \mathcal{H}\right) $ and 
\begin{equation*}
\left\Vert b\left\{ G\right\} \right\Vert _{p}\leq C\left\Vert b\right\Vert
_{L^{p}}\left\Vert G\right\Vert _{1}.
\end{equation*}
\end{proof}

\begin{lemma}
Let $(\mathcal{H},\mathcal{W})$ be an irreducible Weyl system associated to
the symplectic space $\mathfrak{S}$. Then for any $\varphi ,\psi $ in $%
\mathcal{H}$ the map 
\begin{equation*}
\mathfrak{S}\ni \xi \rightarrow \left( \varphi ,\mathcal{W}\left( \xi
\right) \psi \right) \in 
%TCIMACRO{\U{2102} }%
%BeginExpansion
\mathbb{C}
%EndExpansion
,
\end{equation*}%
belongs to $L^{2}\left( \mathfrak{S}\right) \cap \mathcal{C}_{\mathcal{%
\infty }}(\mathfrak{S})$ and 
\begin{equation*}
\int_{\mathfrak{S}}\left\vert \left( \varphi ,\mathcal{W}\left( \xi \right)
\psi \right) \right\vert ^{2}\func{d}\xi =\left\Vert \varphi \right\Vert
^{2}\left\Vert \psi \right\Vert ^{2}\quad \text{and}\quad \left\Vert \left(
\varphi ,\mathcal{W}\left( \cdot \right) \psi \right) \right\Vert _{\infty
}\leq \left\Vert \varphi \right\Vert \left\Vert \psi \right\Vert .
\end{equation*}

For $\varphi ,\psi ,\varphi ^{\prime },\psi ^{\prime }\in \mathcal{H}$ we
have 
\begin{equation*}
\int_{\mathfrak{S}}\overline{\left( \varphi ^{\prime },\mathcal{W}\left( \xi
\right) \psi ^{\prime }\right) }\left( \varphi ,\mathcal{W}\left( \xi
\right) \psi \right) \func{d}\xi =\left( \varphi ,\varphi ^{\prime }\right)
\left( \psi ^{\prime },\psi \right) .
\end{equation*}
\end{lemma}

\begin{proof}
Since a symplectic space has only one irreducible representation (modulo
unitary equivalence), we may assume that $\mathfrak{S}=T^{\ast }\left(
X\right) $ and that $(\mathcal{H},\mathcal{W})$ is the Schr\"{o}dinger
representation $(\mathcal{H}\left( X\right) ,\mathcal{W})$. If $\varphi
,\psi \in \mathcal{H}\left( X\right) =L^{2}\left( X,\limfunc{d}x\right) $,
then%
\begin{eqnarray*}
\left( \varphi ,\mathcal{W}\left( x,p\right) \psi \right) &=&\int_{X}%
\overline{\varphi \left( y\right) }e^{i\langle y-x/2,p\rangle }\psi \left(
y-x\right) \func{d}y \\
&=&e^{-\frac{i}{2}\langle x,p\rangle }\overline{\mathcal{F}}_{X}\left( 
\overline{\varphi }T_{x}\psi \right) \left( p\right) ,
\end{eqnarray*}%
where $T_{x}\psi \left( \cdot \right) =\psi \left( \cdot -x\right) $, and 
\begin{eqnarray*}
\int_{\mathfrak{S}}\left\vert \left( \varphi ,\mathcal{W}\left( \xi \right)
\psi \right) \right\vert ^{2}\func{d}\xi &=&\iint_{X\times X^{\ast
}}\left\vert \left( \varphi ,\mathcal{W}\left( x,p\right) \psi \right)
\right\vert ^{2}\func{d}x\func{d}p \\
&=&\int_{X}\left\Vert \overline{\mathcal{F}}_{X}\left( \overline{\varphi }%
T_{x}\psi \right) \right\Vert _{L^{2}\left( X^{\ast }\right) }^{2}\func{d}%
x=\int_{X}\left\Vert \overline{\varphi }T_{x}\psi \right\Vert _{L^{2}\left(
X\right) }^{2}\func{d}x \\
&=&\int_{X}\left( \int_{X}\left\vert \overline{\varphi \left( y\right) }\psi
\left( y-x\right) \right\vert ^{2}\func{d}y\right) \func{d}x \\
&=&\int_{X}\left( \int_{X}\left\vert \overline{\varphi \left( y\right) }\psi
\left( y-x\right) \right\vert ^{2}\func{d}x\right) \func{d}y=\left\Vert
\varphi \right\Vert ^{2}\left\Vert \psi \right\Vert ^{2}.
\end{eqnarray*}%
Let $\left\{ \varphi _{n}\right\} ,\left\{ \psi _{n}\right\} \subset 
\mathcal{S}\left( X\right) $ such that $\varphi _{n}\rightarrow \varphi $
and $\psi _{n}\rightarrow \psi $ in $\mathcal{H}\left( X\right) $. Since $%
\left\{ \left( \varphi _{n},\mathcal{W}\left( \cdot \right) \psi _{n}\right)
\right\} \subset \mathcal{S}\left( \mathfrak{S}\right) $ and 
\begin{eqnarray*}
\left\Vert \left( \varphi _{n},\mathcal{W}\left( \cdot \right) \psi
_{n}\right) -\left( \varphi ,\mathcal{W}\left( \cdot \right) \psi \right)
\right\Vert _{\infty } &=&\left\Vert \left( \varphi _{n}-\varphi ,\mathcal{W}%
\left( \cdot \right) \psi _{n}\right) +\left( \varphi ,\mathcal{W}\left(
\cdot \right) \left( \psi _{n}-\psi \right) \right) \right\Vert \\
&\leq &\left\Vert \varphi _{n}-\varphi \right\Vert \left\Vert \psi
_{n}\right\Vert +\left\Vert \varphi \right\Vert \left\Vert \psi _{n}-\psi
\right\Vert ,\quad n\in 
%TCIMACRO{\U{2115} }%
%BeginExpansion
\mathbb{N}
%EndExpansion
,
\end{eqnarray*}%
it follows that $\left( \varphi _{n},\mathcal{W}\left( \cdot \right) \psi
_{n}\right) \rightarrow \left( \varphi ,\mathcal{W}\left( \cdot \right) \psi
\right) $ in $\mathcal{C}_{\mathcal{\infty }}(\mathfrak{S})$ as $%
n\rightarrow \infty $. Hence $\left( \varphi ,\mathcal{W}\left( \cdot
\right) \psi \right) $ $\in L^{2}\left( \mathfrak{S}\right) \cap \mathcal{C}%
_{\mathcal{\infty }}(\mathfrak{S})$, $\left\Vert \left( \varphi ,\mathcal{W}%
\left( \cdot \right) \psi \right) \right\Vert _{L^{2}\left( \mathfrak{S}%
\right) }^{2}=\left\Vert \varphi \right\Vert ^{2}\left\Vert \psi \right\Vert
^{2}$ and $\left\Vert \left( \varphi ,\mathcal{W}\left( \cdot \right) \psi
\right) \right\Vert _{\infty }\leq \left\Vert \varphi \right\Vert \left\Vert
\psi \right\Vert $. The last formula is a consequence of polarization
identity.
\end{proof}

\begin{theorem}
\label{KOC2}Let $(\mathcal{H},\mathcal{W})$ be an irreducible Weyl system
associated to the symplectic space $\mathfrak{S}$.

$(\func{a})$ If $b\in L^{\infty }\left( \mathfrak{S}\right) $ and $G\in 
\mathcal{B}_{1}\left( \mathcal{H}\right) $, then the integral%
\begin{equation*}
b\left\{ G\right\} =\int_{\mathfrak{S}}b\left( \xi \right) \mathcal{W}\left(
\xi \right) G\mathcal{W}\left( -\xi \right) \limfunc{d}\xi
\end{equation*}%
is weakly absolutely convergent and defines a bounded operator such that%
\begin{equation*}
\left\Vert b\left\{ G\right\} \right\Vert \leq \left\Vert b\right\Vert
_{L^{\infty }}\left\Vert G\right\Vert _{1}.
\end{equation*}%
Moreover, if $b$ vanishes at $\infty $ in the sense that for any $%
\varepsilon >0$ there is a compact subset $K$ of $\mathfrak{S}$ such that 
\begin{equation*}
\left\Vert b\right\Vert _{L^{\infty }\left( \mathfrak{S}\backslash K\right)
}\leq \varepsilon ,
\end{equation*}%
then $b\left\{ G\right\} $ is a compact operator.

The mapping $\left( b,G\right) \rightarrow b\left\{ G\right\} $ has the
following properties.

$(\func{i})$ $b\geq 0,G\geq 0\Rightarrow b\left\{ G\right\} \geq 0.$

$(\func{ii})$ $1\left\{ G\right\} =\limfunc{Tr}\left( G\right) \limfunc{id}_{%
\mathcal{H}}$.

$\left( \func{iii}\right) $ $\left( b_{1}b_{2}\right) \left\{ G\right\} \ll
\left( \left\vert b_{1}\right\vert ^{2}\left\{ \left\vert G^{\ast
}\right\vert \right\} ;\left\vert b_{2}\right\vert ^{2}\left\{ \left\vert
G\right\vert \right\} \right) $.

$(\func{b})$ If $b\in L^{p}\left( \mathfrak{S}\right) $ with $1\leq p<\infty 
$ and $G\in \mathcal{B}_{1}\left( \mathcal{H}\right) $, then the integral%
\begin{equation*}
b\left\{ G\right\} =\int_{\mathfrak{S}}b\left( \xi \right) \mathcal{W}\left(
\xi \right) G\mathcal{W}\left( -\xi \right) \limfunc{d}\xi
\end{equation*}%
is weakly absolutely convergent and defines an operator $b\left\{ G\right\} $
in $\mathcal{B}_{p}\left( \mathcal{H}\right) $ which satisfies%
\begin{equation*}
\left\Vert b\left\{ G\right\} \right\Vert _{p}\leq \left\Vert b\right\Vert
_{L^{p}}\left\Vert G\right\Vert _{1}.
\end{equation*}
\end{theorem}

\begin{proof}
$(\func{ii})$ Let $G=|\varphi )(\psi |=\left( \psi ,\cdot \right) \varphi ,$ 
$\varphi ,\psi \in \mathcal{H}$. Then 
\begin{equation*}
\mathcal{W}\left( \xi \right) G\mathcal{W}\left( -\xi \right) =|\mathcal{W}%
\left( \xi \right) \varphi )(\mathcal{W}\left( \xi \right) \psi |
\end{equation*}%
and 
\begin{equation*}
\int_{\mathfrak{S}}\left( u,\mathcal{W}\left( \xi \right) \varphi \right)
\left( \mathcal{W}\left( \xi \right) \psi ,v\right) \limfunc{d}\xi =\left(
\psi ,\varphi \right) \left( u,v\right) =\left( u,\limfunc{Tr}\left(
|\varphi )(\psi |\right) v\right) .
\end{equation*}%
So the equality holds for operators of rank $1$. Next we extend this
equality by linearty and continuity.

$\left( \func{iii}\right) $ We have $\mathcal{W}\left( \xi \right) G\mathcal{%
W}\left( -\xi \right) \ll \left( \mathcal{W}\left( \xi \right) \left\vert
G^{\ast }\right\vert \mathcal{W}\left( -\xi \right) ;\mathcal{W}\left( \xi
\right) \left\vert G\right\vert \mathcal{W}\left( -\xi \right) \right) $
which gives 
\begin{multline*}
\left\vert b_{1}\left( \xi \right) b_{2}\left( \xi \right) \left( \varphi ,%
\mathcal{W}\left( \xi \right) G\mathcal{W}\left( -\xi \right) \psi \right)
\right\vert \\
\leq \left( \left\vert b_{1}\left( \xi \right) \right\vert ^{2}\left(
\varphi ,\mathcal{W}\left( \xi \right) \left\vert G^{\ast }\right\vert 
\mathcal{W}\left( -\xi \right) \varphi \right) \right) ^{\frac{1}{2}}\left(
\left\vert b_{2}\left( \xi \right) \right\vert ^{2}\left( \psi ,\mathcal{W}%
\left( \xi \right) \left\vert G\right\vert \mathcal{W}\left( -\xi \right)
\psi \right) \right) ^{\frac{1}{2}}.
\end{multline*}%
Now we just use the Schwarz inequality to conclude that $\left( \func{iii}%
\right) $ is true.
\end{proof}

\section{Some special symbols}

To apply Theorem \ref{KOC2} in combination with Corollary \ref{K4}, we need
some special symbols $g$ for which $g_{X}^{\tau }\left( R\right) $ has an
extension $G\in \mathcal{B}_{1}\left( \mathcal{H}\right) $. Such symbols
have been constructed by Cordes \cite{Cordes}.

Let $\left( E,\left\vert \cdot \right\vert \right) $ be an euclidean space.
If $x\in E$,we set $\left\langle x\right\rangle =\left( 1+\left\vert
x\right\vert ^{2}\right) ^{1/2}$. Sometimes, in order to avoid confusions,
we shall add a subscript specifying the space, e.g. $\left( \cdot ,\cdot
\right) _{E}$, $\left\vert \cdot \right\vert _{E}$ or $\left\langle \cdot
\right\rangle _{E}$.

Let $\mathfrak{S}=T^{\ast }\left( X\right) $ with the standard symplectic
structure and $\left( X,\left\vert \cdot \right\vert _{X}\right) $ an
euclidean space. We shall work in the Schr\"{o}dinger representation $(%
\mathcal{H}\left( X\right) ,\mathcal{W})$.

Let $s>0$. Then $\left( 1-\triangle _{X}\right) ^{\frac{s}{2}}=\mathcal{F}%
_{X}^{-1}M_{\left\langle \cdot \right\rangle _{X^{\ast }}^{\frac{s}{2}}}%
\mathcal{F}_{X}=a_{X}^{0}\left( Q,P\right) $, where $a\left( x,p\right)
=\left\langle p\right\rangle _{X^{\ast }}^{\frac{s}{2}}$, $\left( x,p\right)
\in X\times X^{\ast }.$ Let $\psi _{s}=\psi _{s}^{X}$ be the unique solution
within $\mathcal{S}^{\ast }(X)$ for 
\begin{equation*}
\left( 1-\triangle _{X}\right) ^{\frac{s}{2}}\psi _{s}=\delta .
\end{equation*}%
Similarly, $\left( 1-\triangle _{X^{\ast }}\right) ^{\frac{s}{2}}=\mathcal{F}%
_{X^{\ast }}^{-1}M_{\left\langle \cdot \right\rangle _{X}^{\frac{s}{2}}}%
\mathcal{F}_{X^{\ast }}$. We shall denote by $\chi _{s}=\chi _{s}^{X^{\ast
}} $ the unique solution in $\mathcal{S}^{\ast }(X^{\ast })$ for 
\begin{equation*}
\left( 1-\triangle _{X^{\ast }}\right) ^{\frac{s}{2}}\chi _{s}=\delta .
\end{equation*}%
Let $n=\dim X$. We recall that $\psi _{s}\in L^{1}\left( X\right) \cap 
\mathcal{C}^{\infty }\left( X\backslash \left\{ 0\right\} \right) $, $\psi
_{s}\left( x\right) $ and its derivatives decay exponentially as $\left\vert
x\right\vert \rightarrow \infty $, and that $\partial ^{\alpha }\psi
_{s}\left( x\right) =O\left( 1+\left\vert x\right\vert _{X}^{s-n-\left\vert
\alpha \right\vert }\right) $ as $\left\vert x\right\vert \rightarrow 0$,
except when $s-n-\left\vert \alpha \right\vert =0$ in which case we have $%
\partial ^{\alpha }\psi _{s}\left( x\right) =O\left( 1+\log \frac{1}{%
\left\vert x\right\vert _{X}}\right) $. $\chi _{s}$ has similar properties.

\begin{lemma}[Cordes \protect\cite{Cordes}]
Let $t,s>\frac{n}{2}$ and 
\begin{equation*}
g:X\times X^{\ast }\rightarrow 
%TCIMACRO{\U{211d} }%
%BeginExpansion
\mathbb{R}
%EndExpansion
,\quad g\left( x,p\right) =\psi _{t}\left( x\right) \chi _{s}\left( p\right)
,\quad \left( x,p\right) \in X\times X^{\ast },
\end{equation*}%
i.e. $g=\psi _{t}\otimes \chi _{s}$. Then $g_{X}^{0}\left( Q,P\right) $ has
an extension $G\in \mathcal{B}_{1}\left( \mathcal{H}\left( X\right) \right) $%
. Moreover, for any $\alpha \in 
%TCIMACRO{\U{2115} }%
%BeginExpansion
\mathbb{N}
%EndExpansion
^{n}$, $g_{X}^{0}\left( Q,P\right) P_{\alpha }$ has an extension $G_{\alpha
}\in \mathcal{B}_{1}\left( \mathcal{H}\left( X\right) \right) $.
\end{lemma}

\begin{lemma}[Cordes \protect\cite{Cordes}]
Let $t,s>n$ and 
\begin{equation*}
g:X\times X^{\ast }\rightarrow 
%TCIMACRO{\U{211d} }%
%BeginExpansion
\mathbb{R}
%EndExpansion
,\quad g\left( x,p\right) =\psi _{t}\left( x\right) \chi _{s}\left( p\right)
,\quad \left( x,p\right) \in X\times X^{\ast },
\end{equation*}%
i.e. $g=\psi _{t}\otimes \chi _{s}$. If $\tau \in 
%TCIMACRO{\U{211d} }%
%BeginExpansion
\mathbb{R}
%EndExpansion
\equiv 
%TCIMACRO{\U{211d} }%
%BeginExpansion
\mathbb{R}
%EndExpansion
\cdot 1_{X}$, then $g_{X}^{\tau }\left( Q,P\right) $ has an extension in $%
\mathcal{B}_{1}\left( \mathcal{H}\left( X\right) \right) $ denoted also by $%
g_{X}^{\tau }\left( Q,P\right) $. The mapping%
\begin{equation*}
%TCIMACRO{\U{211d} }%
%BeginExpansion
\mathbb{R}
%EndExpansion
\ni \tau \rightarrow g_{X}^{\tau }\left( Q,P\right) \in \mathcal{B}%
_{1}\left( \mathcal{H}\left( X\right) \right)
\end{equation*}%
is continuous.
\end{lemma}

For the proof of these two lemmas see \cite{Cordes}.

Let $X=X_{1}\oplus ...\oplus X_{k}$ be an orthogonal decomposition with the
canonical injections $j_{1}:X_{1}\hookrightarrow X$,...,$j_{k}:X_{k}%
\hookrightarrow X$ and $\pi _{1}:X\rightarrow X_{1}$,...,$\pi
_{k}:X\rightarrow X_{k}$ the orthogonal projections. Then $X^{\ast
}=X_{1}^{\ast }\oplus ...\oplus X_{k}^{\ast }$ is an orthogonal
decomposition with $\pi _{1}^{\ast }:X_{1}^{\ast }\hookrightarrow X^{\ast }$%
,...,$\pi _{k}^{\ast }:X_{k}^{\ast }\hookrightarrow X^{\ast }$ as the
canonical injections and $j_{1}^{\ast }:X^{\ast }\rightarrow X_{1}^{\ast }$%
,...,$j_{k}^{\ast }:X^{\ast }\rightarrow X_{k}^{\ast }$ as the orthogonal
projections.

The mapping%
\begin{gather*}
J:X\times X^{\ast }\rightarrow \left( X_{1}\times X_{1}^{\ast }\right)
\times ...\times \left( X_{k}\times X_{k}^{\ast }\right) , \\
J\left( x,p\right) =\left( \left( \pi _{1}x,j_{1}^{\ast }p\right)
,...,\left( \pi _{k}x,j_{k}^{\ast }p\right) \right)
\end{gather*}%
is an isometry with the inverse given by 
\begin{gather*}
J^{-1}:\left( X_{1}\times X_{1}^{\ast }\right) \times ...\times \left(
X_{k}\times X_{k}^{\ast }\right) \rightarrow X\times X^{\ast }, \\
J^{-1}\left( \left( x_{1},p_{1}\right) ,...,\left( x_{k},p_{k}\right)
\right) =\left( x_{1}+...+x_{k},p_{1}+...+p_{k}\right) .
\end{gather*}

If $a_{1}\in \mathcal{S}^{\ast }\left( X_{1}\times X_{1}^{\ast }\right) $%
,..., $a_{k}\in \mathcal{S}^{\ast }\left( X_{k}\times X_{k}^{\ast }\right) $%
, then 
\begin{equation*}
a=\left( a_{1}\otimes ...\otimes a_{k}\right) \circ J\in \mathcal{S}^{\ast
}\left( X\times X^{\ast }\right)
\end{equation*}%
and 
\begin{equation*}
a_{X}^{\tau }\left( Q,P\right) =\left( a_{1}\right) _{X_{1}}^{\tau }\left(
Q_{1},P_{1}\right) \otimes ...\otimes \left( a_{k}\right) _{X_{k}}^{\tau
}\left( Q_{k},P_{k}\right) ,
\end{equation*}%
where $Q=\left( Q_{1},...,Q_{k}\right) $ and $P=\left(
P_{1},...,P_{k}\right) $.

Let $s_{1},...,s_{k}>0$ and $\mathbf{s}=\left( s_{1},...,s_{k}\right) $.
Then $\psi _{\mathbf{s}}=\psi _{s_{1}}^{X_{1}}\otimes ...\otimes \psi
_{s_{k}}^{X_{k}}\in \mathcal{S}^{\ast }(X)$ is the unique solution within $%
\mathcal{S}^{\ast }(X)$ for 
\begin{equation*}
\left( 1-\triangle _{X_{1}}\right) ^{\frac{s_{1}}{2}}\otimes ...\otimes
\left( 1-\triangle _{X_{k}}\right) ^{\frac{s_{k}}{2}}\psi =\delta .
\end{equation*}%
Similarly, $\chi _{\mathbf{s}}=\chi _{s_{1}}^{X_{1}^{\ast }}\otimes
...\otimes \chi _{s_{k}}^{X_{k}^{\ast }}\in \mathcal{S}^{\ast }(X^{\ast })$
is the unique solution within $\mathcal{S}^{\ast }(X^{\ast })$ for%
\begin{equation*}
\left( 1-\triangle _{X_{1}^{\ast }}\right) ^{\frac{s_{1}}{2}}\otimes
...\otimes \left( 1-\triangle _{X_{k}^{\ast }}\right) ^{\frac{s_{k}}{2}}\chi
=\delta .
\end{equation*}

For $\mathbf{t}=\left( t_{1},...,t_{k}\right) $, $t_{1},...,t_{k}>0$ and $%
\mathbf{s}=\left( s_{1},...,s_{k}\right) $, $s_{1},...,s_{k}>0$ we introduce
the distributions $g=g_{\mathbf{t,s}}=\psi _{\mathbf{t}}\otimes \chi _{%
\mathbf{s}}\in \mathcal{S}^{\ast }\left( X\times X^{\ast }\right) $ and $%
g_{j}=g_{t_{j}\mathbf{,}s_{j}}=\psi _{t_{j}}\otimes \chi _{s_{j}}\in 
\mathcal{S}^{\ast }\left( X_{j}\times X_{j}^{\ast }\right) $, for $j=1,...,k$%
. Then $g=\left( g_{1}\otimes ...\otimes g_{k}\right) \circ J$ and 
\begin{equation*}
g_{X}^{\tau }\left( Q,P\right) =\left( g_{1}\right) _{X_{1}}^{\tau }\left(
Q_{1},P_{1}\right) \otimes ...\otimes \left( g_{k}\right) _{X_{k}}^{\tau
}\left( Q_{k},P_{k}\right) .
\end{equation*}

We recall that the mapping 
\begin{eqnarray*}
\mathcal{B}_{1}\left( \mathcal{H}\right) \times \mathcal{B}_{1}\left( 
\mathcal{H}\right) &\rightarrow &\mathcal{B}_{1}\left( \mathcal{H}\otimes 
\mathcal{H}\right) , \\
\left( A,B\right) &\rightarrow &A\otimes B,
\end{eqnarray*}%
is well defined and 
\begin{equation*}
\left\Vert A\otimes B\right\Vert _{1}=\left\Vert A\right\Vert _{1}\left\Vert
B\right\Vert _{1},\quad \left( A,B\right) \in \mathcal{B}_{1}\left( \mathcal{%
H}\right) \times \mathcal{B}_{1}\left( \mathcal{H}\right) .
\end{equation*}

\begin{corollary}
\label{sss1}$(\func{a})$ Let $\mathbf{t}=\left( t_{1},...,t_{k}\right) ,$ $%
\mathbf{s}=\left( s_{1},...,s_{k}\right) $ such that $t_{1},s_{1}>\frac{\dim
X_{1}}{4},...$, $t_{k},s_{k}>\frac{\dim X_{k}}{4}$. Let $g=g_{2\mathbf{t},2%
\mathbf{s}}=\psi _{2\mathbf{t}}\otimes \chi _{2\mathbf{s}}\in \mathcal{S}%
^{\ast }\left( X\times X^{\ast }\right) $. Then $g_{X}^{0}\left( Q,P\right) $
has an extension $G\in \mathcal{B}_{1}\left( \mathcal{H}\left( X\right)
\right) $. Moreover, for any $\alpha \in 
%TCIMACRO{\U{2115} }%
%BeginExpansion
\mathbb{N}
%EndExpansion
^{n}$, $g_{X}^{0}\left( Q,P\right) P_{\alpha }$ has an extension $G_{\alpha
}\in \mathcal{B}_{1}\left( \mathcal{H}\left( X\right) \right) $.

$(\func{b})$ Let $\mathbf{t}=\left( t_{1},...,t_{k}\right) ,$ $\mathbf{s}%
=\left( s_{1},...,s_{k}\right) $ such that $t_{1},s_{1}>\frac{\dim X_{1}}{4}%
,...,t_{k},s_{k}>\frac{\dim X_{k}}{4}$. Let $g=g_{4\mathbf{t},4\mathbf{s}%
}=\psi _{4\mathbf{t}}\otimes \chi _{4\mathbf{s}}\in \mathcal{S}^{\ast
}\left( X\times X^{\ast }\right) $. If $\tau \in 
%TCIMACRO{\U{211d} }%
%BeginExpansion
\mathbb{R}
%EndExpansion
\equiv 
%TCIMACRO{\U{211d} }%
%BeginExpansion
\mathbb{R}
%EndExpansion
\cdot 1_{X}$, then $g_{X}^{\tau }\left( Q,P\right) $ has an extension in $%
\mathcal{B}_{1}\left( \mathcal{H}\left( X\right) \right) $ denoted also by $%
g_{X}^{\tau }\left( Q,P\right) $. The mapping%
\begin{equation*}
%TCIMACRO{\U{211d} }%
%BeginExpansion
\mathbb{R}
%EndExpansion
\ni \tau \rightarrow g_{X}^{\tau }\left( Q,P\right) \in \mathcal{B}%
_{1}\left( \mathcal{H}\left( X\right) \right)
\end{equation*}%
is continuous.
\end{corollary}

\section{Schatten-class properties of pseudo-differential operators}

We are now able to consider $L^{2}$-boundedness and Schatten-class
properties of certain pseudo-differential operators. In the notation of the
previous section, $\mathfrak{S}$ will be $T^{\ast }\left( X\right) $ with
the standard symplectic structure and $\left( X,\left\vert \cdot \right\vert
_{X}\right) $ an euclidean space, $(\mathcal{H}\left( X\right) ,\mathcal{W})$
will be the Schr\"{o}dinger representation associated to the symplectic
space $\mathfrak{S}$. Let $X=X_{1}\oplus ...\oplus X_{k}$ be an orthogonal
decomposition and $X^{\ast }=X_{1}^{\ast }\oplus ...\oplus X_{k}^{\ast }$ be
the dual orthogonal decomposition.

\begin{theorem}
\label{tcp2}Let $a\in \mathcal{S}^{\ast }(\mathfrak{S})$ and $1\leq p<\infty 
$.

$(\func{a})$ Assume that there are $\mathbf{t}=\left( t_{1},...,t_{k}\right) 
$, $\mathbf{s}=\left( s_{1},...,s_{k}\right) $ such that $t_{1},s_{1}>\frac{%
\dim X_{1}}{4},...,t_{k},s_{k}>\frac{\dim X_{k}}{4}$ and 
\begin{equation*}
b=\left( 1-\triangle _{X_{1}}\right) ^{t_{1}}\otimes ...\otimes \left(
1-\triangle _{X_{k}}\right) ^{t_{k}}\otimes \left( 1-\triangle _{X_{1}^{\ast
}}\right) ^{s_{1}}\otimes ...\otimes \left( 1-\triangle _{X_{k}^{\ast
}}\right) ^{s_{k}}a\in L^{p}\left( \mathfrak{S}\right) .
\end{equation*}%
Then $a_{X}^{0}\left( Q,P\right) $ has an extension in $\mathcal{B}%
_{p}\left( \mathcal{H}\left( X\right) \right) $ denoted also by $%
a_{X}^{0}\left( Q,P\right) $ and there is $C>0$ such that 
\begin{equation*}
\left\Vert a_{X}^{0}\left( Q,P\right) \right\Vert _{p}\leq C\left\Vert
b\right\Vert _{L^{p}\left( \mathfrak{S}\right) }.
\end{equation*}

$(\func{b})$ Assume that there are $\mathbf{t}=\left( t_{1},...,t_{k}\right) 
$, $\mathbf{s}=\left( s_{1},...,s_{k}\right) $ such that $t_{1},s_{1}>\frac{%
\dim X_{1}}{4},...,t_{k},s_{k}>\frac{\dim X_{k}}{4}$ and 
\begin{equation*}
c=\left( 1-\triangle _{X_{1}}\right) ^{2t_{1}}\otimes ...\otimes \left(
1-\triangle _{X_{k}}\right) ^{2t_{k}}\otimes \left( 1-\triangle
_{X_{1}^{\ast }}\right) ^{2s_{1}}\otimes ...\otimes \left( 1-\triangle
_{X_{k}^{\ast }}\right) ^{2s_{k}}a\in L^{p}\left( \mathfrak{S}\right) .
\end{equation*}%
If $\tau \in 
%TCIMACRO{\U{211d} }%
%BeginExpansion
\mathbb{R}
%EndExpansion
\equiv 
%TCIMACRO{\U{211d} }%
%BeginExpansion
\mathbb{R}
%EndExpansion
\cdot 1_{X}$, then $a_{X}^{\tau }\left( Q,P\right) $ has an extension in $%
\mathcal{B}_{p}\left( \mathcal{H}\left( X\right) \right) $ denoted also by $%
a_{X}^{\tau }\left( Q,P\right) $. The mapping%
\begin{equation*}
%TCIMACRO{\U{211d} }%
%BeginExpansion
\mathbb{R}
%EndExpansion
\ni \tau \rightarrow a_{X}^{\tau }\left( Q,P\right) \in \mathcal{B}%
_{p}\left( \mathcal{H}\left( X\right) \right)
\end{equation*}%
is continuous and for any $K$ a compact subset of $%
%TCIMACRO{\U{211d} }%
%BeginExpansion
\mathbb{R}
%EndExpansion
$, there is $C_{K}>0$ such that 
\begin{equation*}
\left\Vert a_{X}^{\tau }\left( Q,P\right) \right\Vert _{p}\leq
C_{K}\left\Vert c\right\Vert _{L^{p}\left( \mathfrak{S}\right) },
\end{equation*}%
for any $\tau \in K$.
\end{theorem}

\begin{proof}
$(\func{a})$ For $\mathbf{t}=\left( t_{1},...,t_{k}\right) $ and $\mathbf{s}%
=\left( s_{1},...,s_{k}\right) $ we consider $g=g_{2\mathbf{t,}2\mathbf{s}%
}=\psi _{2\mathbf{t}}\otimes \chi _{2\mathbf{s}}\in \mathcal{S}^{\ast
}\left( X\times X^{\ast }\right) $. Recall that we use the notations of the
previous section. Then $g\in L^{1}\left( X\times X^{\ast }\right) $ and $%
a=b\ast g$ because%
\begin{equation*}
\left( 1-\triangle _{X_{1}}\right) ^{t_{1}}\otimes ...\otimes \left(
1-\triangle _{X_{k}}\right) ^{t_{k}}\otimes \left( 1-\triangle _{X_{1}^{\ast
}}\right) ^{s_{1}}\otimes ...\otimes \left( 1-\triangle _{X_{k}^{\ast
}}\right) ^{s_{k}}g=\delta .
\end{equation*}%
It follows from Corollary \ref{sss1} $(\func{a})$ that $g_{X}^{0}\left(
Q,P\right) \subset G\in \mathcal{B}_{1}\left( \mathcal{H}\left( X\right)
\right) $. On the other hand, by using Corollary \ref{K4} and Theorem \ref%
{KOC2} $(\func{b})$ we coclude that 
\begin{equation*}
a_{X}^{0}\left( Q,P\right) =\left( b\ast g\right) _{X}^{0}\left( Q,P\right)
\subset b\left\{ G\right\} \in \mathcal{B}_{p}\left( \mathcal{H}\left(
X\right) \right)
\end{equation*}%
and we have 
\begin{equation*}
\left\Vert a_{X}^{0}\left( Q,P\right) \right\Vert _{p}\leq \left\Vert
b\right\Vert _{L^{p}\left( \mathfrak{S}\right) }\left\Vert g_{X}^{0}\left(
Q,P\right) \right\Vert _{1}.
\end{equation*}

$(\func{b})$ The proof of this point is essentially the same. For $\mathbf{t}%
=\left( t_{1},...,t_{k}\right) $ and $\mathbf{s}=\left(
s_{1},...,s_{k}\right) $ we set $g=g_{4\mathbf{t},4\mathbf{s}}=\psi _{4%
\mathbf{t}}\otimes \chi _{4\mathbf{s}}\in \mathcal{S}^{\ast }\left( X\times
X^{\ast }\right) $. Then $g\in L^{1}\left( X\times X^{\ast }\right) $ and $%
a=c\ast g$. Now we apply Corollary \ref{sss1} to obtain that $g_{X}^{\tau
}\left( Q,P\right) $ has an extension in $\mathcal{B}_{1}\left( \mathcal{H}%
\left( X\right) \right) $ and that the mapping%
\begin{equation*}
%TCIMACRO{\U{211d} }%
%BeginExpansion
\mathbb{R}
%EndExpansion
\ni \tau \rightarrow g_{X}^{\tau }\left( Q,P\right) \in \mathcal{B}%
_{1}\left( \mathcal{H}\left( X\right) \right)
\end{equation*}%
is continuous. From Corollary \ref{K4} and Theorem \ref{KOC2} $(\func{b})$
it follows that 
\begin{equation*}
a_{X}^{\tau }\left( Q,P\right) =\left( c\ast g\right) _{X}^{\tau }\left(
Q,P\right) \subset c\left\{ g_{X}^{\tau }\left( Q,P\right) \right\} \in 
\mathcal{B}_{p}\left( \mathcal{H}\left( X\right) \right)
\end{equation*}%
and we have 
\begin{eqnarray*}
\left\Vert a_{X}^{\tau }\left( Q,P\right) \right\Vert _{p} &\leq &\left\Vert
c\right\Vert _{L^{p}\left( \mathfrak{S}\right) }\left\Vert g_{X}^{\tau
}\left( Q,P\right) \right\Vert _{1}, \\
\left\Vert a_{X}^{\tau }\left( Q,P\right) -a_{X}^{\tau ^{\prime }}\left(
Q,P\right) \right\Vert _{p} &\leq &\left\Vert c\right\Vert _{L^{p}\left( 
\mathfrak{S}\right) }\left\Vert g_{X}^{\tau }\left( Q,P\right) -g_{X}^{\tau
^{\prime }}\left( Q,P\right) \right\Vert _{1},
\end{eqnarray*}%
for any $\tau ,$ $\tau ^{\prime }\in 
%TCIMACRO{\U{211d} }%
%BeginExpansion
\mathbb{R}
%EndExpansion
$. Hence the mapping%
\begin{equation*}
%TCIMACRO{\U{211d} }%
%BeginExpansion
\mathbb{R}
%EndExpansion
\ni \tau \rightarrow a_{X}^{\tau }\left( Q,P\right) \in \mathcal{B}%
_{p}\left( \mathcal{H}\left( X\right) \right)
\end{equation*}%
is continuous and for any $K$ a compact sub set of $%
%TCIMACRO{\U{211d} }%
%BeginExpansion
\mathbb{R}
%EndExpansion
$, there is $C_{K}>0$ such that 
\begin{equation*}
\left\Vert a_{X}^{\tau }\left( Q,P\right) \right\Vert _{p}\leq \sup \left\{
\left\Vert g_{X}^{\tau }\left( Q,P\right) \right\Vert _{1}:\tau \in
K\right\} \left\Vert c\right\Vert _{L^{p}\left( \mathfrak{S}\right) },
\end{equation*}%
for any $\tau \in K$.
\end{proof}

If we replace the $L^{p}$-conditions by $L^{\infty }$-conditions, then we
obtain the theorems due to Cordes.

\begin{theorem}
\label{tcp3}Let $a\in \mathcal{S}^{\ast }(\mathfrak{S})$.

$(\func{a})$ Assume that there are $\mathbf{t}=\left( t_{1},...,t_{k}\right) 
$, $\mathbf{s}=\left( s_{1},...,s_{k}\right) $ such that $t_{1},s_{1}>\frac{%
\dim X_{1}}{4},...,t_{k},s_{k}>\frac{\dim X_{k}}{4}$ and 
\begin{equation*}
b=\left( 1-\triangle _{X_{1}}\right) ^{t_{1}}\otimes ...\otimes \left(
1-\triangle _{X_{k}}\right) ^{t_{k}}\otimes \left( 1-\triangle _{X_{1}^{\ast
}}\right) ^{s_{1}}\otimes ...\otimes \left( 1-\triangle _{X_{k}^{\ast
}}\right) ^{s_{k}}a\in L^{\infty }\left( \mathfrak{S}\right) .
\end{equation*}%
Then $a_{X}^{0}\left( Q,P\right) $ has an extension in $\mathcal{B}\left( 
\mathcal{H}\left( X\right) \right) $ denoted also by $a_{X}^{0}\left(
Q,P\right) $ and there is $C>0$ such that 
\begin{equation*}
\left\Vert a_{X}^{0}\left( Q,P\right) \right\Vert _{\mathcal{B}\left( 
\mathcal{H}\left( X\right) \right) }\leq C\left\Vert b\right\Vert
_{L^{\infty }\left( \mathfrak{S}\right) }.
\end{equation*}

$(\func{b})$ Assume that there are $\mathbf{t}=\left( t_{1},...,t_{k}\right) 
$, $\mathbf{s}=\left( s_{1},...,s_{k}\right) $ such that $t_{1},s_{1}>\frac{%
\dim X_{1}}{4},...,t_{k},s_{k}>\frac{\dim X_{k}}{4}$ and 
\begin{equation*}
c=\left( 1-\triangle _{X_{1}}\right) ^{2t_{1}}\otimes ...\otimes \left(
1-\triangle _{X_{k}}\right) ^{2t_{k}}\otimes \left( 1-\triangle
_{X_{1}^{\ast }}\right) ^{2s_{1}}\otimes ...\otimes \left( 1-\triangle
_{X_{k}^{\ast }}\right) ^{2s_{k}}a\in L^{\infty }\left( \mathfrak{S}\right) .
\end{equation*}%
If $\tau \in 
%TCIMACRO{\U{211d} }%
%BeginExpansion
\mathbb{R}
%EndExpansion
\equiv 
%TCIMACRO{\U{211d} }%
%BeginExpansion
\mathbb{R}
%EndExpansion
\cdot 1_{X}$, then $a_{X}^{\tau }\left( Q,P\right) $ has an extension in $%
\mathcal{B}\left( \mathcal{H}\left( X\right) \right) $ denoted also by $%
a_{X}^{\tau }\left( Q,P\right) $. The mapping%
\begin{equation*}
%TCIMACRO{\U{211d} }%
%BeginExpansion
\mathbb{R}
%EndExpansion
\ni \tau \rightarrow a_{X}^{\tau }\left( Q,P\right) \in \mathcal{B}\left( 
\mathcal{H}\left( X\right) \right)
\end{equation*}%
is continuous and for any $K$ a compact subset of $%
%TCIMACRO{\U{211d} }%
%BeginExpansion
\mathbb{R}
%EndExpansion
$, there is $C_{K}>0$ such that 
\begin{equation*}
\left\Vert a_{X}^{\tau }\left( Q,P\right) \right\Vert _{\mathcal{B}\left( 
\mathcal{H}\left( X\right) \right) }\leq C_{K}\left\Vert c\right\Vert
_{L^{\infty }\left( \mathfrak{S}\right) },
\end{equation*}%
for any $\tau \in K$.
\end{theorem}

The proof of this theorem is essentially the same, the only change we must
do is the reference to part $(\func{a})$ insted of part $(\func{b})$ in
Theorem \ref{KOC2}.

Let $\left( X,\left\vert \cdot \right\vert _{X}\right) $ be an euclidean
space in which we choose an orthonormal basis. Let $m=\left[ \frac{\dim X}{2}%
\right] +1$ and $1\leq p\leq \infty $. Then there are $\tau =\tau \left(
\dim X\right) >\frac{\dim X}{4}$ and $\gamma >0$ such that 
\begin{equation*}
\left\Vert \left( 1-\triangle _{X}\right) ^{\tau }\varphi \right\Vert
_{L^{p}}\leq \gamma \left\Vert \varphi \right\Vert _{p,m},\quad \varphi \in 
\mathcal{S}(X),
\end{equation*}%
where $\Vert \cdot \parallel _{p,m}$ is the norm on $\mathcal{S}(X)$ defined
by 
\begin{equation*}
\left\Vert \varphi \right\Vert _{p,m}=\max \left\{ \left\Vert P_{\beta
}\varphi \right\Vert _{L^{p}}:\left\vert \beta \right\vert \leq m\right\}
=\max \left\{ \left\Vert \partial ^{\beta }\varphi \right\Vert
_{L^{p}}:\left\vert \beta \right\vert \leq m\right\} ,\quad \varphi \in 
\mathcal{S}(X).
\end{equation*}%
For the proof see \cite[p.118-119]{Cordes}.

Let $X=X_{1}\oplus ...\oplus X_{k}$ be an orthogonal decomposition, $m_{1}=%
\left[ \frac{\dim X_{1}}{2}\right] +1,...,m_{k}=\left[ \frac{\dim X_{k}}{2}%
\right] +1$ and $1\leq p\leq \infty $. By an induction argument we conclude
that there are $\tau _{1}>\frac{\dim X_{1}}{4},...,\tau _{k}>\frac{\dim X_{k}%
}{4}$ and $\gamma >0$ such that%
\begin{equation}
\left\Vert \left( 1-\triangle _{X_{1}}\right) ^{\tau _{1}}\otimes ...\otimes
\left( 1-\triangle _{X_{k}}\right) ^{\tau _{k}}\varphi \right\Vert
_{L^{p}}\leq \gamma \Vert \varphi \parallel _{p,m_{1},...,m_{k}},\quad
\varphi \in \mathcal{S}(X),  \label{tcp1}
\end{equation}%
where $\Vert \cdot \parallel _{p,m_{1},...,m_{k}}$ is the norm on $\mathcal{S%
}(X)$ defined by%
\begin{align*}
\Vert \varphi & \parallel _{p,m_{1},...,m_{k}}=\max \left\{ \left\Vert
P_{\beta _{1}}^{X_{1}}...P_{\beta _{k}}^{X_{k}}\varphi \right\Vert
_{L^{p}}:\left\vert \beta _{1}\right\vert \leq m_{1},...,\left\vert \beta
_{k}\right\vert \leq m_{k}\right\} \\
& =\max \left\{ \left\Vert \left( \partial ^{X_{1}}\right) ^{\beta
_{1}}...\left( \partial ^{X_{k}}\right) ^{\beta _{k}}\varphi \right\Vert
_{L^{p}}:\left\vert \beta _{1}\right\vert \leq m_{1},...,\left\vert \beta
_{k}\right\vert \leq m_{k}\right\} ,\quad \varphi \in \mathcal{S}(X).
\end{align*}%
We introduce the space 
\begin{multline*}
\mathcal{C}_{m_{1},...,m_{k}}^{p}=\left\{ u\in \mathcal{S}^{\ast
}(X):P_{\beta _{1}}^{X_{1}}...P_{\beta _{k}}^{X_{k}}u\in L^{p}\left(
X\right) ,\left\vert \beta _{1}\right\vert \leq m_{1},...,\left\vert \beta
_{k}\right\vert \leq m_{k}\right\} \\
=\left\{ u\in \mathcal{S}^{\ast }(X):\left( \partial ^{X_{1}}\right) ^{\beta
_{1}}...\left( \partial ^{X_{k}}\right) ^{\beta _{k}}u\in L^{p}\left(
X\right) ,\left\vert \beta _{1}\right\vert \leq m_{1},...,\left\vert \beta
_{k}\right\vert \leq m_{k}\right\} ,
\end{multline*}%
on which the norm $\Vert \cdot \parallel _{p,m_{1},...,m_{k}}$ has a natural
extension. With this norn $\mathcal{C}_{m_{1},...,m_{k}}^{p}$ becomes a
Banach space.

\begin{lemma}
\label{tcp4}There are $\tau _{1}>\frac{\dim X_{1}}{4},...,\tau _{k}>\frac{%
\dim X_{k}}{4}$ and $\gamma >0$ such that if $u\in \mathcal{C}%
_{m_{1},...,m_{k}}^{p}$, then $\left( 1-\triangle _{X_{1}}\right) ^{\tau
_{1}}\otimes ...\otimes \left( 1-\triangle _{X_{k}}\right) ^{\tau _{k}}u\in
L^{p}\left( X\right) $ and 
\begin{equation*}
\left\Vert \left( 1-\triangle _{X_{1}}\right) ^{\tau _{1}}\otimes ...\otimes
\left( 1-\triangle _{X_{k}}\right) ^{\tau _{k}}u\right\Vert _{L^{p}}\leq
\gamma \Vert u\parallel _{p,m_{1},...,m_{k}},\quad u\in \mathcal{C}%
_{m_{1},...,m_{k}}^{p}.
\end{equation*}
\end{lemma}

\begin{proof}
The constants $\tau _{1}>\frac{\dim X_{1}}{4},...,\tau _{k}>\frac{\dim X_{k}%
}{4}$ and $\gamma >0$ will be those in (\ref{tcp1}). We denote by $L$ the
operator $\left( 1-\triangle _{X_{1}}\right) ^{\tau _{1}}\otimes ...\otimes
\left( 1-\triangle _{X_{k}}\right) ^{\tau _{k}}$. Let $\left\{ \varphi
_{\varepsilon }\right\} $ be a family of smooth functions such that $0\leq
\varphi _{\varepsilon }\in \mathcal{C}_{0}^{\infty }\left( X\right) $, $\int
\varphi _{\varepsilon }\left( x\right) \limfunc{d}x=1$ and $\limfunc{supp}%
\varphi _{\varepsilon }\subset B\left( 0;\varepsilon \right) $. Let $\chi
\in \mathcal{C}_{0}^{\infty }\left( X\right) $ be such that $0\leq \chi \leq
1$, $\limfunc{supp}\chi \subset B\left( 0;2\right) $ and $\chi |B\left(
0;1\right) =1$. For $j\in 
%TCIMACRO{\U{2115} }%
%BeginExpansion
\mathbb{N}
%EndExpansion
$, $j\geq 1$, we set $\chi _{j}=\chi \left( \cdot /j\right) $.

Let $u\in \mathcal{C}_{m_{1},...,m_{k}}^{p}$. Then $u_{j}=\left( \chi
_{j}u\right) \ast \varphi _{1/j}\in \mathcal{C}_{0}^{\infty }\left( X\right) 
$ and $u_{j}\rightarrow u$ in $\mathcal{S}^{\ast }\left( X\right) $.

If $1\leq p<\infty $, then $u_{j}\rightarrow u$ in $\mathcal{C}%
_{m_{1},...,m_{k}}^{p}$, because 
\begin{equation*}
\left\Vert \left( \chi _{j}u\right) \ast \varphi _{1/j}-u\right\Vert
_{p,m_{1},...,m_{k}}\leq \left\Vert \left( \chi _{j}u\right) -u\right\Vert
_{p,m_{1},...,m_{k}}+\left\Vert u\ast \varphi _{1/j}-u\right\Vert
_{p,m_{1},...,m_{k}}.
\end{equation*}
Since 
\begin{equation*}
\left\Vert Lu_{j}\right\Vert _{L^{p}}\leq \gamma \Vert u_{j}\parallel
_{p,m_{1},...,m_{k}},\quad \left\Vert Lu_{j}-Lu_{j^{\prime }}\right\Vert
_{L^{p}}\leq \gamma \Vert u_{j}-u_{j^{\prime }}\parallel
_{p,m_{1},...,m_{k}},
\end{equation*}%
it follows that $Lu_{j}\rightarrow Lu$ in $L^{p}\left( X\right) $ and $%
\left\Vert Lu\right\Vert _{L^{p}}\leq \gamma \Vert u\parallel
_{p,m_{1},...,m_{k}}$.

If $p=\infty $, then $u_{j}\rightarrow u$ only in $\mathcal{S}^{\ast }\left(
X\right) $ and 
\begin{eqnarray*}
\left\vert \left\langle \psi ,Lu_{j}\right\rangle \right\vert &\leq
&\left\Vert \psi \right\Vert _{L^{1}}\left\Vert Lu_{j}\right\Vert
_{L^{\infty }}\leq \gamma \left\Vert \psi \right\Vert _{L^{1}}\Vert
u_{j}\parallel _{\infty ,m_{1},...,m_{k}} \\
&\leq &\gamma \left\Vert \psi \right\Vert _{L^{1}}\Vert \chi _{j}u\parallel
_{\infty ,m_{1},...,m_{k}}\leq \gamma \left( 1+C\left( \chi \right)
j^{-1}\right) \left\Vert \psi \right\Vert _{L^{1}}\Vert u\parallel _{\infty
,m_{1},...,m_{k}},
\end{eqnarray*}%
for all $\psi \in \mathcal{S}\left( X\right) $ and $j\in 
%TCIMACRO{\U{2115} }%
%BeginExpansion
\mathbb{N}
%EndExpansion
$, $j\geq 1$. If we let $j\rightarrow \infty $ we obtain that%
\begin{equation*}
\left\vert \left\langle \psi ,Lu\right\rangle \right\vert \leq \gamma
\left\Vert \psi \right\Vert _{L^{1}}\Vert u\parallel _{\infty
,m_{1},...,m_{k}},\quad \psi \in \mathcal{S}\left( X\right) .
\end{equation*}%
It follows that $Lu\in \left( L^{1}\left( X\right) \right) ^{\ast
}=L^{\infty }\left( X\right) $ and $\left\Vert Lu\right\Vert _{L^{\infty
}}\leq \gamma \Vert u\parallel _{\infty ,m_{1},...,m_{k}}$.
\end{proof}

We return to the symplectic space $\mathfrak{S}=T^{\ast }\left( X\right) $.
An orthogonal decomposition of $X$, $X=X_{1}\oplus ...\oplus X_{k}$, gives
an orthogonal decomposition of $\mathfrak{S}$, $\mathfrak{S}=X_{1}\oplus
...\oplus X_{k}\oplus X_{1}^{\ast }\oplus ...\oplus X_{k}^{\ast }$, if on $%
\mathfrak{S}$ we consider the euclidean norm $\left\Vert \left( x,p\right)
\right\Vert _{\mathfrak{S}}^{2}=\left\Vert x\right\Vert _{X}^{2}+\left\Vert
p\right\Vert _{X^{\ast }}^{2}$. We shall choose an orthonormal basis in each
space $X_{j}$, $j=1,...,k$, while in $X_{j}^{\ast }$, $j=1,...,k$ we shall
consider the dual bases. Then $\partial ^{X}=\left( \partial
^{X_{1}},...,\partial ^{X_{k}}\right) $, $\partial ^{X^{\ast }}=\left(
\partial ^{X_{1}^{\ast }},...,\partial ^{X_{k}^{\ast }}\right) .$

For $1\leq p\leq \infty $ and $\mathbf{t}=\left( t_{1},...,t_{k}\right) \in 
%TCIMACRO{\U{2115} }%
%BeginExpansion
\mathbb{N}
%EndExpansion
^{k}$ we set $\mathcal{M}_{\mathbf{t}}^{p}=\mathcal{M}_{t_{1},...,t_{k}}^{p}$
for the space of all distributions $a\in \mathcal{S}^{\ast }(\mathfrak{S})$
whose derivatives $\left( \partial ^{X_{1}}\right) ^{\alpha _{1}}...\left(
\partial ^{X_{k}}\right) ^{\alpha _{k}}\left( \partial ^{X_{1}^{\ast
}}\right) ^{\beta _{1}}...\left( \partial ^{X_{k}^{\ast }}\right) ^{\beta
_{k}}a$ belong to $L^{p}\left( \mathfrak{S}\right) $ when $\alpha _{j},\beta
_{j}\in 
%TCIMACRO{\U{2115} }%
%BeginExpansion
\mathbb{N}
%EndExpansion
^{\dim X_{j}}$, $\left\vert \alpha _{j}\right\vert ,\left\vert \beta
_{j}\right\vert \leq t_{j}$, $j=1,...,k$. On this space we shall consider
the natural norm $\left\vert \cdot \right\vert _{p,\mathbf{t}}=\left\vert
\cdot \right\vert _{p,t_{1},...,t_{k}}$ defined by 
\begin{equation*}
\left\vert a\right\vert _{p,\mathbf{t}}=\max_{\left\vert \alpha
_{1}\right\vert ,\left\vert \beta _{1}\right\vert \leq t_{1},...,\left\vert
\alpha _{k}\right\vert ,\left\vert \beta _{k}\right\vert \leq
t_{k}}\left\Vert \left( \partial ^{X_{1}}\right) ^{\alpha _{1}}...\left(
\partial ^{X_{k}}\right) ^{\alpha _{k}}\left( \partial ^{X_{1}^{\ast
}}\right) ^{\beta _{1}}...\left( \partial ^{X_{k}^{\ast }}\right) ^{\beta
_{k}}a\right\Vert _{L^{p}}.
\end{equation*}

Let $m_{1}=\left[ \frac{\dim X_{1}}{2}\right] +1,...,m_{k}=\left[ \frac{\dim
X_{k}}{2}\right] +1$ and $\mathbf{m=}\left( m_{1},...,m_{k}\right) $.

A consequence of Theorem \ref{tcp2} and of Lemma \ref{tcp4} is the following

\begin{theorem}
Assume that $1\leq p<\infty $ and let $a\in \mathcal{S}^{\ast }(\mathfrak{S}%
) $.

$(\func{a})$ If $a\in \mathcal{M}_{m_{1},...,m_{k}}^{p}$, then $%
a_{X}^{0}\left( Q,P\right) $ has an extension in $\mathcal{B}_{p}\left( 
\mathcal{H}\left( X\right) \right) $ denoted also by $a_{X}^{0}\left(
Q,P\right) $ and there is $C>0$ such that 
\begin{equation*}
\left\Vert a_{X}^{0}\left( Q,P\right) \right\Vert _{p}\leq C\left\vert
a\right\vert _{p,m_{1},...,m_{k}}.
\end{equation*}

$(\func{b})$ If $a\in \mathcal{M}_{2m_{1},...,2m_{k}}^{p}$, then for any $%
\tau \in 
%TCIMACRO{\U{211d} }%
%BeginExpansion
\mathbb{R}
%EndExpansion
\equiv 
%TCIMACRO{\U{211d} }%
%BeginExpansion
\mathbb{R}
%EndExpansion
\cdot 1_{X}$, $a_{X}^{\tau }\left( Q,P\right) $ has an extension in $%
\mathcal{B}_{p}\left( \mathcal{H}\left( X\right) \right) $ denoted also by $%
a_{X}^{\tau }\left( Q,P\right) $. The mapping 
\begin{equation*}
%TCIMACRO{\U{211d} }%
%BeginExpansion
\mathbb{R}
%EndExpansion
\ni \tau \rightarrow a_{X}^{\tau }\left( Q,P\right) \in \mathcal{B}%
_{p}\left( \mathcal{H}\left( X\right) \right)
\end{equation*}%
is continuous and for any $K$ a compact subset of $%
%TCIMACRO{\U{211d} }%
%BeginExpansion
\mathbb{R}
%EndExpansion
$, there is $C_{K}>0$ such that 
\begin{equation*}
\left\Vert a_{X}^{\tau }\left( Q,P\right) \right\Vert _{p}\leq
C_{K}\left\vert a\right\vert _{p,2m_{1},...,2m_{k}},
\end{equation*}%
for any $\tau \in K$.
\end{theorem}

Similarly, for $p=\infty $, a consequence of Theorem \ref{tcp3} and of Lemma %
\ref{tcp4} is the celebrated Calderon-Vaillancourt Theorem.

\begin{theorem}[Calderon, Vaillancourt]
$(\func{a})$ If $a\in \mathcal{M}_{m_{1},...,m_{k}}^{\infty }$, then $%
a_{X}^{0}\left( Q,P\right) $ is $L^{2}$-bounded and there is $C>0$ such that 
\begin{equation*}
\left\Vert a_{X}^{0}\left( Q,P\right) \right\Vert _{\mathcal{B}\left( 
\mathcal{H}\left( X\right) \right) }\leq C\left\vert a\right\vert _{\infty
,m_{1},...,m_{k}}.
\end{equation*}

$(\func{b})$ If $a\in \mathcal{M}_{2m_{1},...,2m_{k}}^{\infty }$, then $%
a_{X}^{\tau }\left( Q,P\right) $ is $L^{2}$-bounded for any $\tau \in 
%TCIMACRO{\U{211d} }%
%BeginExpansion
\mathbb{R}
%EndExpansion
\equiv 
%TCIMACRO{\U{211d} }%
%BeginExpansion
\mathbb{R}
%EndExpansion
\cdot 1_{X}$. The mapping 
\begin{equation*}
%TCIMACRO{\U{211d} }%
%BeginExpansion
\mathbb{R}
%EndExpansion
\ni \tau \rightarrow a_{X}^{\tau }\left( Q,P\right) \in \mathcal{B}\left( 
\mathcal{H}\left( X\right) \right)
\end{equation*}%
is continuous and for any $K$ a compact subset of $%
%TCIMACRO{\U{211d} }%
%BeginExpansion
\mathbb{R}
%EndExpansion
$, there is $C_{K}>0$ such that 
\begin{equation*}
\left\Vert a_{X}^{\tau }\left( Q,P\right) \right\Vert _{\mathcal{B}\left( 
\mathcal{H}\left( X\right) \right) }\leq C_{K}\left\vert a\right\vert
_{\infty ,2m_{1},...,2m_{k}},
\end{equation*}%
for any $\tau \in K$.
\end{theorem}

To obtain other results we shall use the following

\begin{lemma}
\label{tcp5}Let $X=X_{1}\oplus ...\oplus X_{k}$ be an orthogonal
decomposition. If $s_{1},...,s_{k}\geq 0,\varepsilon >0$, then for any $%
1\leq p\leq \infty $, 
\begin{equation*}
\left( 1-\triangle _{X_{1}}\otimes 1\right) ^{s_{1}}...\left( 1-1\otimes
\triangle _{X_{k}}\right) ^{s_{k}}\left( 1-\triangle _{X}\right) ^{-\left(
s_{1}+...+s_{k}\right) -\varepsilon }\in \mathcal{B}\left( L^{p}\left(
X\right) \right) .
\end{equation*}
\end{lemma}

The proof of this lemma is given the appendix.

The next two theorems are consequences of this lemma, Theorem \ref{tcp2} and
Theorem \ref{tcp3}. Recall that the Sobolev space $H_{p}^{s}\left( \mathfrak{%
S}\right) $, $s\in 
%TCIMACRO{\U{211d} }%
%BeginExpansion
\mathbb{R}
%EndExpansion
$, $1\leq p\leq \infty $, consists of all $a\in \mathcal{S}^{\ast }(%
\mathfrak{S})$ such that $\left( 1-\triangle _{\mathfrak{S}}\right)
^{s/2}a\in L^{p}\left( \mathfrak{S}\right) $, and we set $\left\Vert
a\right\Vert _{H_{p}^{s}\left( \mathfrak{S}\right) }\equiv \left\Vert \left(
1-\triangle _{\mathfrak{S}}\right) ^{s/2}a\right\Vert _{L^{p}\left( 
\mathfrak{S}\right) }$.

\begin{theorem}
Assume that $1\leq p<\infty $.

$(\func{a})$ If $s>\dim X$ and $a\in H_{p}^{s}\left( \mathfrak{S}\right) $,
then $a_{X}^{0}\left( Q,P\right) $ has an extension in $\mathcal{B}%
_{p}\left( \mathcal{H}\left( X\right) \right) $ denoted also by $%
a_{X}^{0}\left( Q,P\right) $ and there is $C>0$ such that 
\begin{equation*}
\left\Vert a_{X}^{0}\left( Q,P\right) \right\Vert _{p}\leq C\left\Vert
a\right\Vert _{H_{p}^{s}\left( \mathfrak{S}\right) }.
\end{equation*}

$(\func{b})$ If $s>2\dim X$ and $a\in H_{p}^{s}\left( \mathfrak{S}\right) $,
then for any $\tau \in 
%TCIMACRO{\U{211d} }%
%BeginExpansion
\mathbb{R}
%EndExpansion
\equiv 
%TCIMACRO{\U{211d} }%
%BeginExpansion
\mathbb{R}
%EndExpansion
\cdot 1_{X}$, $a_{X}^{\tau }\left( Q,P\right) $ has an extension in $%
\mathcal{B}_{p}\left( \mathcal{H}\left( X\right) \right) $ denoted also by $%
a_{X}^{\tau }\left( Q,P\right) $. The mapping 
\begin{equation*}
%TCIMACRO{\U{211d} }%
%BeginExpansion
\mathbb{R}
%EndExpansion
\ni \tau \rightarrow a_{X}^{\tau }\left( Q,P\right) \in \mathcal{B}%
_{p}\left( \mathcal{H}\left( X\right) \right)
\end{equation*}%
is continuous and for any $K$ a compact subset of $%
%TCIMACRO{\U{211d} }%
%BeginExpansion
\mathbb{R}
%EndExpansion
$, there is $C_{K}>0$ such that 
\begin{equation*}
\left\Vert a_{X}^{\tau }\left( Q,P\right) \right\Vert _{p}\leq
C_{K}\left\Vert a\right\Vert _{H_{p}^{s}\left( \mathfrak{S}\right) },
\end{equation*}%
for any $\tau \in K$.
\end{theorem}

\begin{theorem}
$(\func{a})$ If $s>\dim X$ and $a\in H_{\infty }^{s}\left( \mathfrak{S}%
\right) $, then $a_{X}^{0}\left( Q,P\right) $ is $L^{2}$-bounded and there
is $C>0$ such that 
\begin{equation*}
\left\Vert a_{X}^{0}\left( Q,P\right) \right\Vert _{\mathcal{B}\left( 
\mathcal{H}\left( X\right) \right) }\leq C\left\Vert a\right\Vert
_{H_{\infty }^{s}\left( \mathfrak{S}\right) }.
\end{equation*}

$(\func{b})$ If $s>2\dim X$ and $a\in H_{\infty }^{s}\left( \mathfrak{S}%
\right) $, then $a_{X}^{\tau }\left( Q,P\right) $ is $L^{2}$-bounded for any 
$\tau \in 
%TCIMACRO{\U{211d} }%
%BeginExpansion
\mathbb{R}
%EndExpansion
\equiv 
%TCIMACRO{\U{211d} }%
%BeginExpansion
\mathbb{R}
%EndExpansion
\cdot 1_{X}$. The mapping 
\begin{equation*}
%TCIMACRO{\U{211d} }%
%BeginExpansion
\mathbb{R}
%EndExpansion
\ni \tau \rightarrow a_{X}^{\tau }\left( Q,P\right) \in \mathcal{B}\left( 
\mathcal{H}\left( X\right) \right)
\end{equation*}%
is continuous and for any $K$ a compact subset of $%
%TCIMACRO{\U{211d} }%
%BeginExpansion
\mathbb{R}
%EndExpansion
$, there is $C_{K}>0$ such that 
\begin{equation*}
\left\Vert a_{X}^{\tau }\left( Q,P\right) \right\Vert _{\mathcal{B}\left( 
\mathcal{H}\left( X\right) \right) }\leq C_{K}\left\Vert a\right\Vert
_{H_{\infty }^{s}\left( \mathfrak{S}\right) },
\end{equation*}%
for any $\tau \in K$.
\end{theorem}

If we note that $a_{X}^{\tau }\left( Q,P\right) \in \mathcal{B}_{2}\left( 
\mathcal{H}\left( X\right) \right) $ whenever $a\in L^{2}\left( \mathfrak{S}%
\right) =H_{2}^{0}\left( \mathfrak{S}\right) $, then the last two theorems
and standard interpolation results in Sobolev spaces (see \cite[Theorem 6.4.5%
]{Bergh}) give us the following

\begin{theorem}
Let $\mu >1$, $1\leq p<\infty $ and $n=\dim X$ .

$(\func{a})$ If $a\in H_{p}^{\mu n\left\vert 1-2/p\right\vert }\left( 
\mathfrak{S}\right) $, then $a_{X}^{0}\left( Q,P\right) $ has an extension
in $\mathcal{B}_{p}\left( \mathcal{H}\left( X\right) \right) $ denoted also
by $a_{X}^{0}\left( Q,P\right) $ and there is $C>0$ such that 
\begin{equation*}
\left\Vert a_{X}^{0}\left( Q,P\right) \right\Vert _{p}\leq C\left\Vert
a\right\Vert _{H_{p}^{\mu n\left\vert 1-2/p\right\vert }\left( \mathfrak{S}%
\right) }.
\end{equation*}

$(\func{b})$ If $a\in H_{p}^{2\mu n\left\vert 1-2/p\right\vert }\left( 
\mathfrak{S}\right) $, then for any $\tau \in 
%TCIMACRO{\U{211d} }%
%BeginExpansion
\mathbb{R}
%EndExpansion
\equiv 
%TCIMACRO{\U{211d} }%
%BeginExpansion
\mathbb{R}
%EndExpansion
\cdot 1_{X}$, $a_{X}^{\tau }\left( Q,P\right) $ has an extension in $%
\mathcal{B}_{p}\left( \mathcal{H}\left( X\right) \right) $ denoted also by $%
a_{X}^{\tau }\left( Q,P\right) $. The mapping 
\begin{equation*}
%TCIMACRO{\U{211d} }%
%BeginExpansion
\mathbb{R}
%EndExpansion
\ni \tau \rightarrow a_{X}^{\tau }\left( Q,P\right) \in \mathcal{B}%
_{p}\left( \mathcal{H}\left( X\right) \right)
\end{equation*}%
is continuous and for any $K$ a compact subset of $%
%TCIMACRO{\U{211d} }%
%BeginExpansion
\mathbb{R}
%EndExpansion
$, there is $C_{K}>0$ such that 
\begin{equation*}
\left\Vert a_{X}^{\tau }\left( Q,P\right) \right\Vert _{p}\leq
C_{K}\left\Vert a\right\Vert _{H_{p}^{2\mu n\left\vert 1-2/p\right\vert
}\left( \mathfrak{S}\right) },
\end{equation*}%
for any $\tau \in K$.
\end{theorem}

\appendix

\section{A class of Fourier multipliers}

In this appendix we shall prove Lemma \ref{tcp5}.

Let us assume that the lemma has been proved for $k=2$. Let $k\geq 3$. For $%
1\leq l\leq k$, we set $T_{l}=\left( 1-1\otimes \triangle _{X_{l}}\otimes
1\right) ^{s_{l}}\left( 1-\triangle _{X}\right) ^{-s_{l}-\frac{\varepsilon }{%
k}}$. Then $T_{l}\in \mathcal{B}\left( L^{p}\left( X\right) \right) $ and 
\begin{equation*}
\left( 1-\triangle _{X_{1}}\otimes 1\right) ^{s_{1}}...\left( 1-1\otimes
\triangle _{X_{k}}\right) ^{s_{k}}\left( 1-\triangle _{X}\right) ^{-\left(
s_{1}+...+s_{k}\right) -\varepsilon }=T_{1}...T_{k}\in \mathcal{B}\left(
L^{p}\left( X\right) \right) .
\end{equation*}%
This shows that the lema is true provided that it holds for $k=2$.

By choosing an orthonormal basis in each space $X_{1}$ and $X_{2}$ we
identify $X_{1}$ with $%
%TCIMACRO{\U{211d} }%
%BeginExpansion
\mathbb{R}
%EndExpansion
^{n_{1}}$, $X_{2}$ with $%
%TCIMACRO{\U{211d} }%
%BeginExpansion
\mathbb{R}
%EndExpansion
^{n_{2}}$ and $X$ with $%
%TCIMACRO{\U{211d} }%
%BeginExpansion
\mathbb{R}
%EndExpansion
^{n_{1}}\mathbb{\times }%
%TCIMACRO{\U{211d} }%
%BeginExpansion
\mathbb{R}
%EndExpansion
^{n_{2}}$. Then the operator 
\begin{equation*}
\left( 1-\triangle _{1}\otimes 1\right) ^{s_{1}}\left( 1-1\otimes \triangle
_{2}\right) ^{s_{2}}\left( 1-\triangle _{1}\otimes 1-1\otimes \triangle
_{2}\right) ^{-s_{1}-s_{2}-\varepsilon }
\end{equation*}%
has the symbol $a:%
%TCIMACRO{\U{211d} }%
%BeginExpansion
\mathbb{R}
%EndExpansion
^{n_{1}}\mathbb{\times }%
%TCIMACRO{\U{211d} }%
%BeginExpansion
\mathbb{R}
%EndExpansion
^{n_{2}}\rightarrow 
%TCIMACRO{\U{211d} }%
%BeginExpansion
\mathbb{R}
%EndExpansion
$ defined by 
\begin{equation*}
a\left( \xi _{1},\xi _{2}\right) =\left\langle \xi _{1}\right\rangle
^{s_{1}}\left\langle \xi _{2}\right\rangle ^{s_{2}}\left\langle \left( \xi
_{1},\xi _{2}\right) \right\rangle ^{-s_{1}-s_{2}-\varepsilon },\quad \left(
\xi _{1},\xi _{2}\right) \in 
%TCIMACRO{\U{211d} }%
%BeginExpansion
\mathbb{R}
%EndExpansion
^{n_{1}}\mathbb{\times }%
%TCIMACRO{\U{211d} }%
%BeginExpansion
\mathbb{R}
%EndExpansion
^{n_{2}}
\end{equation*}%
where $\left\langle \xi \right\rangle =\left( 1+\left\vert \xi \right\vert
^{2}\right) ^{\frac{1}{2}}$, $\xi \in 
%TCIMACRO{\U{211d} }%
%BeginExpansion
\mathbb{R}
%EndExpansion
^{n}$. If we write $\varepsilon =\varepsilon _{1}+\varepsilon _{2}$ with $%
\varepsilon _{1},\varepsilon _{2}>0$, then it can be easy check that for any 
$\left( \alpha _{1},\alpha _{2}\right) \in 
%TCIMACRO{\U{2115} }%
%BeginExpansion
\mathbb{N}
%EndExpansion
^{n_{1}}\times 
%TCIMACRO{\U{2115} }%
%BeginExpansion
\mathbb{N}
%EndExpansion
^{n_{2}}$, there is $C_{\alpha _{1},\alpha _{2}}=C\left(
n_{1},n_{2},s_{1},s_{2},\alpha _{1},\alpha _{2},\varepsilon \right) >0$ such
that 
\begin{equation*}
\left\vert \partial _{\xi _{1}}^{\alpha _{1}}\partial _{\xi _{2}}^{\alpha
_{2}}a\left( \xi _{1},\xi _{2}\right) \right\vert \leq C_{\alpha _{1},\alpha
_{2}}\left\langle \xi _{1}\right\rangle ^{-\varepsilon _{1-}\left\vert
\alpha _{1}\right\vert }\left\langle \xi _{2}\right\rangle ^{-\epsilon
_{2}-\left\vert \alpha _{2}\right\vert },\quad \left( \xi _{1},\xi
_{2}\right) \in 
%TCIMACRO{\U{211d} }%
%BeginExpansion
\mathbb{R}
%EndExpansion
^{n_{1}}\mathbb{\times }%
%TCIMACRO{\U{211d} }%
%BeginExpansion
\mathbb{R}
%EndExpansion
^{n_{2}}.
\end{equation*}

\begin{definition}
Let $\mathbf{m}=\left( m_{1},m_{2}\right) \in 
%TCIMACRO{\U{211d} }%
%BeginExpansion
\mathbb{R}
%EndExpansion
^{2}$. We shall say that $a:%
%TCIMACRO{\U{211d} }%
%BeginExpansion
\mathbb{R}
%EndExpansion
^{n_{1}}\mathbb{\times }%
%TCIMACRO{\U{211d} }%
%BeginExpansion
\mathbb{R}
%EndExpansion
^{n_{2}}\rightarrow 
%TCIMACRO{\U{2102} }%
%BeginExpansion
\mathbb{C}
%EndExpansion
$ is a symbol of degree $\mathbf{m}$ if $a\in \mathcal{C}^{\infty }\left( 
%TCIMACRO{\U{211d} }%
%BeginExpansion
\mathbb{R}
%EndExpansion
^{n_{1}}\mathbb{\times }%
%TCIMACRO{\U{211d} }%
%BeginExpansion
\mathbb{R}
%EndExpansion
^{n_{2}}\right) $ and for any $\left( \alpha _{1},\alpha _{2}\right) \in 
%TCIMACRO{\U{2115} }%
%BeginExpansion
\mathbb{N}
%EndExpansion
^{n_{1}}\times 
%TCIMACRO{\U{2115} }%
%BeginExpansion
\mathbb{N}
%EndExpansion
^{n_{2}}$, there is $C_{\alpha _{1},\alpha _{2}}>0$ such that 
\begin{equation*}
\left\vert \partial _{\xi _{1}}^{\alpha _{1}}\partial _{\xi _{2}}^{\alpha
_{2}}a\left( \xi _{1},\xi _{2}\right) \right\vert \leq C_{\alpha _{1},\alpha
_{2}}\left\langle \xi _{1}\right\rangle ^{m_{1-}\left\vert \alpha
_{1}\right\vert }\left\langle \xi _{2}\right\rangle ^{m_{2}-\left\vert
\alpha _{2}\right\vert },\quad \left( \xi _{1},\xi _{2}\right) \in 
%TCIMACRO{\U{211d} }%
%BeginExpansion
\mathbb{R}
%EndExpansion
^{n_{1}}\mathbb{\times }%
%TCIMACRO{\U{211d} }%
%BeginExpansion
\mathbb{R}
%EndExpansion
^{n_{2}}.
\end{equation*}
\end{definition}

We denote by $\mathcal{S}^{\mathbf{m}}\left( 
%TCIMACRO{\U{211d} }%
%BeginExpansion
\mathbb{R}
%EndExpansion
^{n_{1}}\mathbb{\times }%
%TCIMACRO{\U{211d} }%
%BeginExpansion
\mathbb{R}
%EndExpansion
^{n_{2}}\right) =\mathcal{S}^{m_{1},m_{2}}\left( 
%TCIMACRO{\U{211d} }%
%BeginExpansion
\mathbb{R}
%EndExpansion
^{n_{1}}\mathbb{\times }%
%TCIMACRO{\U{211d} }%
%BeginExpansion
\mathbb{R}
%EndExpansion
^{n_{2}}\right) $ the vector space of all symbols of degree $\mathbf{m}$ and
observe that 
\begin{gather*}
\mathcal{S}\left( 
%TCIMACRO{\U{211d} }%
%BeginExpansion
\mathbb{R}
%EndExpansion
^{n_{1}}\mathbb{\times }%
%TCIMACRO{\U{211d} }%
%BeginExpansion
\mathbb{R}
%EndExpansion
^{n_{2}}\right) \subset \mathcal{S}^{\mathbf{m}}\left( 
%TCIMACRO{\U{211d} }%
%BeginExpansion
\mathbb{R}
%EndExpansion
^{n_{1}}\mathbb{\times }%
%TCIMACRO{\U{211d} }%
%BeginExpansion
\mathbb{R}
%EndExpansion
^{n_{2}}\right) \subset \mathcal{S}^{\ast }\left( 
%TCIMACRO{\U{211d} }%
%BeginExpansion
\mathbb{R}
%EndExpansion
^{n_{1}}\mathbb{\times }%
%TCIMACRO{\U{211d} }%
%BeginExpansion
\mathbb{R}
%EndExpansion
^{n_{2}}\right) , \\
\mathcal{S}^{\mathbf{m}}\left( 
%TCIMACRO{\U{211d} }%
%BeginExpansion
\mathbb{R}
%EndExpansion
^{n_{1}}\mathbb{\times }%
%TCIMACRO{\U{211d} }%
%BeginExpansion
\mathbb{R}
%EndExpansion
^{n_{2}}\right) \cdot \mathcal{S}^{\overline{\mathbf{m}}}\left( 
%TCIMACRO{\U{211d} }%
%BeginExpansion
\mathbb{R}
%EndExpansion
^{n_{1}}\mathbb{\times }%
%TCIMACRO{\U{211d} }%
%BeginExpansion
\mathbb{R}
%EndExpansion
^{n_{2}}\right) \subset \mathcal{S}^{\mathbf{m+}\overline{\mathbf{m}}}\left( 
%TCIMACRO{\U{211d} }%
%BeginExpansion
\mathbb{R}
%EndExpansion
^{n_{1}}\mathbb{\times }%
%TCIMACRO{\U{211d} }%
%BeginExpansion
\mathbb{R}
%EndExpansion
^{n_{2}}\right) , \\
\mathbf{m\leq }\overline{\mathbf{m}}\Rightarrow \mathcal{S}^{\mathbf{m}%
}\left( 
%TCIMACRO{\U{211d} }%
%BeginExpansion
\mathbb{R}
%EndExpansion
^{n_{1}}\mathbb{\times }%
%TCIMACRO{\U{211d} }%
%BeginExpansion
\mathbb{R}
%EndExpansion
^{n_{2}}\right) \subset \mathcal{S}^{\overline{\mathbf{m}}}\left( 
%TCIMACRO{\U{211d} }%
%BeginExpansion
\mathbb{R}
%EndExpansion
^{n_{1}}\mathbb{\times }%
%TCIMACRO{\U{211d} }%
%BeginExpansion
\mathbb{R}
%EndExpansion
^{n_{2}}\right) ,
\end{gather*}%
where $\mathbf{m\leq }\overline{\mathbf{m}}$ means $m_{1}\leq \overline{m}%
_{1},m_{2}\leq \overline{m}_{2}$.

Let $\varphi \in \mathcal{C}_{0}^{\infty }\left( 
%TCIMACRO{\U{211d} }%
%BeginExpansion
\mathbb{R}
%EndExpansion
^{n}\right) $ such that $0\leq \varphi \leq 1,$ $\varphi \left( \xi \right)
=1$ for $\left\vert \xi \right\vert \leq 1$, $\varphi \left( \xi \right) =0$
for $\left\vert \xi \right\vert \geq 2$. We define $\psi \in \mathcal{C}%
_{0}^{\infty }\left( 
%TCIMACRO{\U{211d} }%
%BeginExpansion
\mathbb{R}
%EndExpansion
^{n}\right) $ by 
\begin{equation*}
\psi \left( \xi \right) =-\xi \cdot \nabla \varphi \left( \xi \right)
=-\sum_{j=1}^{n}\xi _{j}\partial _{j}\varphi \left( \xi \right) ,\quad \xi
\in 
%TCIMACRO{\U{211d} }%
%BeginExpansion
\mathbb{R}
%EndExpansion
^{n}.
\end{equation*}%
Then $\limfunc{supp}\psi \subset \left\{ \xi \in 
%TCIMACRO{\U{211d} }%
%BeginExpansion
\mathbb{R}
%EndExpansion
^{n}:1\leq \left\vert \xi \right\vert \leq 2\right\} $.

If $r>1$ and $\xi \in 
%TCIMACRO{\U{211d} }%
%BeginExpansion
\mathbb{R}
%EndExpansion
^{n}\backslash 0$, then 
\begin{eqnarray*}
\varphi \left( \frac{\xi }{r}\right) -\varphi \left( \xi \right)
&=&\int_{1}^{r}\frac{\func{d}}{\func{d}t}\varphi \left( \frac{\xi }{t}%
\right) \func{d}t=-\int_{1}^{r}\frac{\xi }{t}\cdot \nabla \varphi \left( 
\frac{\xi }{t}\right) \frac{\func{d}t}{t} \\
&=&\int_{1}^{r}\psi \left( \frac{\xi }{t}\right) \frac{\func{d}t}{t}.
\end{eqnarray*}%
We take $r\rightarrow \infty $ to obtain that 
\begin{equation*}
1=\varphi \left( \xi \right) +\int_{1}^{\infty }\psi \left( \frac{\xi }{t}%
\right) \frac{\func{d}t}{t},\quad \xi \in 
%TCIMACRO{\U{211d} }%
%BeginExpansion
\mathbb{R}
%EndExpansion
^{n}.
\end{equation*}%
If we use such identities in each space $%
%TCIMACRO{\U{211d} }%
%BeginExpansion
\mathbb{R}
%EndExpansion
^{n_{1}}$ and $%
%TCIMACRO{\U{211d} }%
%BeginExpansion
\mathbb{R}
%EndExpansion
^{n_{2}}$, we get 
\begin{gather*}
1=\varphi _{1}\left( \xi _{1}\right) \varphi _{2}\left( \xi _{2}\right)
+\int_{1}^{\infty }\psi _{1}\left( \frac{\xi _{1}}{t_{1}}\right) \varphi
_{2}\left( \xi _{2}\right) \frac{\func{d}t_{1}}{t_{1}}+\int_{1}^{\infty
}\varphi _{1}\left( \xi _{1}\right) \psi _{2}\left( \frac{\xi _{2}}{t_{2}}%
\right) \frac{\func{d}t_{2}}{t_{2}} \\
+\int_{1}^{\infty }\int_{1}^{\infty }\psi _{1}\left( \frac{\xi _{1}}{t_{1}}%
\right) \psi _{2}\left( \frac{\xi _{2}}{t_{2}}\right) \frac{\func{d}t_{1}}{%
t_{1}}\frac{\func{d}t_{2}}{t_{2}},\quad \left( \xi _{1},\xi _{2}\right) \in 
%TCIMACRO{\U{211d} }%
%BeginExpansion
\mathbb{R}
%EndExpansion
^{n_{1}}\mathbb{\times }%
%TCIMACRO{\U{211d} }%
%BeginExpansion
\mathbb{R}
%EndExpansion
^{n_{2}}.
\end{gather*}%
where for $j=1,2$ the functions $\varphi _{j},\psi _{j}\in \mathcal{C}%
_{0}^{\infty }\left( 
%TCIMACRO{\U{211d} }%
%BeginExpansion
\mathbb{R}
%EndExpansion
^{n_{j}}\right) $ satisfy $0\leq \varphi _{j}\leq 1$, $\varphi _{j}\left(
\xi _{j}\right) =1$ for $\left\vert \xi _{j}\right\vert \leq 1$, $\varphi
_{j}\left( \xi _{j}\right) =0$ for $\left\vert \xi _{j}\right\vert \geq 2$, $%
\limfunc{supp}\psi _{j}\subset \{\xi _{j}\in 
%TCIMACRO{\U{211d} }%
%BeginExpansion
\mathbb{R}
%EndExpansion
^{n_{j}}:1\leq $ $\left\vert \xi _{j}\right\vert \leq 2\}$, and $\psi _{j}$
is given by $\psi _{j}\left( \xi _{j}\right) =-\xi _{j}\cdot \nabla _{\xi
_{j}}\varphi _{j}\left( \xi _{j}\right) $, $\xi _{j}\in 
%TCIMACRO{\U{211d} }%
%BeginExpansion
\mathbb{R}
%EndExpansion
^{n_{j}}$.

We make use of the following simple but important remark. If $a\in \mathcal{S%
}^{\mathbf{m}}\left( 
%TCIMACRO{\U{211d} }%
%BeginExpansion
\mathbb{R}
%EndExpansion
^{n_{1}}\mathbb{\times }%
%TCIMACRO{\U{211d} }%
%BeginExpansion
\mathbb{R}
%EndExpansion
^{n_{2}}\right) $ $=\mathcal{S}^{m_{1},m_{2}}\left( 
%TCIMACRO{\U{211d} }%
%BeginExpansion
\mathbb{R}
%EndExpansion
^{n_{1}}\mathbb{\times }%
%TCIMACRO{\U{211d} }%
%BeginExpansion
\mathbb{R}
%EndExpansion
^{n_{2}}\right) $, then the families of functions 
\begin{equation*}
\left\{ \left( \psi _{1}\otimes \varphi _{2}\right) a_{t_{1},1}\right\}
_{1\leq t_{1}<\infty },\left\{ \left( \varphi _{1}\otimes \psi _{2}\right)
a_{1,t_{2}}\right\} _{1\leq t_{2}<\infty },\left\{ \left( \psi _{1}\otimes
\psi _{2}\right) a_{t_{1},t_{2}}\right\} _{\substack{ 1\leq t_{1}<\infty  \\ %
1\leq t_{2}<\infty }}
\end{equation*}%
are bounded in $\mathcal{S}\left( 
%TCIMACRO{\U{211d} }%
%BeginExpansion
\mathbb{R}
%EndExpansion
^{n_{1}}\mathbb{\times }%
%TCIMACRO{\U{211d} }%
%BeginExpansion
\mathbb{R}
%EndExpansion
^{n_{2}}\right) $, where for $t_{1},t_{2}>0$ 
\begin{equation*}
a_{t_{1},t_{2}}\left( \xi _{1},\xi _{2}\right)
=t_{1}^{-m_{1}}t_{2}^{-m_{2}}a\left( t_{1}\xi _{1},t_{2}\xi _{2}\right)
,\quad \left( \xi _{1},\xi _{2}\right) \in 
%TCIMACRO{\U{211d} }%
%BeginExpansion
\mathbb{R}
%EndExpansion
^{n_{1}}\mathbb{\times }%
%TCIMACRO{\U{211d} }%
%BeginExpansion
\mathbb{R}
%EndExpansion
^{n_{2}}.
\end{equation*}%
Let $a_{0}=\left( \varphi _{1}\otimes \varphi _{2}\right) a\in \mathcal{C}%
_{0}^{\infty }\left( 
%TCIMACRO{\U{211d} }%
%BeginExpansion
\mathbb{R}
%EndExpansion
^{n_{1}}\mathbb{\times }%
%TCIMACRO{\U{211d} }%
%BeginExpansion
\mathbb{R}
%EndExpansion
^{n_{2}}\right) $. Since 
\begin{equation*}
a\left( \xi _{1},\xi _{2}\right)
=t_{1}^{m_{1}}t_{2}^{m_{2}}a_{t_{1},t_{2}}\left( \frac{\xi _{1}}{t_{1}},%
\frac{\xi _{2}}{t_{2}}\right) ,\quad \left( \xi _{1},\xi _{2}\right) \in 
%TCIMACRO{\U{211d} }%
%BeginExpansion
\mathbb{R}
%EndExpansion
^{n_{1}}\mathbb{\times }%
%TCIMACRO{\U{211d} }%
%BeginExpansion
\mathbb{R}
%EndExpansion
^{n_{2}},1\leq t_{1},t_{2}<\infty ,
\end{equation*}%
we can write%
\begin{multline*}
a\left( \xi _{1},\xi _{2}\right) =a_{0}\left( \xi _{1},\xi _{2}\right)
+\int_{1}^{\infty }t_{1}^{m_{1}}\left( \left( \psi _{1}\otimes \varphi
_{2}\right) a_{t_{1},1}\right) \left( \frac{\xi _{1}}{t_{1}},\xi _{2}\right) 
\frac{\func{d}t_{1}}{t_{1}} \\
+\int_{1}^{\infty }t_{2}^{m_{2}}\left( \left( \varphi _{1}\otimes \psi
_{2}\right) a_{1,t_{2}}\right) \left( \xi _{1},\frac{\xi _{2}}{t_{2}}\right) 
\frac{\func{d}t_{2}}{t_{2}}\qquad \qquad \qquad \\
\qquad +\int_{1}^{\infty }\int_{1}^{\infty }t_{1}^{m_{1}}t_{2}^{m_{2}}\left(
\left( \psi _{1}\otimes \psi _{2}\right) a_{t_{1},t_{2}}\right) \left( \frac{%
\xi _{1}}{t_{1}},\frac{\xi _{2}}{t_{2}}\right) \frac{\func{d}t_{1}}{t_{1}}%
\frac{\func{d}t_{2}}{t_{2}}, \\
\left( \xi _{1},\xi _{2}\right) \in 
%TCIMACRO{\U{211d} }%
%BeginExpansion
\mathbb{R}
%EndExpansion
^{n_{1}}\mathbb{\times }%
%TCIMACRO{\U{211d} }%
%BeginExpansion
\mathbb{R}
%EndExpansion
^{n_{2}},
\end{multline*}

We have 
\begin{multline*}
\left\vert t_{1}^{m_{1}}t_{2}^{m_{2}}\left( \left( \psi _{1}\otimes \psi
_{2}\right) a_{t_{1},t_{2}}\right) \left( \frac{\xi _{1}}{t_{1}},\frac{\xi
_{2}}{t_{2}}\right) \right\vert \\
\leq t_{1}^{m_{1}-2N_{1}}t_{2}^{m_{2}-2N_{2}}\left( \sup \left\vert \left(
\psi _{1}\otimes \psi _{2}\right) a_{t_{1},t_{2}}\right\vert \right)
\left\vert \xi _{1}\right\vert ^{2N_{1}}\left\vert \xi _{2}\right\vert
^{2N_{2}}, \\
\left( \xi _{1},\xi _{2}\right) \in 
%TCIMACRO{\U{211d} }%
%BeginExpansion
\mathbb{R}
%EndExpansion
^{n_{1}}\mathbb{\times }%
%TCIMACRO{\U{211d} }%
%BeginExpansion
\mathbb{R}
%EndExpansion
^{n_{2}},1\leq t_{1},t_{2}<\infty ,
\end{multline*}%
\begin{multline*}
\left\vert t_{1}^{m_{1}}\left( \left( \psi _{1}\otimes \varphi _{2}\right)
a_{t_{1},1}\right) \left( \frac{\xi _{1}}{t_{1}},\xi _{2}\right) \right\vert
\leq t_{1}^{m_{1}-2N_{1}}\left( \sup \left\vert \left( \psi _{1}\otimes
\varphi _{2}\right) a_{t_{1},1}\right\vert \right) \left\vert \xi
_{1}\right\vert ^{2N_{1}}, \\
\left( \xi _{1},\xi _{2}\right) \in 
%TCIMACRO{\U{211d} }%
%BeginExpansion
\mathbb{R}
%EndExpansion
^{n_{1}}\mathbb{\times }%
%TCIMACRO{\U{211d} }%
%BeginExpansion
\mathbb{R}
%EndExpansion
^{n_{2}},1\leq t_{1}<\infty ,
\end{multline*}%
\begin{multline*}
\left\vert t_{2}^{m_{2}}\left( \left( \varphi _{1}\otimes \psi _{2}\right)
a_{1,t_{2}}\right) \left( \xi _{1},\frac{\xi _{2}}{t_{2}}\right) \right\vert
\leq t_{2}^{m_{2}-2N_{2}}\left( \sup \left( \varphi _{1}\otimes \psi
_{2}\right) a_{1,t_{2}}\right) \left\vert \xi _{2}\right\vert ^{2N_{2}}, \\
\left( \xi _{1},\xi _{2}\right) \in 
%TCIMACRO{\U{211d} }%
%BeginExpansion
\mathbb{R}
%EndExpansion
^{n_{1}}\mathbb{\times }%
%TCIMACRO{\U{211d} }%
%BeginExpansion
\mathbb{R}
%EndExpansion
^{n_{2}},1\leq t_{2}<\infty ,
\end{multline*}%
so if we choose $N_{1},N_{2}\in \mathbb{%
%TCIMACRO{\U{2115} }%
%BeginExpansion
\mathbb{N}
%EndExpansion
}$ such that $m_{1}<2N_{1}$ and $m_{2}<2N_{2}$, it follows that the
representation 
\begin{gather*}
a=a_{0}+\int_{1}^{\infty }\left( \left( \psi _{1}\otimes \varphi _{2}\right)
a_{t_{1},1}\right) _{t_{1}^{-1},1}\frac{\func{d}t_{1}}{t_{1}}%
+\int_{1}^{\infty }\left( \left( \varphi _{1}\otimes \psi _{2}\right)
a_{1,t_{2}}\right) _{1,t_{2}^{-1}}\frac{\func{d}t_{2}}{t_{2}} \\
+\int_{1}^{\infty }\int_{1}^{\infty }\left( \left( \psi _{1}\otimes \psi
_{2}\right) a_{t_{1},t_{2}}\right) _{t_{1}^{-1},t_{2}^{-1}}\frac{\func{d}%
t_{1}}{t_{1}}\frac{\func{d}t_{2}}{t_{2}}.
\end{gather*}%
is valid also as equality in $\mathcal{S}^{\ast }\left( 
%TCIMACRO{\U{211d} }%
%BeginExpansion
\mathbb{R}
%EndExpansion
^{n_{1}}\mathbb{\times }%
%TCIMACRO{\U{211d} }%
%BeginExpansion
\mathbb{R}
%EndExpansion
^{n_{2}}\right) $ with the integrals weakly absolutely convergent. If we
apply $\mathcal{F}^{-1}$ to this formula we obtain a decomposition of $%
\mathcal{F}^{-1}a$. Thus 
\begin{multline*}
\mathcal{F}^{-1}a=\mathcal{F}^{-1}a_{0}+\int_{1}^{\infty }\mathcal{F}%
^{-1}\left( \left( \psi _{1}\otimes \varphi _{2}\right) a_{t_{1},1}\right)
_{t_{1}^{-1},1}\frac{\func{d}t_{1}}{t_{1}} \\
+\int_{1}^{\infty }\mathcal{F}^{-1}\left( \left( \varphi _{1}\otimes \psi
_{2}\right) a_{1,t_{2}}\right) _{1,t_{2}^{-1}}\frac{\func{d}t_{2}}{t_{2}} \\
+\int_{1}^{\infty }\int_{1}^{\infty }\mathcal{F}^{-1}\left( \left( \psi
_{1}\otimes \psi _{2}\right) a_{t_{1},t_{2}}\right) _{t_{1}^{-1},t_{2}^{-1}}%
\frac{\func{d}t_{1}}{t_{1}}\frac{\func{d}t_{2}}{t_{2}}
\end{multline*}%
with the integrals weakly absolutely convergent.

We have 
\begin{multline*}
\mathcal{F}^{-1}\left( \left( \left( \psi _{1}\otimes \psi _{2}\right)
a_{t_{1},t_{2}}\right) _{t_{1}^{-1},t_{2}^{-1}}\right) \left(
x_{1},x_{2}\right) \\
=t_{1}^{m_{1}+n_{1}}t_{2}^{m_{2}+n_{2}}\mathcal{F}^{-1}\left( \left( \psi
_{1}\otimes \psi _{2}\right) a_{t_{1},t_{2}}\right) \left(
t_{1}x_{1},t_{2}x_{2}\right) ,
\end{multline*}%
\begin{eqnarray*}
\mathcal{F}^{-1}\left( \left( \left( \psi _{1}\otimes \varphi _{2}\right)
a_{t_{1},1}\right) _{t_{1}^{-1},1}\right) \left( x_{1},x_{2}\right)
&=&t_{1}^{m_{1}+n_{1}}\mathcal{F}^{-1}\left( \left( \psi _{1}\otimes \varphi
_{2}\right) a_{t_{1},1}\right) \left( t_{1}x_{1},x_{2}\right) , \\
\mathcal{F}^{-1}\left( \left( \left( \varphi _{1}\otimes \psi _{2}\right)
a_{1,t_{2}}\right) _{1,t_{2}^{-1}}\right) \left( x_{1},x_{2}\right)
&=&t_{2}^{m_{2}+n_{2}}\mathcal{F}^{-1}\left( \left( \varphi _{1}\otimes \psi
_{2}\right) a_{1,t_{2}}\right) \left( x_{1},t_{2}x_{2}\right) .
\end{eqnarray*}

Since the families of functions $\left\{ \left( \psi _{1}\otimes \varphi
_{2}\right) a_{t_{1},1}\right\} _{1\leq t_{1}<\infty }$, $\left\{ \left(
\varphi _{1}\otimes \psi _{2}\right) a_{1,t_{2}}\right\} _{1\leq
t_{2}<\infty }$, $\left\{ \left( \psi _{1}\otimes \psi _{2}\right)
a_{t_{1},t_{2}}\right\} _{\substack{ 1\leq t_{1}<\infty  \\ 1\leq
t_{2}<\infty }}$ are bounded in $\mathcal{S}\left( 
%TCIMACRO{\U{211d} }%
%BeginExpansion
\mathbb{R}
%EndExpansion
^{n_{1}}\mathbb{\times }%
%TCIMACRO{\U{211d} }%
%BeginExpansion
\mathbb{R}
%EndExpansion
^{n_{2}}\right) $, it follows that the families 
\begin{gather*}
\left\{ \mathcal{F}^{-1}\left( \left( \psi _{1}\otimes \varphi _{2}\right)
a_{t_{1},1}\right) \right\} _{1\leq t_{1}<\infty },\left\{ \mathcal{F}%
^{-1}\left( \left( \varphi _{1}\otimes \psi _{2}\right) a_{1,t_{2}}\right)
\right\} _{1\leq t_{2}<\infty }, \\
\left\{ \mathcal{F}^{-1}\left( \left( \psi _{1}\otimes \psi _{2}\right)
a_{t_{1},t_{2}}\right) \right\} _{\substack{ 1\leq t_{1}<\infty  \\ 1\leq
t_{2}<\infty }},
\end{gather*}%
are also bounded in $\mathcal{S}\left( 
%TCIMACRO{\U{211d} }%
%BeginExpansion
\mathbb{R}
%EndExpansion
^{n_{1}}\mathbb{\times }%
%TCIMACRO{\U{211d} }%
%BeginExpansion
\mathbb{R}
%EndExpansion
^{n_{2}}\right) $. In particular, for any $N\in 
%TCIMACRO{\U{2115} }%
%BeginExpansion
\mathbb{N}
%EndExpansion
$, there is $C_{N}>0$ such that if $\left( x_{1},x_{2}\right) \in 
%TCIMACRO{\U{211d} }%
%BeginExpansion
\mathbb{R}
%EndExpansion
^{n_{1}}\mathbb{\times }%
%TCIMACRO{\U{211d} }%
%BeginExpansion
\mathbb{R}
%EndExpansion
^{n_{2}}$, then%
\begin{eqnarray*}
\left\vert \mathcal{F}^{-1}\left( \left( \psi _{1}\otimes \varphi
_{2}\right) a_{t_{1},1}\right) \left( x_{1},x_{2}\right) \right\vert &\leq
&C_{N}\left\langle \left( x_{1},x_{2}\right) \right\rangle ^{-2N-M}, \\
\left\vert \mathcal{F}^{-1}\left( \left( \varphi _{1}\otimes \psi
_{2}\right) a_{1,t_{2}}\right) \left( x_{1},x_{2}\right) \right\vert &\leq
&C_{N}\left\langle \left( x_{1},x_{2}\right) \right\rangle ^{-2N-M}, \\
\left\vert \mathcal{F}^{-1}\left( \left( \psi _{1}\otimes \psi _{2}\right)
a_{t_{1},t_{2}}\right) \left( x_{1},x_{2}\right) \right\vert &\leq
&C_{N}\left\langle \left( x_{1},x_{2}\right) \right\rangle ^{-2N-2M},
\end{eqnarray*}%
where $M\in 
%TCIMACRO{\U{2115} }%
%BeginExpansion
\mathbb{N}
%EndExpansion
$, $M\geq 1+\max \left\{ 0,m_{1}+n_{1},m_{2}+n_{2}\right\} $ is fixed.

It follows that 
\begin{multline*}
\left\vert \mathcal{F}^{-1}\left( \left( \left( \psi _{1}\otimes \psi
_{2}\right) a_{t_{1},t_{2}}\right) _{t_{1}^{-1},t_{2}^{-1}}\right) \left(
x_{1},x_{2}\right) \right\vert \\
\leq C_{N}t_{1}^{m_{1}+n_{1}}t_{2}^{m_{2}+n_{2}}\left\langle
x_{1}\right\rangle ^{-N}\left\langle x_{2}\right\rangle ^{-N}\left\langle
t_{1}x_{1}\right\rangle ^{-M}\left\langle t_{2}x_{2}\right\rangle ^{-M}, \\
\left( x_{1},x_{2}\right) \in 
%TCIMACRO{\U{211d} }%
%BeginExpansion
\mathbb{R}
%EndExpansion
^{n_{1}}\mathbb{\times }%
%TCIMACRO{\U{211d} }%
%BeginExpansion
\mathbb{R}
%EndExpansion
^{n_{2}},1\leq t_{1},t_{2}<\infty ,
\end{multline*}%
\begin{multline*}
\left\vert \mathcal{F}^{-1}\left( \left( \left( \psi _{1}\otimes \varphi
_{2}\right) a_{t_{1},1}\right) _{t_{1}^{-1},1}\right) \left(
x_{1},x_{2}\right) \right\vert \leq C_{N}t_{1}^{m_{1}+n_{1}}\left\langle
x_{1}\right\rangle ^{-N}\left\langle x_{2}\right\rangle ^{-N}\left\langle
t_{1}x_{1}\right\rangle ^{-M}, \\
\left( x_{1},x_{2}\right) \in 
%TCIMACRO{\U{211d} }%
%BeginExpansion
\mathbb{R}
%EndExpansion
^{n_{1}}\mathbb{\times }%
%TCIMACRO{\U{211d} }%
%BeginExpansion
\mathbb{R}
%EndExpansion
^{n_{2}},1\leq t_{1}<\infty ,
\end{multline*}%
and 
\begin{multline*}
\left\vert \mathcal{F}^{-1}\left( \left( \left( \varphi _{1}\otimes \psi
_{2}\right) a_{1,t_{2}}\right) _{1,t_{2}^{-1}}\right) \left(
x_{1},x_{2}\right) \right\vert \leq C_{N}t_{2}^{m_{2}+n_{2}}\left\langle
x_{1}\right\rangle ^{-N}\left\langle x_{2}\right\rangle ^{-N}\left\langle
t_{2}x_{2}\right\rangle ^{-M}, \\
\left( x_{1},x_{2}\right) \in 
%TCIMACRO{\U{211d} }%
%BeginExpansion
\mathbb{R}
%EndExpansion
^{n_{1}}\mathbb{\times }%
%TCIMACRO{\U{211d} }%
%BeginExpansion
\mathbb{R}
%EndExpansion
^{n_{2}},1\leq t_{2}<\infty ,
\end{multline*}

To calculate the restriction of the temperate distribution $\mathcal{F}%
^{-1}a $ to the complement of the singular subset $\mathcal{M}=\left\{
\left( x_{1},x_{2}\right) \in 
%TCIMACRO{\U{211d} }%
%BeginExpansion
\mathbb{R}
%EndExpansion
^{n_{1}}\mathbb{\times }%
%TCIMACRO{\U{211d} }%
%BeginExpansion
\mathbb{R}
%EndExpansion
^{n_{2}}:\left\vert x_{1}\right\vert \left\vert x_{2}\right\vert =0\right\} $%
, we need the following easy consequence of Fubini theorem.

\begin{lemma}
Let $\left( T,\mu \right) $ be a measure space, $\Omega \subset 
%TCIMACRO{\U{211d} }%
%BeginExpansion
\mathbb{R}
%EndExpansion
^{n}$ an open set and $f:\Omega \times T\rightarrow 
%TCIMACRO{\U{2102} }%
%BeginExpansion
\mathbb{C}
%EndExpansion
$ a measurable function.

$(\func{a})$ If for any $\varphi \in \mathcal{C}_{0}^{\infty }\left( \Omega
\right) $ the function $\Omega \times T\ni \left( x,t\right) \rightarrow
\varphi \left( x\right) f\left( x,t\right) \in 
%TCIMACRO{\U{2102} }%
%BeginExpansion
\mathbb{C}
%EndExpansion
$ belongs to $L^{1}\left( \Omega \times T\right) $, then the mapping%
\begin{equation*}
\mathcal{C}_{0}^{\infty }\left( \Omega \right) \ni \varphi \rightarrow \iint
\varphi \left( x\right) f\left( x,t\right) \limfunc{d}x\limfunc{d}\mu \left(
t\right) \in 
%TCIMACRO{\U{2102}}%
%BeginExpansion
\mathbb{C}%
%EndExpansion
\end{equation*}%
define a distribution, the function $\Omega \ni x\rightarrow \int f\left(
x,t\right) \limfunc{d}\mu \left( t\right) \in 
%TCIMACRO{\U{2102} }%
%BeginExpansion
\mathbb{C}
%EndExpansion
$, defined a.e., belongs to $L_{loc}^{1}\left( \Omega \right) $ and we have 
\begin{eqnarray*}
\left\langle \varphi ,\int f\left( \cdot ,t\right) \limfunc{d}\mu \left(
t\right) \right\rangle &=&\iint \varphi \left( x\right) f\left( x,t\right) 
\limfunc{d}x\limfunc{d}\mu \left( t\right) \\
&=&\int \left( \int \varphi \left( x\right) f\left( x,t\right) \limfunc{d}%
x\right) \limfunc{d}\mu \left( t\right) ,\quad \varphi \in \mathcal{C}%
_{0}^{\infty }\left( \Omega \right) .
\end{eqnarray*}

$(\func{b})$ Assume that $\Omega =%
%TCIMACRO{\U{211d} }%
%BeginExpansion
\mathbb{R}
%EndExpansion
^{n}$. If there is $\tau \in 
%TCIMACRO{\U{211d} }%
%BeginExpansion
\mathbb{R}
%EndExpansion
$ such that the function $%
%TCIMACRO{\U{211d} }%
%BeginExpansion
\mathbb{R}
%EndExpansion
^{n}\times T\ni \left( x,t\right) \rightarrow \left\langle x\right\rangle
^{-\tau }f\left( x,t\right) \in 
%TCIMACRO{\U{2102} }%
%BeginExpansion
\mathbb{C}
%EndExpansion
$ belongs to $L^{1}\left( 
%TCIMACRO{\U{211d} }%
%BeginExpansion
\mathbb{R}
%EndExpansion
^{n}\times T\right) $, then the mapping%
\begin{equation*}
\mathcal{S}\left( 
%TCIMACRO{\U{211d} }%
%BeginExpansion
\mathbb{R}
%EndExpansion
^{n}\right) \ni \varphi \rightarrow \iint \varphi \left( x\right) f\left(
x,t\right) \limfunc{d}x\limfunc{d}\mu \left( t\right) \in 
%TCIMACRO{\U{2102}}%
%BeginExpansion
\mathbb{C}%
%EndExpansion
\end{equation*}%
define a temperate distribution, the function $%
%TCIMACRO{\U{211d} }%
%BeginExpansion
\mathbb{R}
%EndExpansion
^{n}\ni x\rightarrow \int f\left( x,t\right) \limfunc{d}\mu \left( t\right)
\in 
%TCIMACRO{\U{2102} }%
%BeginExpansion
\mathbb{C}
%EndExpansion
$, defined a.e., belongs to $L_{loc}^{1}\left( 
%TCIMACRO{\U{211d} }%
%BeginExpansion
\mathbb{R}
%EndExpansion
^{n}\right) $, $\left\langle \cdot \right\rangle ^{-\tau }\int f\left( \cdot
,t\right) \limfunc{d}\mu \left( t\right) \in L^{1}\left( 
%TCIMACRO{\U{211d} }%
%BeginExpansion
\mathbb{R}
%EndExpansion
^{n}\right) $ and we have%
\begin{eqnarray*}
\left\langle \varphi ,\int f\left( \cdot ,t\right) \limfunc{d}\mu \left(
t\right) \right\rangle &=&\iint \varphi \left( x\right) f\left( x,t\right) 
\limfunc{d}x\limfunc{d}\mu \left( t\right) \\
&=&\int \left( \int \varphi \left( x\right) f\left( x,t\right) \limfunc{d}%
x\right) \limfunc{d}\mu \left( t\right) ,\quad \varphi \in \mathcal{S}\left( 
%TCIMACRO{\U{211d} }%
%BeginExpansion
\mathbb{R}
%EndExpansion
^{n}\right) .
\end{eqnarray*}
\end{lemma}

From the representation formula of $\mathcal{F}^{-1}a$, the above estimates
and part $(\func{a})$ of the previous lemma we conclude that $\mathcal{F}%
^{-1}a$ is a $\mathcal{C}^{\infty }$ function on $%
%TCIMACRO{\U{211d} }%
%BeginExpansion
\mathbb{R}
%EndExpansion
^{n_{1}}\mathbb{\times }%
%TCIMACRO{\U{211d} }%
%BeginExpansion
\mathbb{R}
%EndExpansion
^{n_{2}}\backslash \mathcal{M}$ which decays at infinity, together with all
its derivatives, more rapidly than any power of $\left\langle \left(
x_{1},x_{2}\right) \right\rangle ^{-1}$. We shall denote by \QTR{cal}{S}$%
\left( 
%TCIMACRO{\U{211d} }%
%BeginExpansion
\mathbb{R}
%EndExpansion
^{n_{1}}\mathbb{\times }%
%TCIMACRO{\U{211d} }%
%BeginExpansion
\mathbb{R}
%EndExpansion
^{n_{2}}\backslash \mathcal{M}\right) $ the space of elements of $\mathcal{S}%
^{\ast }\left( 
%TCIMACRO{\U{211d} }%
%BeginExpansion
\mathbb{R}
%EndExpansion
^{n_{1}}\mathbb{\times }%
%TCIMACRO{\U{211d} }%
%BeginExpansion
\mathbb{R}
%EndExpansion
^{n_{2}}\right) $ that are smooth outside the set $\mathcal{M}$ and rapidly
decreasing at infinity. Since on $%
%TCIMACRO{\U{211d} }%
%BeginExpansion
\mathbb{R}
%EndExpansion
^{n_{1}}\mathbb{\times }%
%TCIMACRO{\U{211d} }%
%BeginExpansion
\mathbb{R}
%EndExpansion
^{n_{2}}\backslash \mathcal{M}$ we have 
\begin{multline*}
\mathcal{F}^{-1}a\left( x_{1},x_{2}\right) =\mathcal{F}^{-1}a_{0}\left(
x_{1},x_{2}\right) +\int_{1}^{\infty }\mathcal{F}^{-1}\left( \left( \left(
\psi _{1}\otimes \varphi _{2}\right) a_{t_{1},1}\right)
_{t_{1}^{-1},1}\right) \left( x_{1},x_{2}\right) \frac{\func{d}t_{1}}{t_{1}}
\\
+\int_{1}^{\infty }\mathcal{F}^{-1}\left( \left( \left( \varphi _{1}\otimes
\psi _{2}\right) a_{1,t_{2}}\right) _{1,t_{2}^{-1}}\right) \left(
x_{1},x_{2}\right) \frac{\func{d}t_{2}}{t_{2}} \\
+\int_{1}^{\infty }\int_{1}^{\infty }\mathcal{F}^{-1}\left( \left( \left(
\psi _{1}\otimes \psi _{2}\right) a_{t_{1},t_{2}}\right)
_{t_{1}^{-1},t_{2}^{-1}}\right) \left( x_{1},x_{2}\right) \frac{\func{d}t_{1}%
}{t_{1}}\frac{\func{d}t_{2}}{t_{2}}, \\
\left( x_{1},x_{2}\right) \in 
%TCIMACRO{\U{211d} }%
%BeginExpansion
\mathbb{R}
%EndExpansion
^{n_{1}}\mathbb{\times }%
%TCIMACRO{\U{211d} }%
%BeginExpansion
\mathbb{R}
%EndExpansion
^{n_{2}}\backslash \mathcal{M},
\end{multline*}%
it follows that 
\begin{multline*}
\left\vert \mathcal{F}^{-1}a\left( x_{1},x_{2}\right) \right\vert \leq
\left\vert \mathcal{F}^{-1}a_{0}\left( x_{1},x_{2}\right) \right\vert \\
+C_{N}\left\langle x_{1}\right\rangle ^{-N}\left\langle x_{2}\right\rangle
^{-N}\int_{1}^{\infty }t_{1}^{m_{1}+n_{1}}\left\langle
t_{1}x_{1}\right\rangle ^{-M}\frac{\func{d}t_{1}}{t_{1}} \\
+C_{N}\left\langle x_{1}\right\rangle ^{-N}\left\langle x_{2}\right\rangle
^{-N}\int_{1}^{\infty }t_{2}^{m_{2}+n_{2}}\left\langle
t_{2}x_{2}\right\rangle ^{-M}\frac{\func{d}t_{2}}{t_{2}} \\
+C_{N}\left\langle x_{1}\right\rangle ^{-N}\left\langle x_{2}\right\rangle
^{-N}\int_{1}^{\infty }\int_{1}^{\infty
}t_{1}^{m_{1}+n_{1}}t_{2}^{m_{2}+n_{2}}\left\langle t_{1}x_{1}\right\rangle
^{-M}\left\langle t_{2}x_{2}\right\rangle ^{-M}\frac{\func{d}t_{1}}{t_{1}}%
\frac{\func{d}t_{2}}{t_{2}}, \\
\left( x_{1},x_{2}\right) \in 
%TCIMACRO{\U{211d} }%
%BeginExpansion
\mathbb{R}
%EndExpansion
^{n_{1}}\mathbb{\times }%
%TCIMACRO{\U{211d} }%
%BeginExpansion
\mathbb{R}
%EndExpansion
^{n_{2}}\backslash \mathcal{M},
\end{multline*}%
or equivalently 
\begin{multline*}
\left\langle x_{1}\right\rangle ^{N}\left\langle x_{2}\right\rangle
^{N}\left\vert \mathcal{F}^{-1}a\left( x_{1},x_{2}\right) \right\vert \leq
\left\langle x_{1}\right\rangle ^{N}\left\langle x_{2}\right\rangle
^{N}\left\vert \mathcal{F}^{-1}a_{0}\left( x_{1},x_{2}\right) \right\vert \\
+C_{N}\left\vert x_{1}\right\vert ^{-m_{1}-n_{1}}\int_{\left\vert
x_{1}\right\vert }^{\infty }t_{1}^{m_{1}+n_{1}}\left\langle
t_{1}\right\rangle ^{-M}\frac{\func{d}t_{1}}{t_{1}} \\
+C_{N}\left\vert x_{2}\right\vert ^{-m_{2}-n_{2}}\int_{\left\vert
x_{2}\right\vert }^{\infty }t_{2}^{m_{2}+n_{2}}\left\langle
t_{2}\right\rangle ^{-M}\frac{\func{d}t_{2}}{t_{2}} \\
+C_{N}\left\vert x_{1}\right\vert ^{-m_{1}-n_{1}}\left\vert x_{2}\right\vert
^{-m_{2}-n_{2}}\int_{\left\vert x_{1}\right\vert }^{\infty }\int_{\left\vert
x_{2}\right\vert }^{\infty }t_{1}^{m_{1}+n_{1}}\left\langle
t_{1}\right\rangle ^{-M}t_{2}^{m_{2}+n_{2}}\left\langle t_{2}\right\rangle
^{-M}\frac{\func{d}t_{1}}{t_{1}}\frac{\func{d}t_{2}}{t_{2}}, \\
\left( x_{1},x_{2}\right) \in 
%TCIMACRO{\U{211d} }%
%BeginExpansion
\mathbb{R}
%EndExpansion
^{n_{1}}\mathbb{\times }%
%TCIMACRO{\U{211d} }%
%BeginExpansion
\mathbb{R}
%EndExpansion
^{n_{2}}\backslash \mathcal{M}.
\end{multline*}

If $m_{1}+n_{1}>0$ and $m_{2}+n_{2}>0$, then 
\begin{equation*}
C=\max \left\{ \int_{0}^{\infty }t_{1}^{m_{1}+n_{1}}\left\langle
t_{1}\right\rangle ^{-M}\frac{\func{d}t_{1}}{t_{1}},\int_{0}^{\infty
}t_{2}^{m_{2}+n_{2}}\left\langle t_{2}\right\rangle ^{-M}\frac{\func{d}t_{2}%
}{t_{2}}\right\} <\infty
\end{equation*}
and 
\begin{multline*}
\left\langle x_{1}\right\rangle ^{N}\left\langle x_{2}\right\rangle
^{N}\left\vert \mathcal{F}^{-1}a\left( x_{1},x_{2}\right) \right\vert \leq
\left\langle x_{1}\right\rangle ^{N}\left\langle x_{2}\right\rangle
^{N}\left\vert \mathcal{F}^{-1}a_{0}\left( x_{1},x_{2}\right) \right\vert \\
+CC_{N}\left( \left\vert x_{1}\right\vert ^{-m_{1}-n_{1}}+\left\vert
x_{2}\right\vert ^{-m_{2}-n_{2}}+C\left\vert x_{1}\right\vert
^{-m_{1}-n_{1}}\left\vert x_{2}\right\vert ^{-m_{2}-n_{2}}\right) , \\
\left( x_{1},x_{2}\right) \in 
%TCIMACRO{\U{211d} }%
%BeginExpansion
\mathbb{R}
%EndExpansion
^{n_{1}}\mathbb{\times }%
%TCIMACRO{\U{211d} }%
%BeginExpansion
\mathbb{R}
%EndExpansion
^{n_{2}}\backslash \mathcal{M}.
\end{multline*}

Assume that $-n_{1}<m_{1}<0$, $-n_{2}<m_{2}<0$ and $N\geq \max \left\{
n_{1}+1,n_{2}+1\right\} $. Then using part $(\func{b})$ of the previous
lemma we conclude that $\mathcal{F}^{-1}a\in L^{1}\left( 
%TCIMACRO{\U{211d} }%
%BeginExpansion
\mathbb{R}
%EndExpansion
^{n_{1}}\mathbb{\times }%
%TCIMACRO{\U{211d} }%
%BeginExpansion
\mathbb{R}
%EndExpansion
^{n_{2}}\right) $. Since $\mathbf{m\leq }\overline{\mathbf{m}}\Rightarrow 
\mathcal{S}^{\mathbf{m}}\left( 
%TCIMACRO{\U{211d} }%
%BeginExpansion
\mathbb{R}
%EndExpansion
^{n_{1}}\mathbb{\times }%
%TCIMACRO{\U{211d} }%
%BeginExpansion
\mathbb{R}
%EndExpansion
^{n_{2}}\right) \subset \mathcal{S}^{\overline{\mathbf{m}}}\left( 
%TCIMACRO{\U{211d} }%
%BeginExpansion
\mathbb{R}
%EndExpansion
^{n_{1}}\mathbb{\times }%
%TCIMACRO{\U{211d} }%
%BeginExpansion
\mathbb{R}
%EndExpansion
^{n_{2}}\right) $, it follows that $\mathcal{F}^{-1}a$ belongs to $%
L^{1}\left( 
%TCIMACRO{\U{211d} }%
%BeginExpansion
\mathbb{R}
%EndExpansion
^{n_{1}}\mathbb{\times }%
%TCIMACRO{\U{211d} }%
%BeginExpansion
\mathbb{R}
%EndExpansion
^{n_{2}}\right) $ for any $a\in \mathcal{S}^{\mathbf{m}}\left( 
%TCIMACRO{\U{211d} }%
%BeginExpansion
\mathbb{R}
%EndExpansion
^{n_{1}}\mathbb{\times }%
%TCIMACRO{\U{211d} }%
%BeginExpansion
\mathbb{R}
%EndExpansion
^{n_{2}}\right) $ with $m_{1}<0$ and $m_{2}<0$.

Thus we have proved the following

\begin{proposition}
Let $a\in \mathcal{S}^{\mathbf{m}}\left( 
%TCIMACRO{\U{211d} }%
%BeginExpansion
\mathbb{R}
%EndExpansion
^{n_{1}}\mathbb{\times }%
%TCIMACRO{\U{211d} }%
%BeginExpansion
\mathbb{R}
%EndExpansion
^{n_{2}}\right) =\mathcal{S}^{m_{1},m_{2}}\left( 
%TCIMACRO{\U{211d} }%
%BeginExpansion
\mathbb{R}
%EndExpansion
^{n_{1}}\mathbb{\times }%
%TCIMACRO{\U{211d} }%
%BeginExpansion
\mathbb{R}
%EndExpansion
^{n_{2}}\right) $. We denote by $\mathcal{M}$ the set $\left\{ \left(
x_{1},x_{2}\right) \in 
%TCIMACRO{\U{211d} }%
%BeginExpansion
\mathbb{R}
%EndExpansion
^{n_{1}}\mathbb{\times }%
%TCIMACRO{\U{211d} }%
%BeginExpansion
\mathbb{R}
%EndExpansion
^{n_{2}}:\left\vert x_{1}\right\vert \left\vert x_{2}\right\vert =0\right\} $%
. Then:

$(\func{i})$ $\mathcal{F}^{-1}a\in \mathcal{S}\left( 
%TCIMACRO{\U{211d} }%
%BeginExpansion
\mathbb{R}
%EndExpansion
^{n_{1}}\mathbb{\times }%
%TCIMACRO{\U{211d} }%
%BeginExpansion
\mathbb{R}
%EndExpansion
^{n_{2}}\backslash \mathcal{M}\right) $.

$(\func{ii})$ If $m_{1}+n_{1}>0$ and $m_{2}+n_{2}>0$, then for any $N\in 
%TCIMACRO{\U{2115} }%
%BeginExpansion
\mathbb{N}
%EndExpansion
$, there is $C_{N}>0$ such that 
\begin{multline*}
\left\vert \mathcal{F}^{-1}a\left( x_{1},x_{2}\right) \right\vert \leq
C_{N}\left\langle x_{1}\right\rangle ^{-N}\left\langle x_{2}\right\rangle
^{-N}\left( 1+\left\vert x_{1}\right\vert ^{-m_{1}-n_{1}}\right) \left(
1+\left\vert x_{2}\right\vert ^{-m_{2}-n_{2}}\right) , \\
\left( x_{1},x_{2}\right) \in 
%TCIMACRO{\U{211d} }%
%BeginExpansion
\mathbb{R}
%EndExpansion
^{n_{1}}\mathbb{\times }%
%TCIMACRO{\U{211d} }%
%BeginExpansion
\mathbb{R}
%EndExpansion
^{n_{2}}\backslash \mathcal{M}.
\end{multline*}

$\left( \func{iii}\right) $ If $m_{1}<0$ and $m_{2}<0$, then $\mathcal{F}%
^{-1}a\in L^{1}\left( 
%TCIMACRO{\U{211d} }%
%BeginExpansion
\mathbb{R}
%EndExpansion
^{n_{1}}\mathbb{\times }%
%TCIMACRO{\U{211d} }%
%BeginExpansion
\mathbb{R}
%EndExpansion
^{n_{2}}\right) $.
\end{proposition}

\begin{corollary}
Let $a\in \mathcal{S}^{\mathbf{m}}\left( 
%TCIMACRO{\U{211d} }%
%BeginExpansion
\mathbb{R}
%EndExpansion
^{n_{1}}\mathbb{\times }%
%TCIMACRO{\U{211d} }%
%BeginExpansion
\mathbb{R}
%EndExpansion
^{n_{2}}\right) =\mathcal{S}^{m_{1},m_{2}}\left( 
%TCIMACRO{\U{211d} }%
%BeginExpansion
\mathbb{R}
%EndExpansion
^{n_{1}}\mathbb{\times }%
%TCIMACRO{\U{211d} }%
%BeginExpansion
\mathbb{R}
%EndExpansion
^{n_{2}}\right) $ and $1\leq p\leq \infty $. If $m_{1}<0$ and $m_{2}<0$,
then $a\left( P_{%
%TCIMACRO{\U{211d} }%
%BeginExpansion
\mathbb{R}
%EndExpansion
^{n_{1}}},P_{%
%TCIMACRO{\U{211d} }%
%BeginExpansion
\mathbb{R}
%EndExpansion
^{n_{2}}}\right) \in \mathcal{B}\left( L^{p}\left( 
%TCIMACRO{\U{211d} }%
%BeginExpansion
\mathbb{R}
%EndExpansion
^{n_{1}}\times 
%TCIMACRO{\U{211d} }%
%BeginExpansion
\mathbb{R}
%EndExpansion
^{n_{2}}\right) \right) $.
\end{corollary}

\begin{corollary}
If $s_{1},s_{2}\geq 0,\varepsilon >0$ and $1\leq p\leq \infty $, then 
\begin{multline*}
\left( I-\triangle _{1}\otimes I\right) ^{s_{1}}\left( I-I\otimes \triangle
_{2}\right) ^{s_{2}}\left( I-\triangle _{1}\otimes I-I\otimes \triangle
_{2}\right) ^{-s_{1}-s_{2}-\varepsilon } \\
\in \mathcal{B}\left( L^{p}\left( 
%TCIMACRO{\U{211d} }%
%BeginExpansion
\mathbb{R}
%EndExpansion
^{n_{1}}\times 
%TCIMACRO{\U{211d} }%
%BeginExpansion
\mathbb{R}
%EndExpansion
^{n_{2}}\right) \right) .
\end{multline*}
\end{corollary}

\end{document}